\documentclass[11pt]{article}
\usepackage[dvipsnames]{xcolor}
\usepackage[utf8]{inputenc}
\usepackage{tikz}
\usetikzlibrary{shapes,arrows,positioning}
\usetikzlibrary{automata,arrows,positioning,calc}
\usepackage{amsmath}
\usepackage{amssymb}
\usepackage{amsthm}

\usepackage{dutchcal}
\usepackage{hyperref}

\usepackage{bbm}
\usepackage{blkarray}
\usepackage{algorithmic}
\usepackage{caption}
\usepackage{subcaption}
\usepackage{float}
\usepackage{stmaryrd}
\usepackage[normalem]{ulem}
\usepackage{enumitem}

\usepackage{fancyhdr}
\usepackage{amsopn}

\usepackage{booktabs}

\definecolor{burgundy}{rgb}{0.5,0.0, 0.13}

\hypersetup{colorlinks,
	citecolor=burgundy,
	citebordercolor=burgundy,
	linkcolor=burgundy,
	urlcolor=MidnightBlue
}

\DeclareMathOperator{\mbsP}{{\boldsymbol{P}}}

\usepackage{natbib}
\def\BIBand{and}%
\bibpunct[, ]{[}{]}{,}{n}{}{,}%

\providecommand{\keywords}[1]
{
	\small	
	\textbf{{Keywords.}} #1
}

\providecommand{\AMS}[1]
{
	\small	
	\textbf{{AMS subject classifications.}} #1
}

\usepackage[linesnumbered,ruled,vlined]{algorithm2e}

\numberwithin{equation}{section}
\allowdisplaybreaks

\usepackage{todonotes}

\newtheorem{theorem}{Theorem}[section]
\newtheorem{remark}[theorem]{Remark}
\newtheorem{assumption}[theorem]{Assumption}
\newtheorem{definition}[theorem]{Definition}
\newtheorem{lemma}[theorem]{Lemma}
\newtheorem{proposition}[theorem]{Proposition}

\DeclareMathOperator*{\argmin}{arg\,min}

\newcommand{\pzcN}{N}
\newcommand{\bsbeta}{{\boldsymbol{\beta}}}
\newcommand{\ubsbeta}{{\boldsymbol{\underline\beta}}}
\newcommand{\bsX}{\boldsymbol{X}}
\newcommand{\cC}{\mathcal{C}}
\newcommand{\EE}{\mathbb{E}}

\newcommand{\RR}{\mathbb{R}}
\newcommand{\TT}{\mathbb{T}}
\newcommand{\UU}{\mathbb{U}}
\newcommand{\VV}{\mathbb{V}}

\newcommand{\mcE}{\mathcal{E}}
\newcommand\ff{\mathfrak f}

\newcommand\bN{\mathbf N}

\newcommand\cP{\mathcal{P}}

\title{A Machine Learning Method for Stackelberg Mean Field Games}

\author{G\"ok\c ce Dayan{\i}kl{\i}\footnote{Department of Statistics, 
	University of Illinois at Urbana-Champaign, Champaign, IL 61820, USA
		(\href{mailto:gokced@illinois.edu}{gokced@illinois.edu}).}
	\and Mathieu Lauri\`ere
	\footnote{Shanghai Frontiers Science Center of Artificial Intelligence and Deep Learning; NYU-ECNU Institute of Mathematical Sciences, NYU Shanghai, Shanghai, 200126, People’s Republic of China 
		(\href{mailto:mathieu.lauriere@nyu.edu}{mathieu.lauriere@nyu.edu}).}
}
\date{}

\begin{document}

\maketitle

\begin{abstract}
We propose a single-level numerical approach to solve Stackelberg mean field game (MFG) problems. In Stackelberg MFG, an infinite population of agents play a non-cooperative game and choose their controls to optimize their individual objectives while interacting with the principal and other agents through the population distribution. The principal can influence the mean field Nash equilibrium at the population level through policies, and she optimizes her own objective, which depends on the population distribution. This leads to a bi-level problem between the principal and mean field of agents that cannot be solved using traditional methods for MFGs. We propose a reformulation of this problem as a single-level mean field optimal control problem through a penalization approach. We prove convergence of the reformulated problem to the original problem. We propose a machine learning method based on (feed-forward and recurrent) neural networks and illustrate it on several examples from the literature.
\end{abstract}

\vskip3mm

\vskip3mm
\AMS{
	49N90, 
	91A13, 
	91A15, 
	62M45. 
}

\section{Introduction.}
\label{sec:introduction}

In policy making, finding \textit{optimal} policies to solve socioeconomic problems is the ultimate goal. However, this problem has an underlying complexity: Individualistic nature of humankind prevents policy makers to directly control the behavior of people. Instead, the policymakers should take into account the reaction of the society -- consisting of non-cooperative agents -- to a policy while deciding on the best one. One approach to understand the emergent reactions to any given policy can be to use simulation techniques such as agent-based simulation, which is commonly used for modeling complex interactions among agents. However, it lacks the tractability of the solutions. This lack of tractability issue prevents us from solving for \textit{optimal} policies. Instead, with this approach only the outcomes of different policies can be compared through simulations. Therefore, in order to attain tractability of the solutions when there are many agents interacting with each other, we can utilize game theoretical tools such as \textit{mean field games} (MFGs). Intuitively, in the MFG setup, we focus on a game among a large number of indistinguishable agents (i.e., players) that are interacting symmetrically. Then, we study the equilibrium between a representative player and the distribution of the other players’ states (and possibly actions) instead of focusing on the interactions of every player with each other. 
With this simplification, we can characterize an approximate Nash equilibrium in the society given an policy (i.e., incentive or contract) by using forward-backward differential equations. Even if this forward-backward system could be hard to analyze, mean field equilibria are simpler to identify and compute than equilibria of populations with finite but large number of agents because of the curse of dimensionality that results from the exponentially increasing number of interactions between the agents when the number of agents increases. Therefore, the mean field equilibria provide approximate Nash equilibria for games with a large but finite number of players.

After finding a tractable solution for the non-cooperative game (i.e., Nash equilibrium) of large number of agents given any policy by using MFG approach, our second step is to find the \textit{optimal} policy. This can be done by the addition of a \textit{principal} to the model who has her own objectives (different than the objectives of the agents in the society) and chooses \textit{policies} in order to affect the Nash equilibrium in the society. The addition of this principal changes our problem from a standard MFG to a Stackelberg MFG between a principal and a mean field population of agents. This Stackelberg MFG problem is an intrinsically bi-level problem where we look for the mean field equilibrium in the society given any policy at the population level and optimize the policy of the principal by taking into account this equilibrium reaction from the society at the top level. Mathematically, this corresponds to solving an optimal control problem that is constrained by an MFG equilibrium.

In this paper, we focus on writing the Stackelberg MFG problem as a single level optimization problem which is a McKean-Vlasov (MKV) control problem with a constraint at the terminal time, and proposing a machine learning approach to solve Stackelberg MFG.

\subsection{Literature review.}
\label{subsec:intro_literature}
\subsubsection{Mean field games.}

When there is a large number of agents in a game setting, computing a Nash equilibrium is a daunting task because of the exponentially increasing number of interactions among the agents. The MFG approach was proposed in \cite{huang2006large} and \cite{lasry2006jeux,lasry2006jeux2} simultaneously and independently to overcome theses challenges. In MFGs, it is assumed that there is infinitely many indistinguishable agents interacting in a symmetric way. The agents do not directly see the state or the action of another specific agent. Instead, each agent sees the population distribution of states (or controls or joint distribution of states and controls). In the mean-field regime, Nash equilibria can be characterized by a system of forward-backward equations. In the analytical approach, the system consists of two partial differential equations (PDEs): a Hamilton-Jacobi-Bellman and a Kolmogorov-Fokker-Plank equations. In the probabilistic approach, the system consists of two stochastic differential equations, and the system is referred to as a forward-backward stochastic differential equation (FBSDE) system. Interested readers can refer to the books \cite{bensoussan2013mean,carmona2018probabilistic}. %

\subsubsection{Contract theory and Stackelberg mean field games.}

Contract theory studies the interaction between a principal and an agent. The principal proposes a contract to the agent and the agent reacts to the given contract so as to minimize his cost (or, equivalently, to maximize her reward). For example, the agent can decide how much effort he will put in the work while working for the principal. For a given contract, the agent needs to solve an optimal control problem. The goal of the principal is to find the best contract to minimize her own cost (or to maximize her own reward). Solutions to contract theory problems are generally studied using the notion of Stackelberg equilibrium. In~\cite{holmstrom1987aggregation} continuous time method is used to study this type of problems, and in~\cite{Sannikov2008,Sannikov2013} dynamic programming and martingale optimality principles are used to characterize the solution. These ideas are later generalized in~\cite{Cvitanic2018}. In~\cite{djehiche2014principal} the solution in a general class of principal agent problem is characterized by using the Pontryagin stochastic maximum principle. In~\cite{bensoussan_stackelberg_max}, the authors discuss different information structures in Stackelberg games and derive a maximum principle for a specific information structure. In recent literature, existence of optimal contracts is shown by using compactness arguments in~\cite{alvarez2023optimal} under linearly controlled drift assumptions by using a relaxed formulation of the controls in the principal-agent problem and in~\cite{krsek2023randomisation} by allowing mixed (i.e., measure-valued) controls.

In the model we consider in this paper, instead of having one principal and one agent, we have a population of large number of agents and we will model the large number of agents with a mean field game. In our setup, the principal tries to optimize her policies by taking the mean field Nash equilibrium reaction of the population to the given policy. This creates a Stackelberg equilibrium between the principal and the mean field of agents. In the context of MFGs, problems with a principal and a mean field of agents have been studied in~\cite{Elie_2019} in the continuous state space setting and in~\cite{wang_2020} for finite state spaces. Linear-quadratic models have received particular attention for their tractability, see e.g.~\cite{nourian2012mean,MR3376121,bensoussan2017linearstackelberg,moon2018linearstackelberg,lin2018openstackelberg,oner2018mean,yang2021linear,wang2022leader}. %
Stackelberg MFG approach is used to find equilibrium in many applications such as optimal advertisement decisions \cite{salhab2016,advertisement,SALHAB20221079}, %
epidemic control \cite{StackelbergMFG_epidemics}, optimal energy demand response with accidents \cite{elie2021mean,mastrolia2022agency}, or optimal carbon tax regulations \cite{carbon_StackelbergMFG}. The theory has been extended to cover problems with delay in~\cite{MR3376121}, with a terminal state constraint in~\cite{fu2018meanfield} or problems with a collaborative population~\cite{hubert2022incentives}, which corresponds to Stackelberg equilibrium between a principal and another player with McKean-Vlasov dynamics. A discrete time setting has been analyzed in~\cite{guo2022optimizationstackelberg} using an optimization viewpoint. 
\cite{hu2022principal} considered a model with a single agent and a mean field of principals. The connection with finite-player games has been studied rigorously for instance in~\cite{bergault2023mean} and ~\cite{djete2023stackelberg} using respectively analytical and probabilistic techniques. An extension to hierarchical Stackelberg MFGs is studied in~\cite{oner2018mean}. 
Last, Stackelberg MFG should be distinguished from major-minor MFGs, which usually refer to a situation in which one player interacts with a mean field of infinitesimal agents and the whole situation is described as a Nash equilibrium among all the players. Interested readers can refer to~\cite[Chapter 7.1]{carmona2018probabilistic2}. In~\cite{advertisement} the differences between Stackelberg mean field games and major-minor mean field games are discussed in an advertisement competition application with multiple major players.

\subsubsection{Numerical approaches in Stackelberg MFGs.}

Even though MFGs provide a simplified way to find approximate equilibria in games with finite but large number of agents, the forward-backward differential equation systems characterizing the mean field Nash equilibrium are still difficult to solve. Explicit solutions are available only in very few cases, such as in linear-quadratic MFGs. However, these models are not sufficient to capture the complexity of real life problems. On top of this difficulty, adding the optimization of the principal in the model introduces another level of optimization to our problem and in general the solution cannot be characterized using a classical PDE or FBSDE system. As a consequence, we cannot rely directly on traditional numerical methods for MFGs based on solving PDEs or FBSDEs~\cite{achdou2010mean,achdou2012mean,briceno2018proximal,carlini2014fully,chassagneux2019numerical,angiuli2019cemracs} (see also~\cite{achdou2020mean} for a recent survey). Therefore, new numerical approaches needs to be developed to solve Stackelberg MFG problems.

Recently, a few works have considered numerical approaches for specific classes of Stackelberg MFGs. In~\cite{StackelbergMFG_epidemics}, the authors are motivated by finding the optimal time dependent social distancing policies and path dependent stimulus payment to mitigate an epidemic and model the government as the principal and the people in the society as the mean field game population. In this paper, building on the method proposed in~\cite{carmona2022convergence}, a numerical approach that uses deep learning tools and Monte Carlo simulation is proposed to solve finite state Stackelberg MFG as a single-level optimization problem. However, this approach required a specific relationship between the terminal costs of principal and the agents. In~\cite{Jaimungal}, authors propose a deep learning method to solve Stackelberg MFG problem where the principal chooses only the terminal payment structure among piecewise-linear function set. With this simplification, they solve bi-level Stackelberg MFG where they fix a control for the principal and solve the equilibrium and then optimize over the control set of the principal. In this work, we propose a method that can tackle a much wider class of Stackelberg MFGs, in terms of cost functions and dynamics.

\subsection{Contributions and paper structure.}
\label{subsec:intro_contribution}

The main contributions of this paper are three-fold. As the first contribution, we propose a general approach which is based on rephrasing the Stackelberg MFG as a single-level stochastic optimal control problem driven by forward SDEs instead of the original formulation involving bi-level optimization. Roughly speaking, the proposed approach relies on writing the Stackelberg MFG problem as an MKV stochastic optimal control problem with a constraint at the terminal time, adding a penalty to merge the constraint with the principal's cost, and finally approximating of the controls by using parameterized functions. 
As the second contribution, we provide a theoretical foundation to our approach by proving convergence results in terms of the penalization, empirical distribution approximation, time dicretization, and the function approximation. As the third contribution, we provide numerical experiments on several examples from the literature. In particular, even if we give the continuous space setup in our paper, our approach is suitable for both continuous and finite space models with general cost functions.

In Section~\ref{sec:model}, we introduce the model setting and the analysis of the model. Furthermore, we rewrite the Stackelberg MFG model as a single-level stochastic optimal control problem. In Section~\ref{sec:main_theory}, we prove the convergence of the solution of rewritten optimization problem to the original problem and give results on the convergence of the numerical approach that is introduced in Section~\ref{sec:numerics}. In Section~\ref{sec:numerics}, we give details of our numerical approach and discuss some extensions in Section~\ref{subsec:numerics_extended}. 
In Section~\ref{sec:experiments}, we evaluate the proposed numerical approaches with different models with applications to contract theory and systemic risk problems. Finally, in Section~\ref{sec:conclusion}, we conclude by summarizing our contributions and discuss some future work.

\section{Model.}
\label{sec:model}
In this section, we introduce the mathematical model for the Stackelberg MFG. We analyze the solution of the mean field game given the policies of the principal. In the last subsection, we rewrite the bi-level Stackelberg MFG problem as a single-level problem.

\subsection{Preliminaries and notation.}
\label{subsec:model_prelim}

We adopt the following notation throughout the paper. Let $T>0$ be a finite time horizon. We assume that the problem is set on a complete filtered probability space $(\Omega, \mathcal{F}, \mathbb{F} = (\mathcal F_t)_t, \mathbb{P})$ supporting a $d$-dimensional Wiener process $\boldsymbol W = (W_t)_{t\in [0,T]}$. The representative infinitesimal agent has an $m$-dimensional continuous state process $\bsX = (X_t)_{t\in [0,T]}$ where $ X_t \in \mathbb{R}^{{m}}$ for all $t\in[0,T]$. We endow $\mathcal{P}(\mathbb{R}^m)$, the space of probability measures on $\mathbb{R}^m$, with the 2-Wasserstein distance denoted by $W_2$ and defined as:
\begin{equation}
    \label{eq:wasserstein_definition}
    W_2(\mu,\nu):= \inf_{\pi \in \Pi(\mu, \nu)}\left[\int_{\mathbb{R}^m \times \mathbb{R}^m} d(x,y)^2d\pi(dx,dy)\right]^{1/2},
\end{equation}
where $\Pi(\mu, \nu)$ denotes the set of probability measures in $\mathcal{P}(\mathbb{R}^m \times \mathbb{R}^m)$ with marginals $\mu$ and $\nu$ and $d(x,y)$ denotes the Euclidean distance between vectors $x$ and $y$. Let $M(\mathcal{P}(\mathbb{R}^{{m}}))$ be the set of measurable mappings from $[0,T]$ to $\mathcal{P}(\mathbb{R}^{{m}})$. We denote the collection of square-integrable and $\mathbb{F}$-adapted action processes $\boldsymbol \alpha = (\alpha_t)_{t \in [0,T]}$ with $\mathbb{A}$. The action space is denoted by $A$ is a subset of $\mathbb{R}^k$ i.e., $\alpha_t \in A\subset \mathbb{R}^k$.

Furthermore, given the complete filtered probability space $(\Omega, \mathcal{F}, \mathbb{F}, \mathbb{P})$, we denote by $\mathcal{H}^{2,n}$ the set of $\mathbb{F}$-adapted and real-valued processes 

$$\mathcal{H}^{2,n}=\left\{\boldsymbol Z \in\mathcal{H}^{0,n}: \mathbb{E}\int_0^T|Z_s|^2ds<\infty\right\},$$ 
where $\mathcal{H}^{0,n}$ is the collection of all $\mathbb{R}^n$-valued measurable processes on [0,T]. 

\subsection{General setting for Stackelberg mean field game.}
\label{subsec:model_SMFGsetup}

We consider {\it an infinite population of non-cooperative players} and {\it one principal}. Since the population of non-cooperative players is infinite, each player has a negligible influence on the rest of the population. 
The infinitesimal agents are influenced by the actions of the principal, which can be thought of as a regulator who can influence the Nash equilibrium of the mean field population. We begin by focusing on the game between the members of the population. As is common in the MFG paradigm, we focus on a representative player. Given knowledge of how the population and the principal act over time, the representative player optimizes their cost functional.

Let $\Lambda$ be the set of possible actions for the principal. Although more general sets could be considered, to simplify the presentation, we will use $\Lambda = \RR$ and we consider policies that are elements of $C^0([0,T]{\times \RR^m}; \RR)$, which is the set of real-valued continuous functions of time and the state of the individual agent.\footnote{Under further regularity assumptions, we can actually consider policies that are functions of time, the state of the individual agent, and the state distribution of the population (i.e., $\lambda(t, x, \mu)$) or functions of time and the state distribution of the population (i.e., $\lambda(t, \mu)$) depending on the application interest. %
} Assume the principal has chosen policy $\boldsymbol \lambda = (\lambda(t{, x}))_{t\in[0,T]{, x \in \RR^m}} \in C^0([0,T]{\times \RR^m}; \Lambda)$, and the flow of population distribution is given by $\boldsymbol{\mu} = (\mu_t)_{t \in [0,T]}\in M(\mathcal{P}(\mathbb{R}^{{m}}))$. Consider an infinitesimal agent using control $\boldsymbol\alpha = (\alpha_t)_{t\in[0,T]} \in \mathbb{A}$. 
Then, the state of this agent is described by an $m$-dimensional stochastic process $\bsX$ which has the following dynamics:
\begin{equation}
\label{eq:dynamics_of_minors}
\begin{aligned}
    dX_t =  b(t, X_t, \alpha_t, \mu_t; \lambda(t{, X_t})) dt + \sigma dW_t, \qquad X_0 =\zeta,
\end{aligned}
\end{equation}
where $b : [0,T] \times \mathbb{R}^{m}\times A\times \mathcal{P}(\mathbb{R}^{{m}}) \times \Lambda \to  \mathbb{R}^{{m}}$ is the drift of the agent's state process that depends on the agent's own state $X_t$, her own action $\alpha_t$, the mean field distribution $\mu_t$ and the principal's policy $\lambda(t{, X_t})$. $\sigma >0$ is a constant volatility and $\boldsymbol W = (W_t)_{t\in[0,T]}$ is a $d$-dimensional Wiener process that represents the idiosyncratic noise. The initial condition of the state process is denoted by $\zeta \in L^2(\Omega,\mathcal F_0, \mathbb{P}; \mathbb{R}^{{m}})$ and has distribution $\mu_0 = \mathcal{L}(\zeta)$, which is also the initial state distribution of the population.

The representative player pays the expected total cost:
\begin{equation}
\label{eq:cost_of_minors}
    J^{\boldsymbol\lambda}(\boldsymbol\alpha, \boldsymbol\mu)
    := 
    \mathbb{E}\left[
    \int_0^T f(t, X_{t}, \alpha_t, \mu_t; \lambda(t{,X_t}))dt + g(X_T, \mu_T; \lambda(T{,X_T}))
    \right],
\end{equation}
where $f : [0,T]\times \mathbb{R}^{m}\times A\times \mathcal{P}(\mathbb{R}^{{m}}) \times \Lambda \to \mathbb{R}$ is a running cost which depends on the agent's own state $X_t$, her own action $\alpha_t$, the mean field distribution $\mu_t$, and the principal's policy $\lambda(t{,X_t})$. Finally, $g : \mathbb{R}^m \times \mathcal{P}(\mathbb{R}^{{m}}) \times \Lambda \to \mathbb{R}$ is the terminal cost of the agent and it depends on the agent's own state, the mean field distribution, and the principal's policy at the terminal time.

We stress that the representative player interacts with the population through the mean field distributions $\boldsymbol{\mu} = (\mu_t)_{t\in[0,T]}$ in the running cost, terminal cost and the drift of the state dynamics and interacts with the principal through the policy of the principal $\boldsymbol{\lambda} = (\lambda(t{,x}))_{t\in[0,T]{,x\in \RR^m}}$ in the running cost, terminal cost and the drift of the state dynamics. Now, we can define the equilibrium in the mean field population given the principal's policy selection.
\begin{definition}
\label{def:mfg_nash}
Let $\boldsymbol\lambda \in C^0([0,T]{\times \RR^m}; \RR)$ be a policy of the principal. We say that $(\boldsymbol{\hat \alpha, \hat \mu}) \in \mathbb{A}$ is a mean field Nash equilibrium given the policy $\boldsymbol\lambda$ if the following two conditions are met:
\begin{itemize}
    \item[(i)] $\boldsymbol{\hat\alpha} = \arg\inf_{\boldsymbol{\alpha}\in\mathbb{A}} J^{\boldsymbol\lambda}(\boldsymbol\alpha, \boldsymbol{\hat{\mu}})$;
    \item[(ii)] For all $t\in[0,T]$, $\hat{\mu}_t = \mathcal{L}(X_t)$, where $\bsX$ has dynamics~\eqref{eq:dynamics_of_minors} with control $\boldsymbol{\hat\alpha}$.
\end{itemize}
We denote by $ \mathcal{N}(\boldsymbol\lambda)$ the set of mean field Nash equilibria.
\end{definition}

Intuitively, the first condition in Definition~\ref{def:mfg_nash} states that the representative player finds her best response given the population state distribution and the second condition says that since the players are identical, the distribution of the representative player's state induced by her best response control should correspond to the population state distribution fixed at the beginning.

The principal's problem is to find the policies that create the most favorable reaction of the agents to minimize the principal's own cost. By using policy $\boldsymbol \lambda$ as an incentive, the principal can affect the set of mean-field Nash equilibria $\mathcal{N}(\boldsymbol \lambda)$ and in this way she can influence the population's behavior to minimize her own cost.

\begin{definition}
\label{def:principal_policy_admissibility}
    A policy of the principal $\boldsymbol{\lambda}$ is admissible if, given the deterministic mapping $\boldsymbol\lambda $, $\mathcal{N}(\boldsymbol\lambda)$ is a singleton. We denote by $\boldsymbol\Lambda \subseteq C^0([0,T]\times \RR^m; \RR)$ the set of admissible policies. 
\end{definition}
\begin{remark}
\label{rem:MFG_uniqueness}
    One of the ways to ensure uniqueness is to assume Lasry-Lions monotonicity~\cite{lasry2006jeux,lasry2006jeux2}. 
    Interested readers can refer to \cite[Theorem 3.32]{carmona2018probabilistic} for assumptions that yield uniqueness of the solution in MFG problems.
	In our case, uniqueness of the Nash equilibrium is guaranteed provided the functions $\tilde{b}: (t,x,\alpha,\mu) \mapsto b(t,x,\alpha,\mu;\lambda(t,x))$ and $\tilde{f}: (t,x,\alpha,\mu) \mapsto f(t,x,\alpha,\mu;\lambda(t,x))$ satisfy the ``\textbf{L-monotonicity}'' condition on p. 175 of ~\cite{carmona2018probabilistic}.  A sufficient assumption to satisfy this condition is:
    \begin{itemize}
        \item $\lambda$ is of the form: $\lambda(t,x) = \lambda_0(t) + \lambda_1(t)x$ for some bounded measurable functions $\lambda_0, \lambda_1$.
        \item $b$ is of the form: $b(t,x,\alpha,\mu;\lambda) = b_0(t) + b_1(t)x + b_2(t)\alpha + b_3(t)\lambda$ for some bounded measurable functions $b_i, i=1,\dots,3$.
        \item $f$ has a separated structure: $f(t,x,\alpha,\mu;\lambda) = f_0(t,x,\mu;\lambda) + f_1(t,x,\alpha;\lambda)$ with $f_0$ and $f_1$ measurable.
        \item  Letting $\tilde{f}_0(t,x,\mu) = f_0(t,x,\mu;\lambda(t,x))$, and $\tilde{g}(x,\mu) = g(x,\mu;\lambda(T,x))$, the functions  $\tilde{f}_0(t,\cdot,\mu)$ and $\tilde{g}(\cdot,\mu)$ are $\cC^1$, with $\partial_x \tilde{f}_0(t,\cdot,\mu)$ and $\partial_x \tilde{g}(\cdot,\mu)$ being at most of linear growth, uniformly in $(t,\mu)$; letting   $\tilde{f}_1(t,x,\alpha) = f_1(t,x,\alpha;\lambda(t,x))$, the function $\tilde{f}_1(t,\cdot,\cdot)$ is $\cC^1$, with $\partial_{(x,\alpha)} \tilde{f}_1(t,\cdot,\cdot)$  being at most of linear growth in $(x,\alpha)$, uniformly in $t$; the function $t \mapsto \tilde{f}_1(t,0,0)$ is bounded, and $t \mapsto \tilde{f}_0(t,0,\mu)$ is bounded on bounded subsets.
        \item The functions $\partial_x \tilde{f}_0(t,\cdot,\cdot)$ and $\partial_x \tilde{g}$ are L-monotone in the sense of~\cite[Definition 3.31]{carmona2018probabilistic}. 
        Moreover, the function $\tilde{f}_1$ satisfies the strong form of convexity: $\tilde{f}_1(t,x',\alpha') - \tilde{f}_1(t,x,\alpha) - (x-x', \alpha-\alpha') \cdot \partial_{(x,\alpha)} \tilde{f}_1(t,x,\alpha) \ge C|\alpha' - \alpha|^2$ for some $C>0$. 
\end{itemize} 
    When $\lambda$ depends only on $t$ (i.e., $\lambda_1=0$ above), these conditions are satisfied for broad class of models.
\end{remark}

In order to use an admissible policy $\boldsymbol \lambda \in \boldsymbol\Lambda$ the principal pays the cost:
\begin{equation}
\label{eq:cost_of_principal}
    J^0(\boldsymbol \lambda) := {\EE\Big[}\int_0^T f_0(t, \hat \mu^{\boldsymbol \lambda}_t, \lambda(t{,X_t}))dt +g_0(\hat \mu_T^{\boldsymbol \lambda}, \lambda(T{,X_T})){\Big]}
\end{equation}
where $(\boldsymbol{\hat \alpha}^{\boldsymbol\lambda},\boldsymbol{\hat \mu}^{\boldsymbol\lambda})=\mathcal{N}(\boldsymbol \lambda)$. 
In this way, the principal's optimization problem is:
\begin{equation}
\label{eq:principal-threshold-kappa}
    \inf_{\boldsymbol\lambda\in \boldsymbol\Lambda}
    J^0(\boldsymbol\lambda).
\end{equation}

\begin{remark}
    The principal's objective is not just merely minimizing her own cost. The principal also needs to satisfy the constraint that the population is reacting the policy with a Nash equilibrium behavior. Because of this, our problem becomes a dynamic optimal control problem that is constrained by a Nash equilibrium condition. As we will see below, this constraint can be expressed in terms of a forward-backward system of stochastic differential equations (SDEs). The full problem is given in Subsection~\ref{subsec:model_SMFGrewrite}.
\end{remark}

\subsection{Analysis of the mean field game.}
\label{subsec:model_mfganalysis}
In this section, before discussing the full problem between the principal and the mean field population, we first analyze the mean field game in the population. We characterize the Nash equilibrium in the mean field population given a policy of the principal, $\boldsymbol \lambda$, and give an existence result for the mean field Nash equilibrium.

Given a policy $\boldsymbol \lambda\in \boldsymbol\Lambda$ and the mean field interactions $\boldsymbol{\mu}$, the Hamiltonian for the representative player's optimal control problem is given by the function $H : [0,T]\times \mathbb{R}^{{m}}\times\mathcal{P}(\mathbb{R}^{{m}})\times\mathbb{R}^{{m}}\times A \times \Lambda \mapsto \mathbb{R},$
\begin{equation}
    \label{eq:Hamiltonian}
    H(t,x,\mu,\gamma,\alpha;\lambda) := b(t,x, \alpha,\mu;\lambda) \cdot \gamma + f(t,x,\alpha,\mu;\lambda). 
\end{equation}
The minimizer of the Hamiltonian is denoted by $\hat \alpha_t(x,\mu,\gamma;\lambda)$. In other words, we have:
\begin{equation}
    \label{eq:Hamiltonian_minimizer}
    \hat \alpha_t(x,\mu,\gamma;\lambda) = \argmin_{\alpha} H(t, x,\mu, \gamma, \alpha;\lambda).
\end{equation}

In order to characterize the Nash equilibrium in the MFG, we follow the analysis in~\cite[Chapter 4.4]{carmona2018probabilistic}.

\begin{assumption}%
\label{assu:mfg_existence_value_rep}
We assume:
\begin{enumerate}[label=\textbf{(A\theassumption.\arabic*)},align=left]
    \item $b$ and $f$ are measurable mappings from $[0,T]\times\mathbb{R}^m\times A\times \mathcal{P}(\mathbb{R}^m)\times\Lambda$ into $\mathbb{R}^m$ and $\mathbb{R}$, respectively. $g$ is a measurable function from $\mathbb{R}^m \times \mathcal{P}(\mathbb{R}^m)\times \Lambda$ into $\mathbb{R}$. 
    \item \label{assu:mfg_existence_value_rep_bound-alpha} There exists a constant $L\geq0$ such that, 
    \begin{equation}
        \begin{aligned}
            |b(t,x,\alpha, \mu; \lambda)|&\leq L (1+|\alpha|),\\
            |f(t,x,\alpha, \mu; \lambda)|&\leq L (1+|\alpha|^2),\\
            |g(x,\mu;\lambda)|&\leq L,
        \end{aligned}
    \end{equation}
    for all $(t, x, \alpha, \mu, \lambda) \in [0,T]\times\mathbb{R}^m\times A\times \mathcal{P}(\mathbb{R}^m)\times\Lambda$
    \item \label{assu:mfg_existence_value_rep_cont-diff-alpha} For any $t \in[0,T]$, $x\in \mathbb{R}^m$, $\mu \in \mathcal{P}(\mathbb{R}^m)$ and $\lambda \in \Lambda$, the functions $b(t,x,\cdot,\mu; \lambda)$ and $f(t,x,\cdot,\mu; \lambda)$ are continuously differentiable in $\alpha$.
    \item \label{assu:mfg_existence_value_rep_lip-x} For any $t \in[0,T]$, $\alpha\in A$, $\mu \in \mathcal{P}(\mathbb{R}^m)$ and $\lambda \in \Lambda$, the functions $b(t,\cdot,\alpha,\mu; \lambda)$, $f(t,\cdot,\alpha,\mu; \lambda)$, $g(\cdot,\mu; \lambda)$, and $\lambda(t, \cdot)$ are L-Lipschitz continuous (in $x$). 
    \item \label{assu:mfg_existence_value_rep_cont-mu} For any $t \in[0,T]$, $x\in \mathbb{R}^m$, $\alpha\in A$, and $\lambda \in \Lambda$, the functions $b(t,x,\alpha,\cdot; \lambda)$, $f(t,x,\alpha,\cdot; \lambda)$, and $g(x,\cdot; \lambda)$ are continuous in the measure argument with respect to the 2-Wasserstein distance.
    \item \label{assu:mfg_existence_value_rep_bdd-diff-alpha} For the same constant $L$ and for all $(t, x, \alpha, \mu, \lambda) \in [0,T]\times\mathbb{R}^m\times A\times \mathcal{P}(\mathbb{R}^m)\times\Lambda$,
    \begin{equation}
        \begin{aligned}
            |\partial_\alpha b(t,x,\alpha, \mu; \lambda)| &\leq L,\\
            |\partial_\alpha f(t,x,\alpha, \mu; \lambda)| &\leq L.
        \end{aligned}
    \end{equation}
    \item \label{assu:mfg_existence_value_rep_cont-hat-alpha} There exists a unique minimizer $\hat \alpha_t(x,\mu,\gamma;\lambda)$ of the Hamiltonian $H(t,x,\mu,\gamma,\alpha; \lambda)$ given in \eqref{eq:Hamiltonian}. It is L-Lipschitz continuous in $(x,\gamma)$ and continuous in $\mu$ and satisfies: for all $(t, x, \alpha, \mu, \lambda) \in [0,T]\times\mathbb{R}^m\times A\times \mathcal{P}(\mathbb{R}^m)\times\Lambda$, 
    \begin{equation}
        |\hat \alpha_t(x,\mu,\gamma;\lambda)| \leq L(1+|\gamma|).
    \end{equation}
    \item\label{assu:mfg_existence_value_rep_cont-lambda} For any $t \in[0,T]$, $x\in \mathbb{R}^m$, $\alpha\in A$, $\mu \in \mathcal{P}(\RR^m)$ and $\gamma \in \RR^m$, the functions $b(t,x,\alpha,\mu; \cdot)$, $f(t,x,\alpha,\mu; \cdot)$, $g(x,\mu; \cdot)$, and $\hat\alpha_t(x,\mu,\gamma; \cdot)$ are $L$-Lipschitz continuous (in $\lambda$).
\end{enumerate}
\end{assumption}

\begin{theorem}
\label{the:MFG_existence}
    If Assumption~\ref{assu:mfg_existence_value_rep} holds, then for any initial condition $\zeta \in L^2(\Omega,\mathcal F_0, \mathbb{P}; \mathbb{R}^{{m}})$, the FBSDE system
\begin{equation}
    \begin{aligned}
        \label{eq:Nash_FBSDE_sytem}
        X_t &= \zeta +\int_0^t b(s, X_s, \hat \alpha_s(X_s, \mu_s, {\sigma^{-1}}^{\top} Z_s;\lambda(s{, X_s})), \mu_s; \lambda(s{, X_s})) ds + \int_0^t \sigma dW_s\\
        Y_t &= g(X_T,  \mu_T;\lambda(T{, X_T}))+ \int_t^T f(s, X_s, \hat \alpha_s (X_s, \mu_s, {\sigma^{-1}}^{\top} Z_s;\lambda(s{, X_s})), \mu_s; \lambda(s{, X_s}))ds -\int_t^T Z_s \cdot dW_s
    \end{aligned}
\end{equation}
with $\mu_s = \mathcal{L}(X_s),\ \forall s\in[0,T]$ is solvable and the flow given by the marginal distributions forward component of any solution is an equilibrium for the associated mean field game problem.
\end{theorem}

\proof[Proof of Theorem~\ref{the:MFG_existence}. ] Introduce 
\begin{equation*}
    \begin{aligned}
        B(t,x,\alpha,\mu) &:= b(t,x,\alpha,\mu;\lambda(t{, x})),\\
        F(t,x,\alpha,\mu) &:= F(t,x,\alpha,\mu;\lambda(t{, x})),\\
        \hat a_t(x,\mu, \gamma) &:= \hat \alpha_t(x,\mu, \gamma;\lambda(t{, x})),\\
        G(x,\mu) &:= g(x, \mu; \lambda(T{, x})). 
    \end{aligned}
\end{equation*}
Then the forward backward system in \eqref{eq:Nash_FBSDE_sytem} can be written as
\begin{equation}
    \begin{aligned}
        \label{eq:Nash_FBSDE_sytem_woutLambda}
        X_t &= \zeta +\int_0^t B(s, X_s, \hat a_s(X_s, \mu_s, {\sigma^{-1}}^{\top} Z_s), \mu_s) ds + \int_0^t \sigma dW_s\\
        Y_t &= G(x,\mu)+ \int_t^T F(s, X_s, \hat a_s (X_s, \mu_s, {\sigma^{-1}}^{\top} Z_s), \mu_s)ds -\int_t^T Z_s \cdot dW_s
    \end{aligned}
\end{equation}
The proof follows the same ideas given in \cite[Chapter 4.4]{carmona2018probabilistic}.
\endproof

\begin{remark}
We stress that the process $\boldsymbol{Y}$ in the FBSDE system~\eqref{eq:Nash_FBSDE_sytem} represents the value function of a representative agent.
\end{remark}

Theorem~\ref{the:MFG_existence} implies that the 
mean field Nash equilibrium given the principal's incentive can be found by solving the FBSDE. Furthermore, it deduces the existence of the mean field Nash equilibrium given the principal's incentive through the existence of the solution of the FBSDE. The uniqueness result can be inferred following the Remark~\ref{rem:MFG_uniqueness}.

\subsection{Rewriting of the Stackelberg mean field game problem.}
\label{subsec:model_SMFGrewrite}

As introduced in Subsection~\ref{subsec:model_SMFGsetup}, we have a principal that minimizes her own cost by choosing policies for the infinitesimal agents in the population and by taking their reaction to these policies into account. In order to find the optimal policies of the principal, we need to find the mean field Nash equilibrium in the population given any chosen policy by the principal. Then the principal should optimize her policies by taking into account the equilibrium behavior of the society. This creates a bi-level problem where at the higher level the principal minimizes her own cost and at the lower level the mean field Nash equilibrium in the population is found. However, in order to implement an efficient numerical approach, we need to reduce this bi-level problem to a single level problem.
\newpage
The bi-level problem is given as follows:  
\begin{equation*}
\left.\inf_{\substack{\boldsymbol\lambda\in\boldsymbol\Lambda}} {\EE\Big[}\int_0^T f_0(t, \mu_t^{\boldsymbol \lambda}, \lambda(t{,X^{\boldsymbol\lambda}_t}))dt +g_0( \mu_T^{\boldsymbol \lambda},\lambda(T{,X^{\boldsymbol\lambda}_T})){\Big]}\right\}\; \parbox{7em}{\textcolor{black}{\small Optimization of the Principal}}
\end{equation*}
\begin{equation*}
\left. \begin{aligned}
        X_t^{\boldsymbol \lambda} &= \zeta + \int_0^t b(s, X_s^{\boldsymbol \lambda}, \hat \alpha^{\boldsymbol \lambda}_s,  \mu_s^{\boldsymbol \lambda}; \lambda(s{,X^{\boldsymbol\lambda}_s})) ds + \int_0^t \sigma dW_s\\
        Y_t^{\boldsymbol \lambda} &= g( X^{\boldsymbol \lambda}_T, \mu^{\boldsymbol \lambda}_T; \lambda(t{,X^{\boldsymbol\lambda}_t})) + \int_t^T f(s,  X_s^{\boldsymbol \lambda}, \hat \alpha^{\boldsymbol \lambda}_s,  \mu_s^{\boldsymbol \lambda}; \lambda(s{,X^{\boldsymbol\lambda}_s}))ds - \int_t^T Z_s^{\boldsymbol \lambda}\cdot dW_s
    \end{aligned}\right\}\; \parbox{7em}{\textcolor{black}{\small Equilibrium in\\ the Population}}
\end{equation*}
where $ \mu_t^{\boldsymbol{\lambda}} = \mathcal{L}( X_t^{\boldsymbol{\lambda}})$ and $\hat{\alpha}_s^{\boldsymbol{\lambda}} := \hat{\alpha}_s(X_s^{\boldsymbol\lambda}, \mu_s^{\boldsymbol\lambda}, {\sigma^{-1}}^{\top}Z_s^{\boldsymbol \lambda}; \lambda(s{,X^{\boldsymbol\lambda}_s}))$.

As our first step, we write the backward SDE as a forward one and add a constraint to match the terminal condition of the backward SDE. In this way our problem can be written as
\begin{equation}
    \begin{aligned}
\label{eq:regulator_cost_original}   \inf_{\substack{\boldsymbol\lambda\in\boldsymbol\Lambda\\\boldsymbol Z\in \mathcal{H}^{2,m}\\ Y_0}} {\EE\Big[}\int_0^T f_0(t, \mu_t^{\boldsymbol\lambda, \boldsymbol Z, Y_0}, \lambda(t{,X^{\boldsymbol\lambda,\boldsymbol Z, Y_0}_t}))dt +g_0(\mu_T^{\boldsymbol \lambda, \boldsymbol Z, Y_0},\lambda(T{,X^{\boldsymbol\lambda,\boldsymbol Z, Y_0}_T})){\Big]},
    \end{aligned}
\end{equation}
with the constraint
\begin{equation}
    \begin{aligned}
    \label{eq:constraint_original}
        Y_T^{\boldsymbol \lambda, \boldsymbol Z, Y_0} = g(X^{\boldsymbol \lambda, \boldsymbol Z, Y_0}_T, \mu_T^{\boldsymbol \lambda, \boldsymbol Z, Y_0}; \lambda(T{,X^{\boldsymbol\lambda,\boldsymbol Z, Y_0}_T})),
    \end{aligned}
\end{equation}
where the trajectories of $\bsX^{\boldsymbol \lambda, \boldsymbol Z, Y_0}$ and $\boldsymbol Y^{\boldsymbol \lambda, \boldsymbol Z, Y_0}$ are determined by the following ``forward-forward" stochastic differential equations\footnote{Here, we use the term \textit{forward-forward} to emphasize that the backward stochastic differential equation is rewritten in forward direction. However, in order to avoid repetition, below we will just use the term \textit{forward stochastic differential equation}.}:
\begin{equation}
    \begin{aligned}
        \label{eq:FFSDE_original}
        X_t^{\boldsymbol \lambda, \boldsymbol Z, Y_0} &= \zeta + \int_0^t b(s, X_s^{\boldsymbol \lambda, \boldsymbol Z, Y_0}, \hat \alpha^{\boldsymbol \lambda, \boldsymbol Z, Y_0}_s, \mu_s^{\boldsymbol \lambda, \boldsymbol Z, Y_0}; \lambda(s{,X^{\boldsymbol\lambda,\boldsymbol Z, Y_0}_s})) ds + \int_0^t \sigma dW_s\\
        Y_t^{\boldsymbol \lambda, \boldsymbol Z, Y_0} &= Y_0^{\boldsymbol \lambda, \boldsymbol Z, Y_0} - \int_0^t f(s, X_s^{\boldsymbol \lambda, \boldsymbol Z, Y_0}, \hat \alpha^{\boldsymbol \lambda, \boldsymbol Z, Y_0}_s,\mu_s^{\boldsymbol \lambda, \boldsymbol Z, Y_0}; \lambda(s{,X^{\boldsymbol\lambda,\boldsymbol Z, Y_0}_s}))ds + \int_0^t Z_s \cdot dW_s,
    \end{aligned}
\end{equation}
where $\mu_t^{\boldsymbol \lambda, \boldsymbol Z, Y_0}= \mathcal{L} (X_t^{\boldsymbol\lambda, \boldsymbol Z^{\boldsymbol \lambda, \boldsymbol Z, Y_0}, Y_0^{\boldsymbol \lambda, \boldsymbol Z, Y_0}})$ and $\hat{\alpha}_s^{\boldsymbol{\lambda}, \boldsymbol Z, Y_0} := \hat{\alpha}_s(X_s^{\boldsymbol\lambda, \boldsymbol Z, Y_0}, \mu_s^{\boldsymbol\lambda, \boldsymbol Z, Y_0}, {\sigma^{-1}}^{\top}Z_s; \lambda(s{,X^{\boldsymbol\lambda,\boldsymbol Z, Y_0}_s}))$.

The system given in \eqref{eq:FFSDE_original} are the same equations given in \eqref{eq:Nash_FBSDE_sytem} where the dynamic of $\boldsymbol Y$ is written in the forward direction in time. Hence, the terminal point of the adjoint process, $Y_T$, is going to be determined by the initial point chosen $Y_0$ and the coefficient $\boldsymbol Z$ in front of the Wiener process. In order to penalize any deviation from the terminal condition, $Y_T = g(X_T, \mu_T;\lambda(T{, X_T}))$, a penalty term is added to the cost function as a second step of the rewriting of the problem. When this penalty term is equal to 0, it is implied that the FBSDE is solved and in turn the solution of the FBSDE characterizes the mean field Nash equilibrium in the agent population.

Previously, the principal had $\boldsymbol \lambda$ as their control. However, as it is stated before, now the terminal point of the adjoint process $Y_T$ is determined by $Y_0$ and $\boldsymbol Z$. Therefore, in the new stochastic optimal control problem, the controls are going to be $(\boldsymbol \lambda, \boldsymbol Z, Y_0)$.

Now, we introduce a penalty for not satisfying the constraint and we rewrite the above \textit{constrained} problem as follows:
\begin{equation}
    \begin{aligned}
    \label{eq:regulator_cost_with_penalty}
    & \inf_{\substack{\boldsymbol\lambda\in\boldsymbol\Lambda\\\boldsymbol Z\in \mathcal{H}^{2,m}\\ Y_0}} {\EE\Big[}\int_0^T f_0(t, \mu_t^{\boldsymbol \lambda, \boldsymbol Z, Y_0}, \lambda(t{,X^{\boldsymbol \lambda, \boldsymbol Z, Y_0}_t}))dt +g_0(\mu_T^{\boldsymbol \lambda, \boldsymbol Z, Y_0},\lambda(T{,X^{\boldsymbol \lambda, \boldsymbol Z, Y_0}_T})){\Big]}\\ &\hskip2cm + \nu \mathbb{E}\left[\mbsP\Big(Y^{\boldsymbol \lambda, \boldsymbol Z, Y_0}_T-g\big(X^{\boldsymbol \lambda, \boldsymbol Z, Y_0}_T, \mu_T^{\boldsymbol \lambda, \boldsymbol Z, Y_0}; \lambda(T{,X^{\boldsymbol \lambda, \boldsymbol Z, Y_0}_T})\big)\Big)\right],
    \end{aligned}
\end{equation}
where the trajectories of $\bsX^{\boldsymbol \lambda, \boldsymbol Z, Y_0}$ and $\boldsymbol Y^{\boldsymbol \lambda, \boldsymbol Z, Y_0}$ are determined by the forward SDE system given in \eqref{eq:FFSDE_original} and $\mu_t^{\boldsymbol \lambda, \boldsymbol Z, Y_0}= \mathcal{L} (X_t^{\boldsymbol\lambda, \boldsymbol Z, Y_0})$. Here $\mbsP(\cdot)$ is the penalty function such that $\mbsP(0)=0$ and $\mbsP(x)>0,$ for all $x\neq 0$.

In the next section, we discuss that the solution of the penalized optimal control problem given in \eqref{eq:regulator_cost_with_penalty} converges to the solution of the stochastic optimal control problem given in \eqref{eq:regulator_cost_original} with the constraint \eqref{eq:constraint_original}.

\section{Main theoretical results.}
\label{sec:main_theory}

In this section, we discuss the convergence results for the solution of the problem where a penalty term is added. Then, we discuss the convergence results for the \textit{parameterized} problem on which we will build the proposed numerical approach in Section~\ref{sec:numerics}.

\subsection{Reformulation of the original problem} 

In the rewritten problem, we will denote by $ \bsbeta:= (\boldsymbol\lambda, \boldsymbol z, y_0)$ the control. Let $\VV$ be the set of controls for the penalized problem, defined as: 
\begin{align*}
    \VV &= \Big\{ \bsbeta = (\boldsymbol \lambda, \boldsymbol z, y_0) \in C^0([0,T]{\times \RR^m}; \RR) \times C^0([0,T] \times \RR^m; \RR^d)  \times C^0(\RR^m; \RR) \,:\, 
    \\
    &\qquad\qquad\qquad{\hbox{for every $t \in [0,T]$, } \lambda(t, \cdot) \hbox{ is $L_{\lambda}$-Lipschitz continuous}}\\
    &\qquad\qquad\qquad  {\hbox { and has linear growth } |\lambda(t,x)|<L_{\lambda}(1+|x|), x \in \RR^m, }
    \\
    &\qquad\qquad\qquad  \hbox{for every $t \in [0,T]$, } z(t,\cdot) \hbox{ is $L_z$-Lipschitz continuous }
    \\
    &\qquad\qquad\qquad  \hbox { and has linear growth } |z(t,x)|<L_z(1+|x|), x \in \RR^m  \Big\}.
\end{align*}
We note that in the formulation of solution characterization of the mean field game, $Z_t$ behaves as the derivative of the value function (with respect to state $x$) of the representative player. Therefore $z(t,x)$ is a function of time $t$ and state $x$.
We will sometimes view $\bsbeta\in \VV$ as an element of $C^0([0,T] \times \RR^m \times \RR^m; \RR \times \RR^d \times \RR)$ by writing $\beta(t,x_0,x) = (\lambda(t{,x}), z(t,x), y_0(x_0))$. We endow $\VV$ with the sup norm $\|\cdot\|_{\VV} = \|\cdot\|_{\cC^0([0,T] \times \RR^m \times \RR^m; \RR \times \RR^d \times \RR)}$.

For the simplicity in notation, we also define 
\begin{equation}
    \label{eq:def-calb-calf}
    \mathcal{b}(t,x,\beta, \mu):= b(t, x, \hat{\alpha}_t(x, \mu, {\sigma^{-1}}^{\top}z; \lambda), \mu;\lambda), \qquad 
    \mathcal{f}(t, x, \beta, \mu):=f(t, x, \hat{\alpha}_t(x, \mu, {\sigma^{-1}}^{\top}z; \lambda), \mu;\lambda).
\end{equation}
For a given $\bsbeta = (\boldsymbol \lambda, \boldsymbol z, y_0)$, we consider the forward SDE system: 
\begin{equation}
    \begin{aligned}
        \label{eq:FFSDE_y0z}
        X_t^{\bsbeta} &=  X_0+ \int_0^t \mathcal{b}(s, X_s^{\bsbeta}, \beta_s,  \mu^{\bsbeta}_s) ds + \int_0^t \sigma dW_s, \qquad X_0 \sim \mu_0\\
        Y_t^{\bsbeta} &= y_0(X_0) - \int_0^t \mathcal{f}(s, X_s^{\bsbeta}, \beta_s, \mu_s^{\bsbeta})ds + \int_0^t z(s, X_s^{\bsbeta}) \cdot dW_s,
    \end{aligned}
\end{equation}
where $\mu^{\bsbeta}_s = \mathcal{L}({X^{\bsbeta}_s})$ and $\beta_s = \beta(s, X_0, X_s^{\bsbeta})=(\lambda(s{,X_s^{\bsbeta}}), z(s, X_s^{\bsbeta}), y_0(X_0))$.

We define the following penalty function as a function of $\bsbeta$:
$$
    \overline{\mbsP}(\bsbeta) =  \EE \left[ {\mbsP} \Big(Y_T^{\bsbeta} - {g}\big(X_T^{\bsbeta}, \mu_T^{\bsbeta}; \lambda(T{,X_T^{\bsbeta}})\big)\Big)\right].
$$
For $\bsbeta = (\boldsymbol \lambda, \boldsymbol z, y_0)$, we will still denote by $J^0$ the principal's cost in the new formulation, i.e.,
$$
    J^0(\bsbeta) = {\EE\Big[}\int_0^T f_0(t, \mu_t^{\bsbeta}, \lambda(t{,X_t^{\bsbeta}}))dt +g_0(\mu_T^{\bsbeta},\lambda(T{,X_T^{\bsbeta}})){\Big]}
$$
subject to the dynamics~\eqref{eq:FFSDE_y0z}.

We first start by reformulating the initial problem, with the controlled two forward equations for $X$ and $Y$. 

\noindent\textbf{Original Problem:} %
\begin{equation}
\label{eq:original-pb-J0}
    \begin{aligned}
        &\inf_{\bsbeta \in \VV} J^0\big(\bsbeta\big)\\
        &\text{s.t. \eqref{eq:FFSDE_y0z} holds and } Y_T^{\bsbeta} = {g}\big(X_T^{\bsbeta}, \mu_T^{\bsbeta};\lambda(T{,X_T^{\bsbeta}})\big).
    \end{aligned}
\end{equation}

We then introduce another problem, in which there is no terminal constraint but, instead, the total cost incorporates an extra term, weighted by a coefficient $\nu>0$, which penalizes the discrepancy between $Y_T^{\bsbeta}$ and ${g}\big(X_T^{\bsbeta}, \mu_T^{\bsbeta};\lambda(T{,X_T^{\bsbeta}})\big)$. We define the notation: 
$$
    J^0_\nu(\bsbeta) = J^0\big(\bsbeta\big) + \nu \overline{\mbsP}(\bsbeta)
$$
and we consider the following problem.

\noindent \textbf{Penalized Problem:}%
\begin{equation}
\label{eq:penalized-pb-J0P}
    \begin{aligned}
    &\inf_{\bsbeta\in \VV} J^0_\nu(\bsbeta)  \\ %
        &\text{s.t. \eqref{eq:FFSDE_y0z} holds}.
    \end{aligned}
\end{equation}

We let $\UU$ be the subset of $\VV$ of feasible controls for the original problem~\eqref{eq:original-pb-J0}, i.e., controls such that the terminal constraint is satisfied exactly.
\begin{remark}
By definition of the admissible policies $\boldsymbol\Lambda$ (see Definition~\ref{def:principal_policy_admissibility}), for every policy $\lambda$, there exists a unique Nash equilibrium for the MFG. Furthermore, under suitable regularity assumptions on the model, $y_0$ and $z$ obtained by the solution to the MKV FBSDE system are continuous and there exists a constant $L_z$ such that, for every $\lambda \in \Lambda$, $z$ is $L_z$-Lipschitz continuous in the spatial variable, hence there exists a pair $(y_0,z)$ which makes the penalty equal to $0$, which means that the MFG is solved exactly. 
\end{remark}

\vskip 6pt
\noindent\paragraph{Analysis of the penalized problem. } 
We now show that, if the coefficient of the penalty for the terminal condition matching increases to infinity, the optimal value of the penalized problem converges to the optimal value of the original problem, namely:
\begin{equation}
    \inf (J^0_{\nu_k}) \xrightarrow{\nu_k \rightarrow \infty} \inf (J^0),
\end{equation}
in a sense that we make precise below.

We will work under the following assumptions.

\begin{assumption}
\label{hyp2}
We assume:
\begin{enumerate}[label=\textbf{(A\theassumption.\arabic*)},align=left]
\item\label{hyp2:mu0} The fourth moment of the initial distribution is bounded i.e., $\int |x|^4d\mu_0(x)<\infty$. 
\item\label{hyp2:drift-lip} Remember we defined $\mathcal{b}(t, x,  \beta, \mu) := b(t, x, \hat\alpha_t(x, \mu, {\sigma^{-1}}^\top z), \mu; \lambda)$ where $\beta = (\lambda, z, y_0)$.  The function $\mathcal{b}$ is Lipschitz continuous with respect to $(x,\beta,\mu)$:
    \begin{equation}
        |\mathcal{b}(t, x^{\prime}, \beta^{\prime}, \mu^{\prime}) - \mathcal{b}(t, x, \beta, \mu)| 
        \leq  L_3 \big( |x^{\prime}-x| + |\beta^{\prime} - \beta| + W_2(\mu^{\prime}, \mu) \big).
    \end{equation}
\item\label{hyp2:cost0-lip} The running cost function, $f_0$, and the terminal cost function, $g_0$, of the principal are bounded from below by some constant and they are Lipschitz continuous with respect to $(\mu, \lambda)$:
    \begin{equation}
    \begin{aligned}
        |f_0(t,\mu^{\prime}, \lambda^{\prime}) - f_0(t,\mu, \lambda)| &\leq L_4 \big(W_2(\mu^{\prime}, \mu) + |\lambda^{\prime} - \lambda|\big)\\[2mm]
        |g_0(\mu^{\prime}, \lambda^{\prime}) - g_0(\mu, \lambda)| &\leq L_5 \big(W_2(\mu^{\prime}, \mu) + |\lambda^{\prime} - \lambda|\big)
    \end{aligned}
    \end{equation}
\item \label{hyp2:Pbar} The penalty function is of the form: $\mbsP(x)=x^2$ if $|x| \le C_P$ and $\mbsP(x)=C_P|x|$ otherwise, where $C_P$ is a positive constant.
\item \label{hyp2:g-lip} The terminal cost, $g$, of the agent is Lipschitz continuous in $(x,\mu,\lambda)$ such that
\begin{equation}
    |g(x^\prime,\mu^{\prime}; \lambda^{\prime}) - g(x, \mu; \lambda)| \leq L_6 \big(|x^\prime - x|+W_2(\mu^{\prime}, \mu) + |\lambda^{\prime} - \lambda|\big).
\end{equation}
\item\label{hyp2:f-lip} The function $\mathcal{f}$ is Lipschitz continuous with respect to $(x,\beta,\mu)$ such that: 
    \begin{equation}
        |\mathcal{f}(t, x^{\prime}, \beta^{\prime}, \mu^{\prime}) - \mathcal{f}(t, x, \beta, \mu)| 
        \leq  L_7 \big( |x^{\prime}-x| + |\beta^{\prime} - \beta| + W_2(\mu^{\prime}, \mu) \big).
    \end{equation}
\end{enumerate}
\end{assumption}

\begin{remark}
    In fact, for penalty function $\mbsP(\cdot)$, one can take any Lipschitz continuous function. However, later in our experiments (see Section~\ref{sec:experiments}), we will be using quadratic penalty function; therefore, in our theoretical analysis, we decided to focus on a penalty function which is quadratic on a compact set.
\end{remark}

\begin{remark}
    The Lipschitz property of $\mathcal{b}$ and $\mathcal{f}$, is satisfied for example when $b$ and $f$ are Lipschitz, as well as $\hat\alpha$. 
\end{remark}

\begin{theorem}[Main result 1: Convergence of penalized cost]
\label{thm:main-cv-thm-penalization}
	Suppose Assumption~\ref{hyp2} holds. Let $(\nu_k)_k$ be a sequence of positive numbers such that $\nu_k \to +\infty$. Let $(\epsilon_k)_k$ be a sequence of positive numbers such that $\epsilon_k \to \epsilon$. Let $(\bsbeta_k)_k$ be a sequence in $\VV$ such that, for every $k$, $\bsbeta_k$ is an $\epsilon_k$-minimizer for $J^0_{\nu_k}$. Then, every limit point of $(\bsbeta_k)_k$ is an $\epsilon$-minimizer for problem~\eqref{eq:original-pb-J0}. 
\end{theorem}
\begin{remark} In the above statement, \textit{limit point} is understood in the following sense of the sup-norm: $\bsbeta$ is a limit point of $(\bsbeta_k)_k$ if there is a subsequence $(\bsbeta_{k_j})_j$ such that
$$
    \|\bsbeta_{k_j} - \bsbeta\|_{\cC^0([0,T] \times \RR^m \times \RR^m)} \xrightarrow[j \to +\infty]{} 0.
$$

Furthermore, note that for a control $\bsbeta$, being an $\epsilon$-minimizer for problem~\eqref{eq:original-pb-J0} means, in particular, being an element of $\UU$ so that the terminal constraint is satisfied and, among the set of controls satisfying this constraint, $\bsbeta$ is $\epsilon$-optimal.
\end{remark}

\begin{remark}
The existence of a limit point can be guaranteed under suitable assumptions. For example, if the controls $(\bsbeta_k)_k$ are uniformly Lipschitz (with a common Lipschitz constant) and uniformly bounded, then the existence of a limit point can be obtained by applying Arzela-Ascoli theorem. 
\end{remark}

We now present the proof of Theorem~\ref{thm:main-cv-thm-penalization} by using two auxiliary results that are given and proved in Appendix~\ref{sec:proofs-lemma-penalization} for the sake of readability. 

\proof[Proof of Theorem~\ref{thm:main-cv-thm-penalization}. ]
Let $\overline{\bsbeta}$ be a limit point of $(\bsbeta_k)_k$. 
Up to the extraction of subsequences, we can assume that $\nu_k$ is an increasing sequence and $\epsilon_k$ is a decreasing sequence. 

\textbf{Step 1. } We first show that $\overline{\bsbeta} \in \UU$, i.e.,  the process controlled by $\overline{\bsbeta}$ satisfies the terminal constraint in the original problem~\eqref{eq:original-pb-J0}. 
Let $\bsbeta \in\UU$. By Lemma~\ref{lem:eps-min-ineq}, for every $k$, $J^0(\bsbeta) +\epsilon_{k}\geq J^{0}_{\nu_{k}}(\bsbeta_{{k}})\geq J^0(\bsbeta_{k})$. So, in particular,
\begin{equation}
\label{eq:ineq-nuk-J0k-Pk}
    \frac{1}{\nu_k}\left(J^0(\bsbeta) +\epsilon_{k} - J^0(\bsbeta_k)\right) 
    \ge \overline{\mbsP}(\bsbeta_k) \ge 0.
\end{equation}

Using the fact that $J^0$ is bounded from below (Assumption~\ref{hyp2:cost0-lip}) and using  Lemma~\ref{lem:cont-J-P}, which gives the continuity of $J^0$ and $\overline{\mbsP}$, we have that, as $k \to +\infty$, $\nu_k\to \infty$, $J^0(\bsbeta) +\epsilon_{k} - J^0(\bsbeta_k) \to J^0(\bsbeta) + \epsilon - J^0(\overline{\bsbeta})$, and $\overline{\mbsP}(\bsbeta_k) \to \overline{\mbsP}(\overline{\bsbeta})$. Hence $\overline{\mbsP}(\overline{\bsbeta}) = 0$, which means that $\overline{\bsbeta}$ satisfies the constraint of the original problem~\eqref{eq:original-pb-J0}.

\textbf{Step 2. } We now show that $\overline{\bsbeta}$ is an  $\epsilon$-optimizer of the original problem~\eqref{eq:original-pb-J0}. 

By Lemma~\ref{lem:eps-min-ineq} again, we have: for every $k$, $J^0(\bsbeta) +\epsilon_{k}\geq J^{0}_{\nu_{k}}(\bsbeta_{{k}}) \geq J^0(\bsbeta_{k})$. Passing to the limit in the left-hand side and the right-hand side, and using the continuity of $J^0$ provided by Lemma~\ref{lem:cont-J-P}, we obtain that: $J^0(\bsbeta) +\epsilon \geq J^0(\overline{\bsbeta})$. Since this holds for every feasible $\bsbeta$, $\overline{\bsbeta}$ is $\epsilon$-optimal.

\endproof

\vskip 6pt
\noindent
\paragraph{Finite population approximation. } 
Let $N$ be a positive integer. We consider an approximation of the MKV dynamics~\eqref{eq:FFSDE_y0z} with an interacting particle system. Consider the family of processes $(X^{N,i}_t,Y^{N,i}_t)_{t \in [0,T]}$, $i=1,\dots,N$, given by:
     \begin{equation}
     \label{eq:FFSDE-y0z-N}
     \left\{
    \begin{aligned}
        X_t^{N,i} &=  X_0^{N,i} + \int_0^t \mathcal{b}(s, X_s^{N,i}, \beta_s^{N,i},  \mu^{N,i}_s) ds + \int_0^t \sigma dW_s^i\\
        Y_t^{N,i} &= y_0(X^{N,i}_0) - \int_0^t \mathcal{f}(s, X_s^{N,i}, \beta_s^{N,i}, \mu_s^{N})ds + \int_0^t z(s, X_s^{N,i}) \cdot dW_s^i,
    \end{aligned}
    \right.
    \end{equation}
    with $X_0^{N,i} \sim \mu_0$ i.i.d. Here $\mu^{N}_t = \frac{1}{N} \sum_{i=1}^N \delta_{X^i_t}$ is the empirical distribution of $(X^{N,i}_t)_{i=1,\dots,N}$, with $\delta_x$ denoting a Dirac mass at $x$. Moreover, $\beta_s^{N,i} = \beta(s, X_0^{N,i}, X_s^{N,i})=(\lambda(s{,X_s^{N,i}}), z(s, X_s^{N,i}), y_0(X_0^{N,i}))$, and we recall that are $\mathcal{b}$ and $\mathcal{f}$ are defined in~\eqref{eq:def-calb-calf}.

We recall that the cost for the penalized problem is defined in~\eqref{eq:regulator_cost_with_penalty}. We define an analogous cost with the particle system approximation as: 
    \begin{equation}
    \label{eq:cost-JNnu}
    \begin{split}
        J^{0}_{\nu,N}(\beta) 
        &= \frac{1}{N} \sum_{i=1}^N \Big({\EE\Big[}\int_0^T f_0(t, \mu_t^{N}, \lambda(t{,X^{N,i}_t}))dt +g_0(\mu_T^{N},\lambda(T{,X^{N,i}_T})){\Big]}\\ &\hskip2cm + \nu \mathbb{E}\left[\mbsP\Big(Y^{N,i}_T-g\big(X^{N,i}_T, \mu_T^{N}; \lambda(T{,X^{N,i}_T})\big)\Big)\right] \Big)
    \end{split}
    \end{equation}
    We interpret this cost as a social cost over the population of particles, which justifies the average over $i=1,\dots,N$. However, note that removing the average over $i$ and considering only the cost for the particle $i=1$ does not change the value, since the particles are i.i.d.

We will use the following notation, which is related to propagation of chaos results (see \cite[Theorem 2.12]{carmona2018probabilistic2}): for any positive integer $d$, 
    \begin{equation} 
    \label{eq:chaos-epsilon-N}
        \epsilon_{d}(N) = 
        \begin{cases}
            N^{-1/2}, &\hbox{if $d<4$},\\
            N^{-1/2} \log(N), &\hbox{if $d=4$},\\
            N^{-2/d}, &\hbox{if $d>4$}.
        \end{cases}
    \end{equation}

\begin{theorem}[Main result 2: Finite population approximation]
\label{thm:finite-pop-approx}
    Suppose Assumption~\ref{hyp2} holds. 
    For all $\nu>0$, there exists $C$ depending only on the problem's data, $\nu$ and $L_\lambda$, $L_z$ in the definition of $\VV$, such that for all $\beta \in \VV$:
    \[
        |J^0_{\nu,N}(\beta) - J^0_{\nu}(\beta)| \le C \sqrt{\epsilon_{m+1}(N)}.
    \]
\end{theorem}
The proof is provided in Appendix~\ref{sec:proof-finite-pop-approx} and it relies on two main steps. First, we use a propagation of chaos type result to quantify the difference between trajectories of pairs $(X,Y)$ in the mean field dynamics and in the finite-population particle approximation. Then, we study the difference on the total expected cost induced by this approximation.

\vskip 6pt
\noindent
\paragraph{Discrete time approximation. } 

Next, since we want to use Monte Carlo simulations to estimate the cost, we discretize time using an Euler-Maruyama scheme. More complicated schemes could be used but we focus on this one for the sake of simplicity. We define the following discrete-time approximation of the $N$-particle system~\eqref{eq:FFSDE-y0z-N}:

    Let $N_T$ be a positive integer and consider a time grid $t_n = n \Delta t$, with $\Delta t = T/N_T$. 
    For a given discrete-time control $\ubsbeta = (\beta_{t_n})_{n=0,\dots,N_T}$, let $(X^{N,i}_0,Y^{N,i}_0)$ be i.i.d. with distribution $\nu_0 := \mathcal{L}(X_0,y_0(X_0))$ with $X_0 \sim \mu_0$, and for $n = 0 ,\dots,N_T-1$, let $(X^{N,i}_{t_{n+1}},Y^{N,i}_{t_{n+1}})$ be defined as:
    \begin{equation}
     \label{eq:FFSDE-y0z-N-NT}
     \left\{
    \begin{aligned}
        X_{t_{n+1}}^{N,i} &=  X_{t_{n}}^{N,i} + \mathcal{b}(t_n, X_{t_n}^{N,i}, \beta_{t_n}^{N,i},  \mu^{N,i}_{t_n}) \Delta t + \sigma \Delta W_{t_n}^i 
        \\[2mm]
        Y_{t_{n+1}}^{N,i} &= Y_{t_{n}}^{N,i} - \mathcal{f}(t_n, X_{t_{n}}^{N,i}, \beta_{t_{n}}^{N,i}, \mu_{t_{n}}^{N}) \Delta t + z(t_n, X_{t_{n}}^{N,i}) \cdot \Delta W_{t_{n}}^i,
    \end{aligned}
    \right.
    \end{equation}
    where $\Delta W^i_{t_n} = W^i_{t_{n+1}} - W^i_{t_n}$,  $i \in \llbracket N \rrbracket$,  are i.i.d. with distribution $\pzcN(0,\Delta t)$, and $\beta_{t_{n}}^{N,i} = \beta_{t_n}(X_0^{N,i}, X_{t_n}^{N,i})$. 
    We recall that $J^{0}_{\nu,N}$ is defined in~\eqref{eq:cost-JNnu}
and we define the discrete-time counterpart by:
\begin{align*}
    \check{J}^{0}_{\nu,N,N_T}(\ubsbeta) 
    &= \frac{1}{N} \sum_{i=1}^N \Big({\EE\Big[}\sum_{n=0}^{N_T-1} f_0(t_n, \check{\mu}_{t_n}^{N}, \lambda(t_n{,\check{X}^{N,i}_{t_n}})) \Delta t +g_0(\check{\mu}_T^{N},\lambda(T{,\check{X}^{N,i}_T})){\Big]}\\ &\hskip2cm + \nu \mathbb{E}\left[\mbsP\Big(\check{Y}^{N,i}_T-g\big(\check{X}^{N,i}_T, \check{\mu}_T^{N}; \lambda(T{,\check{X}^{N,i}_T})\big)\Big)\right] \Big).
\end{align*}

We will view this problem as an optimal control problem for $N$ particles with states $Q^{N,i}_t=(X^{N,i}_t,Y^{N,i}_t)$. 
We will use the following notation: for $t \in [0,T], x \in \RR^m, y \in \RR^m, \nu \in \mathcal{P}_2(\RR^{m+1}), \lambda \in \RR$:
\begin{equation}
    \label{eq:def-fC0-gC0}
    f^Q_0(t,\nu,\lambda) = f_0(t,\nu_x,\lambda), \quad
    g^Q_0((x,y), \nu,\lambda)
    = g_0(\nu_x,\lambda) + \nu \mbsP(y - g(x, \nu_x; \lambda)),
\end{equation}
where $\nu_x$ is the first marginal of $\nu$. 
Furthermore, we will denote by $
M_p(\mu) = \left( \int_{\RR^d} |x|^p d \mu(x) \right)^{1/p}$ the $p$-th moment of a probability measure $\mu$ (assuming it exists).
In order to analyze the impact of the time discretization, we will use the following assumptions, which are an adaptation to our setting of the ones used in the proof of~\cite[Proposition 15]{carmona2022convergence}: 
\begin{assumption}
\label{assm:fC0-gC0}
We have:
\begin{enumerate}[label=\textbf{(A\theassumption.\arabic*)},align=left]
    \item $\tilde{g}^Q_0(\nu,\lambda) = g^Q_0(\nu,\lambda)$ is differentiable with respect to both $\nu$ and $\lambda$, and its derivatives $\partial_x \tilde{g}^Q_0(q, \nu)$ and $\partial_\nu \tilde{g}^Q_0(q, \nu)(q')$ are of (at most) linear growth.\vskip1.5mm
    \item $\tilde{f}^Q_0(t, q, \nu) = f^Q_0(t, \nu, \Lambda^Q(t,q))$ is differentiable in all its variables, and its derivatives $\partial_x \tilde{f}^Q_0(t, q, \nu)$ and $\partial_\nu \tilde{g}^Q_0(t, q, \nu)(q')$ are of (at most) linear growth.\vskip1.5mm
    \item Denoting $\Theta = (t,\nu,\lambda) \in [0,T] \times \cP_2(\RR^{m+1}) \times \RR$, there holds:  
        \begin{itemize}
            \item $\Theta \mapsto |\partial_t \tilde{f}^Q_0 \bigl(\Theta\bigr)|$ has at most quadratic growth, 
            \item $\Theta \mapsto |\partial_\lambda \tilde{f}^Q_0\bigl(\Theta\bigr)|$, 
            and 
            $(\Theta,q) \mapsto |\partial_{\nu} \tilde{f}^Q_0\bigl(\Theta\bigr)(q)|$ have at most linear growth, and 
            \item 
            $\Theta \mapsto |\partial^2_{\lambda\lambda} \tilde{f}^Q_0 \bigl(\Theta\bigr)|$,  
            $(\Theta ,q) \mapsto |\partial_{\lambda}\partial_\mu \tilde{f}^Q_0\bigl(\Theta\bigr)(q)|$, 
            $(\Theta, q) \mapsto |\partial_v\partial_{\mu}\tilde{f}^Q_0\bigl(\Theta\bigr)(v)|_{v=q}|$, 
            $(\Theta, q) \mapsto |\partial^2_{\mu}f\bigl(\Theta\bigr)(q,q)|$ are bounded by a constant.\vskip1.5mm
        \end{itemize}
        \item $\mathbf{B}$ and $\mathbf{\Sigma}$ have at most linear growth: $|\mathbf{B}(t, q,\nu)| \le C (|t| + |q| + M_2(\nu) + |z(t,x)| + |\lambda(t,x)|)$ and $|\mathbf{\Sigma}(t, q)| \le C (|t| + |q| + |z(t,x)|)$.\vskip1.5mm
        \item\label{hyp:pontryagin-opt-fg-extra-mu0} The initial distribution $\mu_0$ is in $\cP_4(\RR^d)$, i.e., there exists a (finite) constant $C_{\mu_0}$ such that: 
        $
            M_4(\mu_0) = \int_{\RR^d} |x|^4 d \mu_0(x) \le C_{\mu_0}.
        $
\end{enumerate}
\end{assumption}
We refer to~\cite{carmona2022convergence} for discussions about these assumptions and how they relate to the MFG literature. Next, we state a convergence result for the time discretization.

\begin{theorem}[Main result 3: discrete time analysis]
\label{thm:discrete-time-approximation}
    Let $N$ and $N_T$ be positive integers. Let $\bsbeta \in \VV$. Let $\ubsbeta = (\beta(t_n))_{n =0,\dots,N_T}$ be the corresponding discrete time control obtained by evaluating $\bsbeta$ on the grid of time steps. 
    Suppose Assumption~\ref{assm:fC0-gC0} holds. We further assume that the control $\bsbeta$ satisfies: $\|\partial_t \bsbeta(t,x_0,x)\| \le C \|(t,x_0,x)\|$, $\|\partial_{(x_0,x)} \bsbeta(t,x_0,x)\| \le C$, $\|\partial^2_{(x_0,x),(x_0,x)} \bsbeta(t,x_0,x)\| \le C$. 
    Then:
    \[
        |J^0_{\nu,N,N_T}(\ubsbeta) - J^0_{\nu,N}(\bsbeta)|
        \le C \sqrt{\Delta t},
    \] 
    where the constant $C$ depends on the data of the problem, $\nu$, $\beta(0,0,0)$,  $L_\lambda$ and $L_z$ in the definition of $\VV$, and the constant appearing in the above assumption on $\bsbeta$ (but it is independent of $N$ and $N_T$). 
    
\end{theorem}

Note that the assumption on $\beta$ is satisfied under suitable assumptions on the MFG model which ensure, for instance, existence of a smooth enough classical solution to the Master equation; see \cite{carmona2022convergence} for more details. The proof of Theorem~\ref{thm:discrete-time-approximation} is provided in Appendix~\ref{sec:proof-discrete-time-approx} and consists of two steps: the analysis of the error induced by discretizing time in the trajectories of $Q=(X,Y)$ and then the impact on the cost.

\vskip 6pt
\noindent\paragraph{Neural network approximation. } 

The final step towards the implementation of our approach is to replace the set of controls by parameterized controls, with finite-dimensional parameters. Here, we will use neural networks, motivated by their performance in complex problems. Since the time has been discretized in the previous paragraph, we will approximate the control at each time step by a different neural network, function of the spatial variables only. The advantage of this approach is that we do not need to study the regularity of the control with respect to time. Alternatively, we could also consider neural networks functions of time and space. We will come back to this point in the numerical experiments (Section~\ref{sec:numerics}).

In the rest of this section, we will denote by $\bN^\psi_{d, n_{\mathrm{in}}, k}$ the class of feedforward neural networks with one hidden layer composed of $n_{\mathrm{in}}$ neurons, with activation function $\psi: \RR \to \RR$, taking a $k$-dimensional input and returning a one-dimensional output. For simplicity, we assume that the same activation function $\psi$ is used at every layer, and for the sake of definiteness, we shall assume that the activation function 
$\psi:\RR \to \RR$ is a $2\pi-$periodic function of class $\cC^3$, that is three times continuously differentiable, satisfying:
$
	\hat\psi_1:=\int_{-\pi}^\pi \psi(x) e^{-i  x} dx \neq 0.
$
This choice is due to the technical arguments that we use in the proof but similar arguments could work with other activation functions, yielding an approximation of the control function on any compact set. Furthermore, these restrictions do not seem to be necessary from the numerical viewpoint.

\begin{theorem}[Main result 4: Neural network approximation] \label{thm:aprpox-NN-timestep}

Suppose Assumption~\ref{hyp2} holds. Let $\bsbeta$ be as in Theorem~\ref{thm:discrete-time-approximation} and define the discrete-time control $\ubsbeta = (\beta(t_n,\cdot,\cdot))_{n=0,\dots,N_T}$. 
There exists three constants $n_0$, $K_1$ and $K_2$ depending on the data of the problem, on $N_T$, and on $\psi$ through $\hat\psi_1$ and $\|\psi'\|_{\cC^0(\TT^{m+1})}$ such that for each integer $n_{\mathrm{in}} > n_0$, there exists $\ubsbeta_\theta = (\beta_{t_n,\theta_n})_{n=0,\dots,N_T} \in (\bN^{\psi}_{{m+1}, n_{\mathrm{in}}, {d+2}})^{N_T+1}$ with parameters $\theta$ such that the Lipschitz constant of $\ubsbeta$ is bounded by $K_1$, and which satisfies: 
\begin{equation}
\label{eq:J0_P_bound_compact-discretetime-approx}
    |J^0_{\nu,N,N_T}(\ubsbeta)
    - J^0_{\nu,N,N_T}(\ubsbeta_\theta)| 
    \leq K_2 n_{\mathrm{in}}^{-1/(2(m+1))}.
\end{equation}
\end{theorem}
Notice that, in contrast with universal approximation results~\cite{cybenko1989approximation,hornik1990universal,hornik1991approximation,pinkus1999approximation}, this result provides a rate of convergence and not simply the existence of a neural network (with unbounded width) satisfying a given accuracy.   
The proof is provided in Appendix~\ref{sec:proof-main-thm-approx-NN}. Here again, it relies on two main steps. First, we show that at each time, the control can be approximated  on any compact set by a neural network. We then analyze how the error induced by this approximation propagates in the cost function.

\section{Numerical approach.}
\label{sec:numerics}
In this section, we propose the generalized single-level numerical method that make use of Monte Carlo simulation and machine learning tools to solve the Stackelberg MFG problem. Finding the optimal policies of the principal requires finding the corresponding Nash equilibrium since each policy choice affects the behavior of the population i.e., the Nash equilibrium. In the numerical approach, we generate Monte Carlo simulated particles and use the empirical distribution of these particles while training neural networks to approximate the optimal controls in the rewritten Stackelberg MFG problem \eqref{eq:regulator_cost_with_penalty}. The advantage of the proposed numerical approach is three-fold. Firstly, it is highly scalable when the dimensions of controls and states increase. Secondly, neural networks are good at approximating nonlinear functions which will help us to use more realistic complex models. Thirdly, it solves the bi-level problem as a single-level problem to increase efficiency.

After we discuss the algorithm to solve the rewritten version of the Stackelberg MFG problem, we will also discuss its extension to more complex models and then we will briefly consider a similar numerical approach for a special case that was introduced in \cite{StackelbergMFG_epidemics}.

\subsection{Penalty method.}
\label{subsec:numerics_penalty}

In order to solve the original Stackelberg MFG problem, we need to solve the optimization problem given in \eqref{eq:regulator_cost_with_penalty} under the dynamic constraints \eqref{eq:FFSDE_original}. Therefore, in this section, we focus on solving the optimal control problem given in \eqref{eq:regulator_cost_with_penalty} under the dynamic constraints \eqref{eq:FFSDE_original}.

\subsubsection{Monte carlo simulation.}
\label{subsubsec:numerics_penalty_MC}

 The optimization problem involves the state distribution of the population as the mean field interactions and this state distribution will be replaced by the empirical distribution obtained by simulating a system of interacting particles. In order to simulate the process trajectories of particles, we discretize the time integral. Therefore, in order to simulate the Monte Carlo trajectories of $(\bsX, \boldsymbol Y)$, we are going to use the time discretized versions of the expressions given in \eqref{eq:FFSDE_original}.

In the Monte Carlo simulation, on top of time discretization, we are using a finite number of particles to represent the agents in the population in order to obtain the empirical distribution in the population. The number of particles is taken as $N>0$ and we denote with $\llbracket N \rrbracket = \{1,\dots,N\}$ the set of indices of the particles. We consider $N_T$ time steps and denote $t_n = n \Delta t$ with $\Delta t = T/N_T$. We construct the Monte Carlo trajectories $(X^i_{t_n}, Y^i_{t_n})_{n,i\in \llbracket N\rrbracket}$ given measurable control functions $\lambda : [0,T]{\times \mathbb{R}^m}  \mapsto \mathbb{R}$, $z : [0,T] \times \mathbb{R}^m \mapsto \mathbb{R}^m$, $y_0 : \mathbb{R}^m\mapsto \mathbb{R}$ and the empirical state distribution $\bar \mu^N_{t_n}$. The empirical state distribution is calculated by using the sample and given as $\bar \mu^N_{t_n} = \frac{1}{N}\sum_{i=1}^N \delta_{X^i_{t_n}}$. After initialization, the iterations continue for every $n=1, \dots, N_T-1$. The updates over time steps are done by using the following discretized equations:

\begin{equation}
    \begin{aligned}
        \label{eq:discretized_FFSDE}
        X^i_{t_{n+1}} &= X^i_{t_n} + \Delta t \times b(t_n, X^i_{t_n}, \hat \alpha^i_{t_n}, \bar \mu^N_{t_n}; \lambda(t_n{,X^i_{t_n} }))+\sigma \Delta W^i_{t_n},\qquad && X^i_0 \sim \mu_0  \\[3mm]
        Y^i_{t_{n+1}} &= Y^i_{t_n} - \Delta t \times f(t_n, X^i_{t_n}, \hat \alpha^i_{t_n}, \bar \mu^N_{t_n}; \lambda(t_n{,X^i_{t_n} }))+z(t_n, X_{t_n}^i)\cdot \Delta W^i_{t_n}, && Y_0^i = y_0(X_0^i),
    \end{aligned}
\end{equation}
where $\hat \alpha^i_{t} = \hat{\alpha}(t, X^i_{t}, \bar\mu^N_t, {\sigma^{-1}}^{\top} z(t,X^i_t);\lambda(t{,X^i_{t} }))$ and $\Delta W^i_{t_n}\sim \pzcN(0,\Delta t)$ with $\Delta t = t_{n+1}-t_n$. The details of the Monte Carlo simulation of interacting particles can be seen in Algorithm~\ref{algo:simu-forward-XY}. We want to stress that this approach can be utilized for \textit{extended} mean field games where the mean field interactions come through the joint distribution of the control and state or the distribution of the control. In Section~\ref{subsubsec:experiments_contract_controlmean}, we showcase one such example.

\subsubsection{Approximation based on neural networks.}
\label{subsubsec:numerics_penalty_NNapprox}
In order to treat our problem numerically, we replace the controls $(\boldsymbol \lambda, \boldsymbol Z, Y_0)$ by parameterized functions $\lambda_{\theta_1}: [0,T] {\times \mathbb{R}^m}\to \mathbb{R}$, $z_{\theta_2}: [0,T] \times \mathbb{R}^m \to \mathbb{R}^d$, and $y_{0,\theta_3}: \mathbb{R}^m \to \mathbb{R}$ with respective parameters $\theta_1, \theta_2, \theta_3$. We extend the approach  similar to the numerical approach given in \cite{StackelbergMFG_epidemics}. Since the dimension of states can be potentially large and since we would like to work on complex models, we choose to use neural networks in the function approximation. To be more specific, in the implementation we take feed-forward fully connected deep neural networks. (In Section~\ref{subsec:numerics_extended}, we are going to consider the extension of our numerical approach to more complex models where we will be also using different neural network types such as recurrent neural networks.) Then, taking into account the time discretization and the approximation using a finite number of particles as described above, instead of~\eqref{eq:regulator_cost_with_penalty}, the goal is now to minimize over $\theta = (\theta_1,\theta_2,\theta_3)$ the following objective function:  
\begin{equation}
\label{eq:regulator_cost_with_penalty-NN}
\begin{aligned}
    {\mathbb{J}}^N(\theta)
    =
    & \, \EE \Bigg[
    \sum_{n=0}^{N_T-1} f_0\left(t_n, \bar{\mu}_{t_n}^{N, \theta}, \lambda_{\theta_1}(t_n{,X^{i,\theta}_{t_n}}) \right) \Delta t
    + g_0\left(\bar{\mu}_T^{N, \theta}, \lambda_{\theta_1}(T{,X^{i,\theta}_{T}})\right) \\&\hskip5cm +  \frac{1}{N}\sum_{i=1}^N \nu \left(Y^{i, \theta}_T- g\left(X^{i,\theta}_T, \bar{\mu}^{N,\theta}_T, \lambda_{\theta_1}(T{,X^{i,\theta}_{T}})\right)\right)^2 
    \Bigg],
\end{aligned}
\end{equation}
where $(\boldsymbol{X}^{i, \theta})_{i \in \llbracket N \rrbracket}$, $(\boldsymbol{Y}^{i, \theta})_{i \in \llbracket N \rrbracket}$, and $\boldsymbol{\bar{\mu}}^{N, \theta}$  are constructed by Algorithm~\ref{algo:simu-forward-XY} using $(\boldsymbol\lambda,\boldsymbol{z},y_0) = (\lambda_{\theta_1}, z_{\theta_2}, y_{0,\theta_3})$.  Here, the expected cost of the principal is estimated empirically by averaging over $M$ populations of size $N$.  
As a consequence, we will use a version of stochastic gradient descent in which one sample represents one population. Furthermore, in the numerical approach, we implement quadratic penalty function, which is similar to the penalty function described in Assumption~\ref{hyp2:Pbar} provided the constant $C_P$ is large enough.

To optimize over $\theta = (\theta_1,\theta_2,\theta_3)$, we rely on a variant of stochastic gradient descent (namely the Adaptive Moment Estimation algorithm) by sampling one population of size $N$ at each epoch. This kind of methods is well suited to the minimization of ${\mathbb{J}}^N$ in~\eqref{eq:regulator_cost_with_penalty-NN} since on the one hand the number of parameters in deep neural networks is potentially large and on the other hand this cost is written as an expectation and can thus be computed using Monte Carlo samples. This Monte Carlo simulation can be viewed as an extension of the simulation in~\cite{StackelbergMFG_epidemics} to the continuous state space. More precisely, we introduce for a sample $S = (X^i_{t_n},Y^i_{t_n},Z^i_{t_n})_{n=0,\dots,N_T, i \in \llbracket N \rrbracket}$ of a population of $N$ particles: 
\begin{equation}
\label{eq:regulator_cost_with_penalty-NN-sample}
\begin{aligned}
    {\mathbb{J}}_S^N(\theta)
    =
    & 
    \sum_{n=0}^{N_T-1} f_0\left(t_n, \bar{\mu}_{t_n}^{N, \theta}, \lambda_{\theta_1}(t_n{,X^{i,\theta}_{t_n}}) \right) \Delta t
    + g_0\left(\bar{\mu}_T^{N, \theta}, \lambda_{\theta_1}(T{,X^{i,\theta}_{T}})\right) \\&\hskip3cm + \frac{1}{N}\sum_{i=1}^N \nu \left(Y^{i, \theta}_T - g\left(X_T^{i,\theta}, \bar{\mu}_T^{N,\theta}, \lambda_{\theta_1}(T{,X^{i,\theta}_{T}})\right)\right)^2 
\end{aligned}
\end{equation}
where $\bar{\mu}^{N}_{t_n} = \frac{1}{N} \sum_{i=1}^N \delta_{X^i_{t_n}}$. Note that~\eqref{eq:regulator_cost_with_penalty-NN} is simply the average of $\mathbb{J}^N_{S}$ over samples (here one sample corresponds to the trajectories of one population). See Algorithm~\ref{algo:SGD-SMFG} for more details on the stochastic gradient descent (SGD) method applied in the context of Stackelberg MFG when the cost of the principal is rewritten.

In this approach, the goal is to minimize a cost that is the mixture of the principal's cost (first term in equation~\eqref{eq:regulator_cost_with_penalty}) and also the error of not matching the terminal condition of the backward equation (second term in equation~\eqref{eq:regulator_cost_with_penalty}). As it shown in Section~\ref{sec:main_theory}, when the weight given to second term (i.e., $\nu$) is increased, we end up with having this term converging to 0 i.e., we manage to match the terminal condition of the backward equation. This results in solving for the equilibrium in the population. With the minimization of the first term the algorithm finds the regulator's optimal cost while solving for the equilibrium simultaneously in a single-level optimization.

\newpage

\subsubsection{Algorithms.}
\label{subsubsec:numerics_penalty_algo}

\begin{algorithm}[H]
\caption{Monte Carlo simulation of an interacting batch \label{algo:simu-forward-XY}}

\textbf{Input:} number of particles $N$; time horizon $T$; time increment $\Delta t$; initial distribution $\mu_{0}$; control functions $\lambda, z, y_0$

\textbf{Output:} Approximate sampled trajectories of $(\bsX,\boldsymbol Y,\boldsymbol Z)$ solving~\eqref{eq:discretized_FFSDE}\vskip1mm
\begin{algorithmic}[1]
\STATE{Let $n=0$, $t_0 = 0$; pick $X^i_0 \sim \mu_{0}$ i.i.d. and set $Y^i_0 = y_0(X^i_0), \forall i \in \llbracket N \rrbracket$} \vskip1mm
\WHILE{$k \times \Delta t = t_n < T$}\vskip1mm
    \STATE{Set $ Z^i_{t_n} = z({t_n}, X^i_{t_n}),\ \Lambda_{t_n}= \lambda(t_n{, X^i_{t_n}}),\ i \in \llbracket N \rrbracket$}\vskip1mm
    \STATE{Let $\bar \mu^N_{t_n} = \frac{1}{N} \sum_{i=1}^N \delta_{X^i_{t_n}}$ }\vskip1mm
    \STATE{Set $\hat \alpha^i_{t_n} = \hat{\alpha}({t_n},X^i_{t_n},\bar{\mu}^N_{t_n}, \sigma^{{-1}^{\top}}Z^i_{t_n}; \Lambda_{t_n}),\ i \in \llbracket N \rrbracket$}\vskip1mm
    \STATE{Sample $\Delta W^i_{t_n} \sim \pzcN(0,\Delta t),\ i \in \llbracket N \rrbracket$ }
    \vskip1mm
    \STATE{Let 
        $
            X^{i}_{t_{n+1}} = X^{i}_{t_n} + b(t_n, X^{i}_{t_n}, \hat\alpha^{i}_{t_n}, \bar\mu^N_{t_n}; \Lambda_{t_n})\Delta t + \sigma \Delta W^{i}_{t_n}
        $, $ i \in \llbracket N \rrbracket$}\vskip1mm
    \STATE{Let 
        $
            Y^{i}_{t_{n+1}} = Y^{i}_{t_n} - f(t_n, X^{i}_{t_n}, \hat \alpha^{i}_{t_n}, \bar\mu^N_{t_n}; \Lambda_{t_n})\Delta t + (Z^{i}_{t_n})^{\top} \Delta W^{i}_{t_n}
        $, $ i \in \llbracket N \rrbracket$
    }\vskip1mm
    \STATE{Set $t_n = t_n+\Delta t$ and $n=n+1$}\vskip1mm
\ENDWHILE\vskip1mm
\STATE{Set $N_T = k,\ t_{N_T} = T$}\vskip1mm
\RETURN $(X^i_{t_n},Y^i_{t_n},Z^i_{t_n})_{n =0,\dots,N_T,\ i \in \llbracket N \rrbracket}$ and $(t_n)_{n =0,\dots,N_T}$\vskip1mm
\end{algorithmic}
\end{algorithm}

\begin{algorithm}[H]
\caption{SGD for Stackelberg MFG\label{algo:SGD-SMFG}}

\textbf{Input:} Initial parameter $\theta_0$; number of iterations ${K}$; sequence $(\beta_{{k}})_{{k} = 0, \dots, {K}-1}$ of learning rates; number of particles $N$; time horizon $T$; time increment $\Delta t$; initial distribution $\mu_0$

\textbf{Output:} Approximation of $\theta^*$ minimizing $\mathbb{J}^N$ defined by ~\eqref{eq:regulator_cost_with_penalty-NN}\vskip1mm
\begin{algorithmic}[1]

\FOR{${k} = 0, 1, 2, \dots, {K}-1$}\vskip1mm
    \STATE{Sample $S = (X^i_{t_{{n}}},Y^i_{t_{{n}}},Z^i_{t_{{n}}})_{{n} =0,\dots,{N_T},\ i \in \llbracket N \rrbracket}$ and $(t_{{n}})_{{n}=0,\dots,{N_T}}$ using Algorithm~\ref{algo:simu-forward-XY} with control functions $(\boldsymbol\lambda, \boldsymbol z, y_0) = (\lambda_{\theta_{{k},1}}, z_{\theta_{{k},2}}, y_{0,\theta_{{k},3}})$ and parameters: number of particles $N$; time horizon $T$; time increment $\Delta t$; initial distribution $\mu_{0}$ }\vskip1mm
    \STATE{Compute the gradient $\nabla \mathbb{J}^N_{S}(\theta_{{k}})$ of $\mathbb{J}^N_{S}(\theta_{{k}})$ defined by~\eqref{eq:regulator_cost_with_penalty-NN-sample} }\vskip1mm
    \STATE{Set $\theta_{{k}+1} =  \theta_{{k}} -\beta_{{k}} \nabla \mathbb{J}^N_{S}(\theta_{{k}})$ }\vskip1mm
\ENDFOR

\RETURN $\theta_{{K}}$
\end{algorithmic}
\end{algorithm}

\subsection{Extension to more complex models.}
\label{subsec:numerics_extended}
\subsubsection{General case.}
\label{subsubsec:numerics_extended_general}
In this section, we solve stylized models, some of which are beyond the scope covered by the theoretical analysis presented in the previous section, to show the flexibility of our numerical approach. We will illustrate our numerical method on models that are widely used in contract theory where the principal proposes a contract with a time-dependent policy and also a path-dependent terminal incentive. In other words, the principal has an additional $\mathcal{F}_T$-measurable terminal incentive $\xi$ in addition to the time dependent Markovian policy $\boldsymbol\lambda=(\lambda(t{, X_t}))_{t \in [0,T]}$ that is introduced before. The cost functions for the representative agent and the principal in the new model are given below.
The representative agent pays the following cost $J^{\boldsymbol\lambda,\xi}(\boldsymbol\alpha, \boldsymbol\mu)$ when she chooses the control $\boldsymbol\alpha \in \mathbb{A}$ and interacts with the population through the mean field $\boldsymbol\mu \in M(\mathcal{P}(\mathbb{R}^m))$:
\begin{equation}
\label{eq:cost_of_minors_xi}
    J^{\boldsymbol\lambda,\xi}(\boldsymbol\alpha, \boldsymbol\mu)
    := 
    \mathbb{E}\left[
    \int_0^T f(t, X_{t}, \alpha_t, \mu_t; \lambda(t{,X_t}))dt + g(X_T, \mu_T; \lambda(T{,X_T}))- U(\xi)
    \right],
\end{equation}
where $f$ and $g$ are the running cost and the terminal cost of the agent that are defined in Subsection~\ref{subsec:model_SMFGsetup} and $U : \mathbb{R}\mapsto\mathbb{R}$ is a utility function that represents the agent's utility from the terminal payment. The representative agent's state dynamics are the same as in \eqref{eq:dynamics_of_minors}.

The principal pays the following cost when she chooses the admissible contract $(\boldsymbol\lambda, \xi)$:
\begin{equation}
\label{eq:cost_of_principal_xi}
    J^0(\boldsymbol \lambda, \xi) := \mathbb{E}\left[\int_0^T f_0(t, \hat \mu^{\boldsymbol \lambda, \xi}_t, \lambda(t{,X^{\boldsymbol \lambda, \xi}_t}))dt +g_0(\hat \mu_T^{\boldsymbol \lambda, \xi}, \lambda(T{,X^{\boldsymbol \lambda, \xi}_T})) +\xi\right],
\end{equation}
where $f_0$ and $g_0$ are the running and terminal cost of the agent that are defined in Subsection~\ref{subsec:model_SMFGsetup} and $(\hat{\boldsymbol\alpha}^{\boldsymbol\lambda, \xi}, \hat{\boldsymbol\mu}^{\boldsymbol\lambda, \xi}) \in \mathcal{N}(\boldsymbol\lambda, \xi)$. Here, $\mathcal{N}(\boldsymbol\lambda, \xi)$ denotes the set of Nash equilibrium given the principal's contract choice $(\boldsymbol\lambda, \xi)$. Additionally, now the principal pays a $\mathcal{F}_T$-measurable terminal payment $\xi$ to the agent. 

The last aspect of the extended problem is giving a walk-away option to the infinitesimal agents that is commonly implemented in the contract theory literature. In other words, we require that the cost of the agents should not be higher than a preset reservation threshold denoted by $\kappa$. Therefore, all Nash equilibria in which the representative agent's expected total cost is higher than $\kappa$ are disregarded. We explain how to implement this reservation threshold in our numerical approach shortly. By following the ideas we used in Section~\ref{subsec:model_SMFGrewrite}, we can rewrite the problem as follows:

\begin{equation}
    \begin{aligned}
    \label{eq:rewritten_cost_Wpenalty_complex}
    & \inf_{Y_0:\mathbb{E}[Y_0]\leq\kappa}\inf_{\substack{\boldsymbol Z\in \mathcal{H}^{2,m}\\ (\boldsymbol\lambda, \xi)}} \mathbb{E}\left[\int_0^T f_0(t, \mu_t^{\boldsymbol \lambda, \boldsymbol Z, Y_0, \xi}, \lambda(t{,X^{\boldsymbol \lambda, \boldsymbol Z, Y_0, \xi}_t}))dt +g_0(\mu_T^{\boldsymbol \lambda, \boldsymbol Z, Y_0, \xi},\lambda(T{,X^{\boldsymbol \lambda, \boldsymbol Z, Y_0, \xi}_T}))+ \xi\right]\\ &\hskip2.5cm + \nu \mathbb{E}\left[\Big(Y^{\boldsymbol \lambda, \boldsymbol Z, Y_0, \xi}_T-g\big(X^{\boldsymbol \lambda, \boldsymbol Z, Y_0, \xi}_T, \mu_T^{\boldsymbol \lambda, \boldsymbol Z, Y_0, \xi}; \lambda(T{,X^{\boldsymbol \lambda, \boldsymbol Z, Y_0, \xi}_T})\big) + U(\xi)\Big)^2 \right],
    \end{aligned}
\end{equation}
under the dynamic constraints \eqref{eq:FFSDE_original} and $\mu_t^{\boldsymbol\lambda, \boldsymbol Z, Y_0, \xi}=\mathcal{L}(X_t^{\boldsymbol\lambda, \boldsymbol Z, Y_0, \xi})$. Realize that in the rewritten problem, the controls are $(\boldsymbol \lambda, \boldsymbol Z, Y_0, \xi)$.

Similar to the numerical algorithm proposed in Section~\ref{subsec:numerics_penalty}, we simulate the trajectories of the dynamical constraint \eqref{eq:FFSDE_original} by discretizing the time, and for the distribution, we use the empirical distribution obtained by simulating a system of finite number of interacting particles. In order to solve the coupled forward backward equations, we write the backward equation forward in time. Our aim is to approximate control functions $(\boldsymbol \lambda, \boldsymbol Z, Y_0, \xi)$ by minimizing over $\theta = (\theta_1, \theta_2, \theta_3, \theta_4)$ the discretized version of the cost \eqref{eq:rewritten_cost_Wpenalty_complex}:
\begin{equation}
\label{eq:regulator_cost_with_penalty-NN_complex}
\begin{aligned}
    {\check{\mathbb{J}}}^N(\theta)
    =
    & \frac{1}{M} \sum_{j=1}^M \Bigg[ \frac{1}{N}\sum_{i=1}^N\Big[
    \sum_{n=0}^{N_T-1} f_0\left(t_n, \bar{\mu}_{t_n}^{j, N, \theta}, \lambda^{j}_{\theta_1}(t_n{,X^{j,i,\theta}_{t_n}}) \right) \Delta t
    + g_0\left(\bar{\mu}_T^{j, N, \theta}, \lambda^{j}_{\theta_1}(T{,X^{j,i,\theta}_{T}})\right) + \xi^{j,i}_{\theta_4} \\&\hskip3cm +   \nu \left(Y^{j,i, \theta}_T- g\left(X^{j,i,\theta}_T, \bar{\mu}^{j,N,\theta}_T, \lambda^{j}_{\theta_1}(T{,X^{j,i,\theta}_{T}})\right)+U(\xi^{j,i}_{\theta_4})\right)^2 \Big]
    \Bigg],
\end{aligned}
\end{equation}
where $\bar \mu^{j,N,\theta}_{t_n} = \frac{1}{N}\sum_{i=1}^N \delta_{X^{j,i,\theta}_{t_n}}$.

In order to approximate a function for $\xi$, we make use of recurrent neural networks (RNNs) that use the whole trajectory of the state of agents $(\bsX^{\boldsymbol\lambda, \boldsymbol Z, Y_0, \xi})$ i.e., $\xi_{\theta_4} : (\RR^m)^{N_T+1} \mapsto \mathbb{R}$. To optimize over $\theta = (\theta_1,\theta_2,\theta_3, \theta_4)$, we rely on a variant of stochastic gradient descent (namely the Adaptive Moment Estimation algorithm) by sampling one population of size $N$ at each epoch as we did in Section~\ref{subsubsec:numerics_penalty_NNapprox}. We note that RNNs have been used in the context of MFGs in~\cite{fouque2020deep,han2021recurrent,gomes2022machine} to represent controls that can depend on trajectories. However, here we are using an RNN for the principal's policy. To be specific, in our implementation, we use long-short term memory (LSTM) neural networks, which have are better than plain RNNs to avoid vanishing or exploding gradients. Its input is a realization of the trajectory of $\bsX^{\boldsymbol\lambda, \boldsymbol Z, Y_0, \xi} \in (\RR^m)^{N_T+1}$, and its output is the value of the terminal payment $\xi \in \RR$. We implement the reservation threshold by constraining the neural network output for $y_{0,\theta_3}$ to be smaller than the threshold $\kappa$. In the implementation, this can be done by introducing a proper activation function for the output layer of the neural network. For example, for getting outputs smaller than or equal to 0, we can use rectified linear unit (ReLU) activation function in the output layer and we can scale it with -1.

\begin{remark}
The problem described above is more complex than the setting covered by the theoretical analysis provided in the previous sections. The reason why we choose the present problem in our numerical experiments is to show that our approach can easily be adapted to tackle more complex models, beyond the scope of the assumptions made for the sake of the theoretical analysis. 
We leave the analysis of the present setting for future research. 
\end{remark}

\begin{remark}
    The proposed numerical approach can be extended to even more complex settings. For example, if the population is subject to common noise, then one can still use an empirical population to approximate the \emph{conditional} distribution. The main technical difficulty is that, in general, the players' controls need to take the distribution as an input and not just the individual state. In some classes of models, the optimal control is a function of sufficient statistics, see e.g.~\cite{carmona2022convergence,min2021signatured}. But in general, one needs to design a suitable finite-dimensional approximation of the mean-field distribution. Typical approximations include histograms, moments, or symmetric functions of the empirical distributions, see e.g.~\cite{dayanikli2023deep} for more details on these approximations in the context of mean field control with common noise. More generally, the question of learning population-dependent controls (even without common noise) has recently attracted interest~\cite{perrin2022generalization,germain2022deepsets}, and these methods could also be adapted to tackle Stackelberg MFGs with common noise. 
\end{remark}

\subsubsection{Special case.}
\label{subsubsec:numerics_extended_special}

In this section, we focus on a special case of our problem and present a numerical approach different than the penalty method. This numerical approach is introduced in~\cite{StackelbergMFG_epidemics} for finite state Stackelberg MFGs and here, we present its continuous state implementation.

\begin{assumption}[Regularity of the representative player's terminal cost] 
\label{assu:sannikov_trick}
\item[i.]  $g(X_T, \mu_T; \lambda(T{, X_T}))=0$
\item[ii.] The utility function $U:\mathbb{R}\mapsto \mathbb{R}$ is invertible.
\end{assumption}

With the Assumption~\ref{assu:sannikov_trick}, the terminal condition of the adjoint process, $Y_T$ is equal to $-U(\xi)$. Since the utility function $U(\cdot)$ is assumed to be invertible, $\xi = U(-Y_T)$. Now, we consider the following optimal control problem:

\begin{equation}
    \begin{aligned}
    \label{eq:rewritten_regulator_cost_sannikov}
    \inf_{Y_0:\mathbb{E}[Y_0]\leq\kappa}\inf_{\substack{\boldsymbol Z\in \mathcal{H}^{2,m}\\ \boldsymbol\lambda\in\boldsymbol\Lambda}} \mathbb{E}\left[\int_0^T f_0(t,  \mu^{\boldsymbol \lambda, \boldsymbol Z, Y_0}_t, \lambda(t{, X^{\boldsymbol \lambda, \boldsymbol Z, Y_0}_t}))dt +g_0( \mu_T^{\boldsymbol \lambda, \boldsymbol Z, Y_0},\lambda(T{, X^{\boldsymbol \lambda, \boldsymbol Z, Y_0}_T})) + U^{-1}(-Y^{\boldsymbol \lambda, \boldsymbol Z, Y_0}_T)\right],
    \end{aligned}
\end{equation}
under the dynamic constraints \eqref{eq:FFSDE_original}.

Previously, the principal had $(\boldsymbol \lambda, \xi)$ as their controls. However, now $\xi$ is going to be determined by $U^{-1}(-Y^{\boldsymbol \lambda, \boldsymbol Z, Y_0}_T)$. Furthermore, the terminal point $Y^{\boldsymbol \lambda, \boldsymbol Z, Y_0}_T$ is going to be determined by the initial point chosen $Y_0$ and the coefficient $\boldsymbol Z$ in front of the Wiener process. Therefore, in the new optimal control problem, the controls are going to be $(\boldsymbol \lambda, \boldsymbol Z, Y_0)$. Implementation of the algorithm follows similarly the ideas introduced in the previous sections. Our aim is to approximate functions $(\boldsymbol \lambda, \boldsymbol Z, Y_0)$ by minimizing over $\theta = (\theta_1, \theta_2, \theta_3)$ the discretized version of the cost \eqref{eq:rewritten_regulator_cost_sannikov}:

\begin{equation}
\label{eq:regulator_cost_with_penalty-NN_sannikov}
\begin{aligned}
    {\tilde{\mathbb{J}}}^N(\theta)
    =&
    \frac{1}{M} \sum_{j=1}^M \Bigg[{\frac{1}{N} \sum_{i=1}^N\Big[}
    \sum_{n=0}^{N_T-1} f_0\left(t_n, \bar{\mu}_{t_n}^{j, N, \theta}, \lambda^{j}_{\theta_1}(t_n{,X_{t_n}^{j,i, \theta}}) \right) \Delta t
    \\&\hskip4cm + g_0\left(\bar{\mu}_T^{j, N, \theta}, \lambda^{j}_{\theta_1}(T{,X_{T}^{j,i, \theta}})\right) + U^{-1}(Y_T^{j,i, \theta}) {\Big]}\Bigg]
\end{aligned}    
\end{equation}

Again, to optimize over $\theta = (\theta_1,\theta_2,\theta_3)$, we rely on a variant of stochastic gradient descent (namely the Adaptive Moment Estimation algorithm) by sampling one population of size $N$ at each epoch as we did in Section~\ref{subsubsec:numerics_penalty_NNapprox}. 

Since the $\xi$ will be determined by $Y_0$ and $\boldsymbol Z$ instead of being approximated with a neural network as it is implemented in the generalized penalty method (Section~\ref{subsubsec:numerics_extended_general}), this algorithm is expected to converge faster. However, the limitation of this algorithm is that it can be only used for the problems in which, at time $T$, the principal pays $\xi$ and the agent receives $U(\xi)$ for some invertible function $U$ (see Assumption~\ref{assu:sannikov_trick}). In the next section, we will consider an example where the model is in this special form to show the results of this special algorithm and also the generalized penalty algorithm.

\section{Experimental results.}
\label{sec:experiments}
In this section, we test numerical approaches presented in the previous section against the models in the literature with explicit solutions. First, we focus on the three versions of the explicitly solvable example analyzed in \cite{Elie_2019}. Later, we extend the systemic risk model analyzed in \cite{carmona2018probabilistic} by adding a principal and present our numerical results.

\subsection{Optimal contract between a principal and many many agents.}
\label{subsec:experiments_contract}

We first start by presenting the explicitly solvable contract theory model analyzed in \cite{Elie_2019} between a principal and a mean field of agents. Agents control their effort level $\alpha_t \in \mathbb{R}_+$ and their state is the value of a project $X_t \in \mathbb {R}$. The dynamics of the value of the project is given as 

\begin{equation}
    \begin{aligned}
        dX_t = \Big(\alpha_t + aX_t+ \beta_1 \int_\mathbb{R} x d \mu_t(x) + \beta_2 \int_\mathbb{R} \alpha d q_t(\alpha) - \gamma V_\mu (t)\Big) dt + dW_t
    \end{aligned}
\end{equation}

where $V_\mu (t) = \int_\mathbb{R} |x|^2 d \mu_t(x) - \left|\int_\mathbb{R} x d \mu_t(x)\right|^2$. The principal only has the $\mathcal{F}_T$-measurable incentive that is paid at the terminal time to agents, she does not have a time dependent policy.

The running cost of the representative agent is given as $f(t, X_t, \alpha_t, \mu_t; \lambda(t)) = k \frac{\alpha_t^2}{2}$ where $k>0$, the terminal cost is 0. The representative agent is taken to be risk neutral; in other words, his terminal utility is linear i.e., $U(\xi)=\xi$. Therefore, the representative agent's cost function is as follows:
\begin{equation}
    \mathbb{E}\left[\int_0^T k \frac{\alpha_t^2}{2} dt - \xi\right].
\end{equation}

We further assume the principal is risk neutral and has the following objective:
\begin{equation}
    \inf_{\xi}\mathbb{E}\big[\xi-X_T\big].
\end{equation}

As we stated before, in this example the principal only has one control, the terminal payment $\xi$. Later, in our next example, we consider a case where we have the time dependent policy term, $\boldsymbol\lambda$.

We write the FBSDE that characterize mean field equilibrium in the agent population by following the ideas introduced in Section~\ref{subsec:model_mfganalysis}. The Hamiltonian for the mean field game of agents is 

\begin{equation*}
    \check{H}(t, X, \mu, q, Z, \alpha) = \Big(\alpha + a X+ \beta_1 \int_\mathbb{R} x d \mu_t(x) + \beta_2 \int_\mathbb{R} \alpha d q_t(\alpha) - \gamma V_\mu (t)\Big)Z+ \frac{k}{2} \alpha^2,
\end{equation*}

and the minimizer of the Hamiltonian is
\begin{equation}
    \hat \alpha_t = -\dfrac{Z_t}{k}.
\end{equation}

Here, we stress that the arguments of the Hamiltonian of this example are slightly different than the Hamiltonian introduced in Section~\ref{subsec:model_mfganalysis}. Now in addition we have interactions through the control distribution; however, we don't have the time dependent policy of the principal.

The system of FBSDE that characterize the equilibrium in the mean field population is given as:

\begin{equation*}
    \begin{aligned}
        dX_t &= \Big(-\dfrac{Z_t}{k}+ a X_t +\beta_1 \int_\mathbb{R} x d\mu_t(x) -\beta_2 \dfrac{\bar Z_t}{k}- \gamma V_\mu(t)\Big)dt + dW_t\\
        dY_t &= -\dfrac{k}{2}\left(-\dfrac{Z_t}{k}\right)^2 dt + Z_tdW_t.
    \end{aligned}
\end{equation*}

\begin{remark}
Optimal effort of the agent is deterministic and given by:
\begin{equation}
    \hat\alpha_t = (1+\beta_2) \dfrac{e^{(a+\beta_1)(T-t)}}{k}.
\end{equation}
Interested readers can refer to~\cite{Elie_2019} for the analysis of explicit solution.
\end{remark}

\subsubsection{Experiment 1: interactions through the variance of the states.}
\label{subsubsec:experiments_contract_statevar}

First, we consider a case where the interactions among agents are through the variance of the states, i.e., we take $\beta_1=\beta_2=0$. In this case, the optimal effort of an agent is given as
\begin{equation}
    \hat \alpha_t = \dfrac{e^{a(T-t)}}{k}.
\end{equation}

\begin{figure}[h]
    \centering
    \begin{subfigure}{1\textwidth}
        \centering
        \includegraphics[width=1\textwidth]{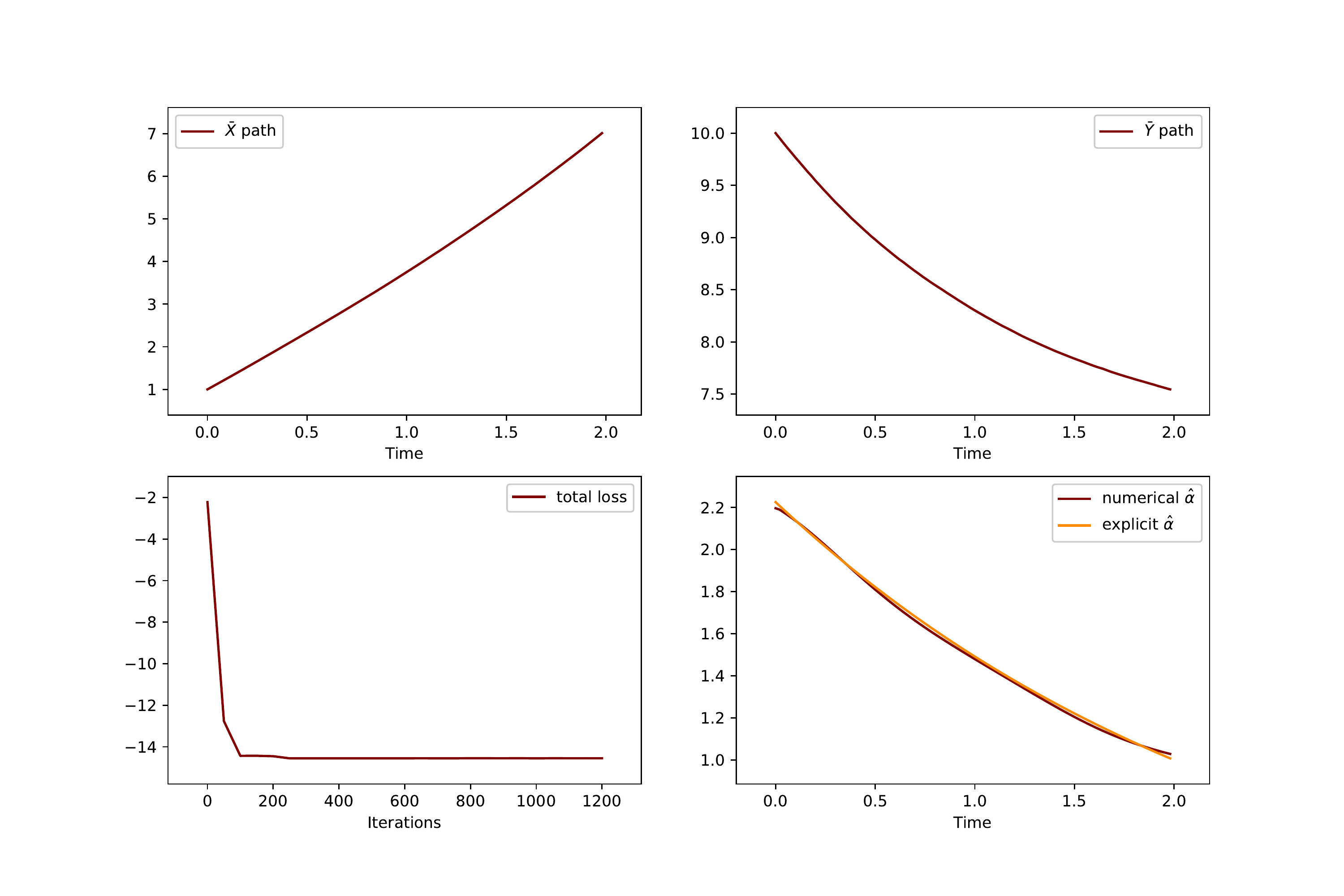}
    \end{subfigure}
    \caption{Results of the contract theory experiment where we have interactions through variance of the state and where we use the algorithm proposed for the special case in Section~\ref{subsubsec:numerics_extended_special}. Top: the paths of $\bar X$ (left) and $\bar Y$ (right), bottom: total loss (left) and the comparison of the explicit effort level to the mean of the effort levels found by using the special numerical approach presented in Section~\ref{subsubsec:numerics_extended_special} (right). It is seen that the special numerical approach is very good at mimicking the explicit results.}
    \label{fig:contract_exp1_sannikov}
\end{figure}

\begin{figure}[h]
     \centering
     \begin{subfigure}{1\textwidth}
         \centering
         \includegraphics[width=\textwidth]{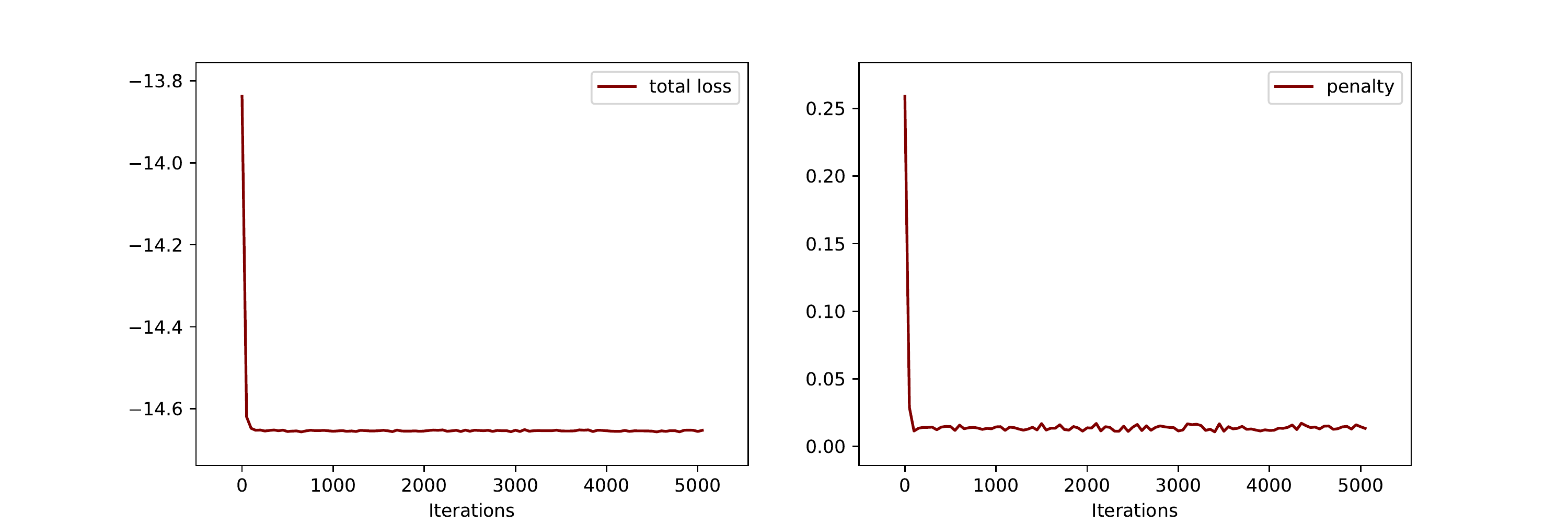}
     \end{subfigure}
     \hfill
     \begin{subfigure}{1\textwidth}
         \centering
         \includegraphics[width=\textwidth]{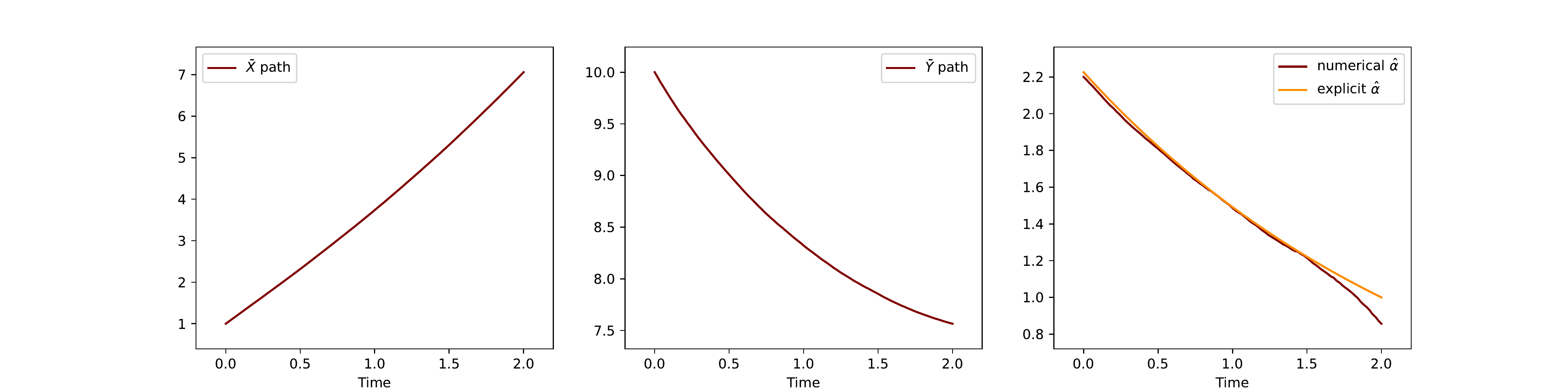}
     \end{subfigure}
        \caption{Results of the contract theory experiment
where we have interactions through variance of the states and where we used the generalized penalty method. Top: total loss including the penalty cost (left), penalty cost (right), bottom: the paths of $\bar X$ (left) and $\bar Y$ (middle), and the comparison of the explicit effort level to the mean of the effort levels found by using the generalized penalty method (right). We see that the penalty cost is converging to levels close to $0$. It is also seen that the numerical approach is
very good at mimicking the explicit results when the generalized numerical approach is used.}
        \label{fig:contract_exp1_penalty}
\end{figure}

Since our model fits in the special case introduced in Section~\ref{subsubsec:numerics_extended_special}, we test two numerical approaches, the generalized penalty method and the special numerical method introduced in Section~\ref{subsubsec:numerics_extended_special}, against the explicit results. In Figure~\ref{fig:contract_exp1_sannikov}, we can see the results of the numerical approach proposed for the special case and in Figure~\ref{fig:contract_exp1_penalty}, we can see the results of the generalized penalty method. In both numerical experiments we plot $\bar X$ and $\bar Y$ paths, total losses, as well as the comparison of the explicit optimal effort levels to the mean optimal effort levels found numerically. From the comparison of explicit effort levels to the numerical effort levels, we can see that both numerical approaches are good at mimicking the explicit results. In addition, in Figure~\ref{fig:contract_exp1_penalty} we also show that the penalty converges to levels close to $0$. We also stress that the special case algorithm is expected to work better as explained in Section~\ref{subsubsec:numerics_extended_special}.

\begin{table}[H]
	\centering
	{\small
		\begin{tabular}{@{}ccccccccc@{}}
			\toprule
			$T$ & $\Delta t$ & $\mu_0$ & $\sigma$ & $k$ & $a$ & $\beta_1$ & $\beta_2$& $\gamma$ \\
			\midrule
			$2.0$ & $0.02$ & $\delta_1$ & $1.0$ & $1.0$ & $0.4$ & $0.0$ & $0.0$ & $0.5$ \\
			\bottomrule
	\end{tabular}}
	\caption{\small Parameters used in experiment 1.}\label{table:var_inter_params}
\end{table}

\subsubsection{Experiment 2: interactions through the mean of the states.}
\label{subsubsec:experiments_contract_statemean}

Second, we consider a case where the interactions among agents are through the mean of the states. In other words, we take $\beta_2=\gamma=0$. In this case, we expect the optimal effort of an agent is given as

\begin{equation}
    \hat\alpha_t = \dfrac{e^{(a+\beta_1)(T-t)}}{k}.
\end{equation}

\begin{figure}[h]
     \centering
     \begin{subfigure}{1\textwidth}
         \centering
    \includegraphics[width=1\textwidth]{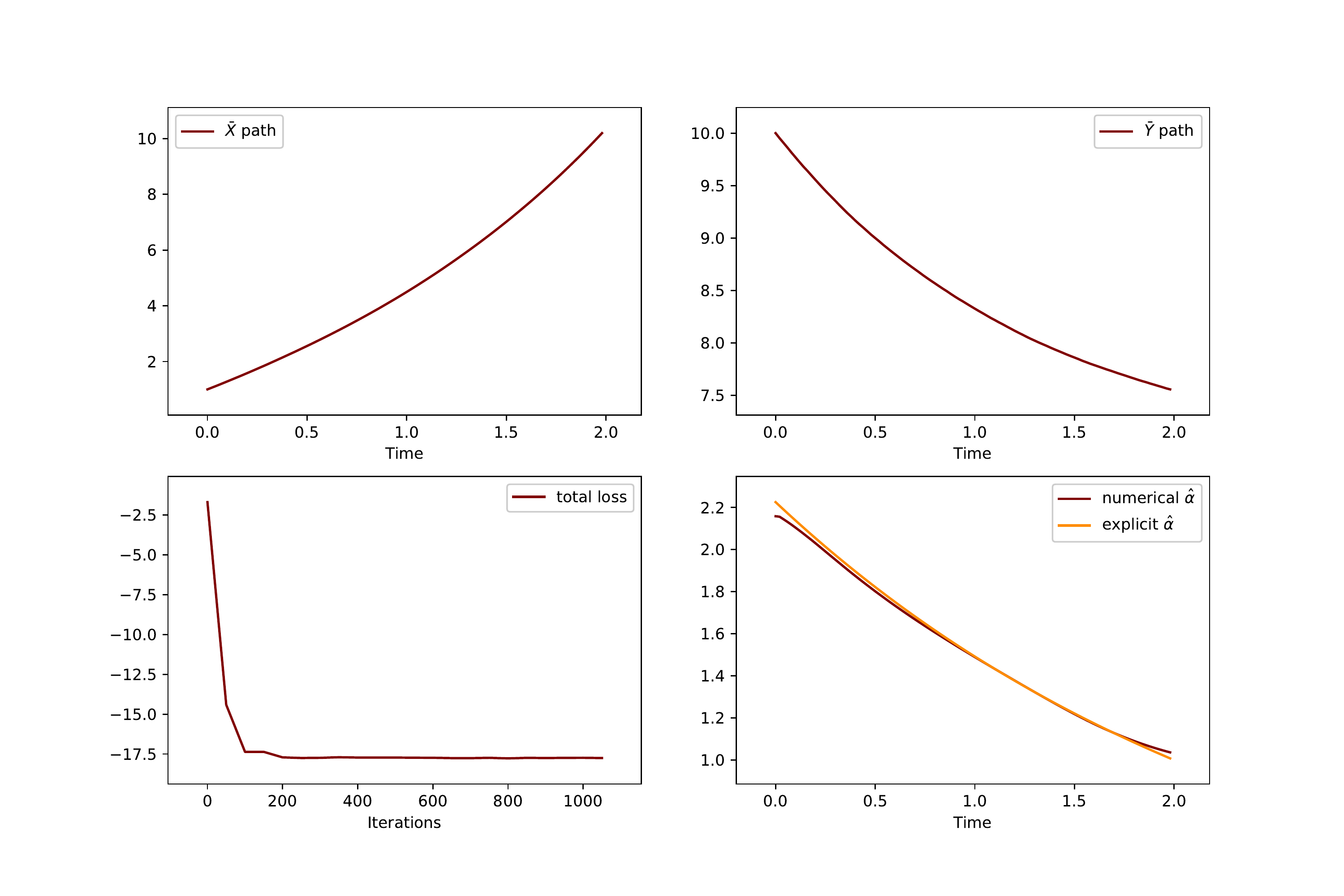}
    \end{subfigure}
    \caption{Results of the contract theory experiment where we have interactions through the mean of the state and where we use the algorithm proposed for the special case in Section~\ref{subsubsec:numerics_extended_special}. Top: the paths of $\bar X$ (left) and $\bar Y$ (right), bottom: total loss (left) and the comparison of the explicit effort level to the mean of the effort levels found by using the special numerical approach presented in Section~\ref{subsubsec:numerics_extended_special} (right). It is seen that the special numerical approach is very good at mimicking the explicit results.}
    \label{fig:contract_exp2_sannikov}
\end{figure}

\begin{figure}[h]
     \centering
     \begin{subfigure}{1\textwidth}
         \centering
         \includegraphics[width=\textwidth]{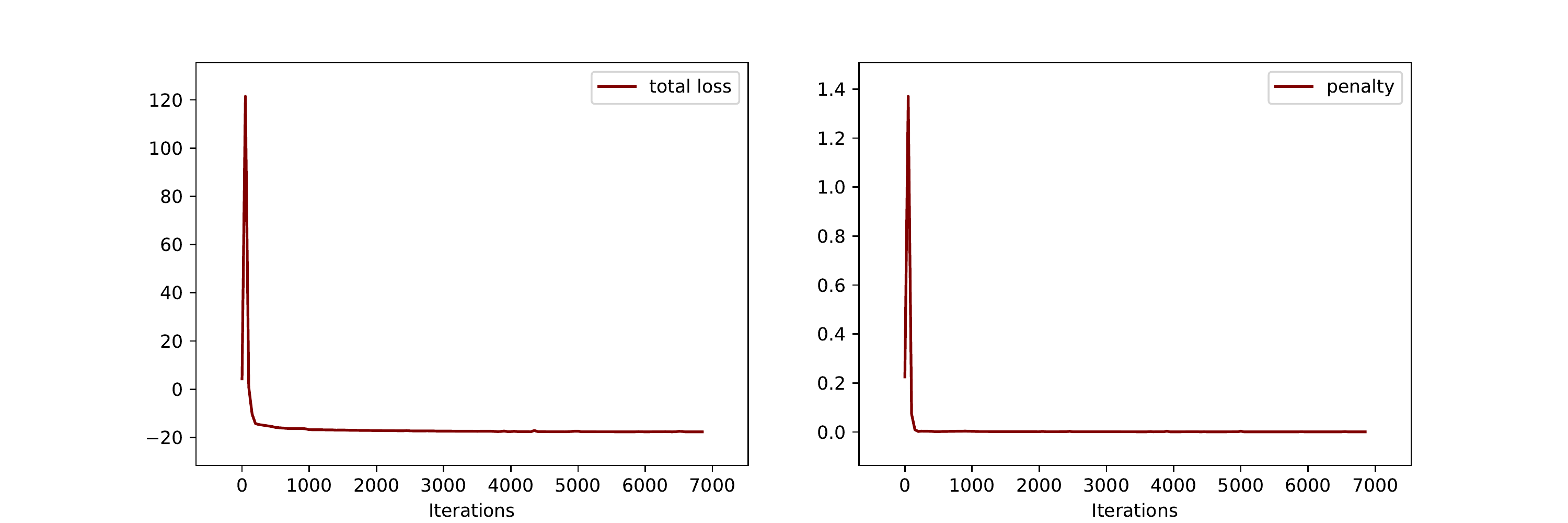}
     \end{subfigure}
     \hfill
     \begin{subfigure}{1\textwidth}
         \centering
         \includegraphics[width=\textwidth]{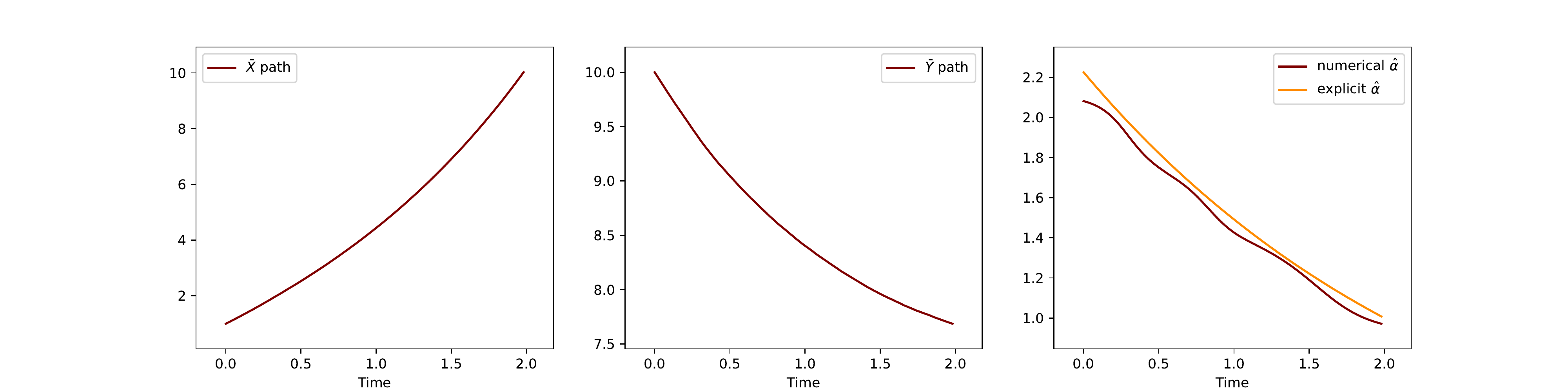}
     \end{subfigure}
        \caption{Results of the contract theory experiment
where we have interactions through the mean of the states and where we used the generalized penalty method where $\xi$ is implemented as a recurrent neural network. Top: total loss including the penalty cost (left), penalty cost (right), bottom: the paths of $\bar X$ (left) and $\bar Y$ (middle), and the comparison of the explicit effort level to the mean of the effort levels found by using the generalized penalty method (right). We see that the penalty cost is converging to levels close to $0$. It is also seen that the numerical approach is good at mimicking the explicit results when the generalized numerical approach is used.}
        \label{fig:contract_exp2_penalty}
\end{figure}

Again we test both numerical approaches (the one for the special case introduced in Section~\ref{subsubsec:numerics_extended_special} and the generalized penalty method) against the explicit solutions and the results can be seen in Figures~\ref{fig:contract_exp2_sannikov} and~\ref{fig:contract_exp2_penalty}, respectively. Similar to the results of Experiment 1, we see that both methods are good at mimicking the explicit results and the penalty cost converges to the levels around $0$. Different than Experiment 1 (and Experiment 3 introduced below), in the generalized penalty method function approximator for $\xi$ is implemented with a recurrent neural network to allow it to be implemented as a function of the whole path of $\bsX$. 

\begin{table}[H]
	\centering
	{\small
		\begin{tabular}{@{}ccccccccc@{}}
			\toprule
			$T$ & $\Delta t$ & $\mu_0$ & $\sigma$ & $k$ & $a$ & $\beta_1$ & $\beta_2$& $\gamma$ \\
			\midrule
			$2.0$ & $0.02$ & $\delta_1$ & $1.0$ & $1.0$ & $0.4$ & $0.25$ & $0.0$ & $0.0$ \\
			\bottomrule
	\end{tabular}}
	\caption{\small Parameters used in experiment 2.}
	\label{table:meanstate_inter_params}
\end{table}

\subsubsection{Experiment 3: interactions through the mean of the controls.}
\label{subsubsec:experiments_contract_controlmean}

In this section, in order to show the applicability of our numerical approach to \textit{extended} MFGs, we consider a case where the interactions among agents are through the mean of the controls. In other words, we take $\beta_1 = \gamma =0$. In this case, we expect the optimal effort of an agent is given as

\begin{equation}
    \hat\alpha_t = (1+\beta_2)\dfrac{e^{\alpha(T-t)}}{k}.
\end{equation}

\begin{table}[H]
	\centering
	{\small
		\begin{tabular}{@{}ccccccccc@{}}
			\toprule
			$T$ & $\Delta t$ & $\mu_0$ & $\sigma$ & $k$ & $a$ & $\beta_1$ & $\beta_2$& $\gamma$ \\
			\midrule
			$2.0$ & $0.02$ & $\delta_1$ & $1.0$ & $1.0$ & $0.4$ & $0.0$ & $0.5$ & $0.0$ \\
			\bottomrule
	\end{tabular}}
	\caption{\small  Parameters used in experiment 3.}
	\label{table:meancontrol_inter_params}
\end{table}

\begin{figure}[h]
    \centering
     \begin{subfigure}{1\textwidth}
         \centering
    \includegraphics[width=1\textwidth]{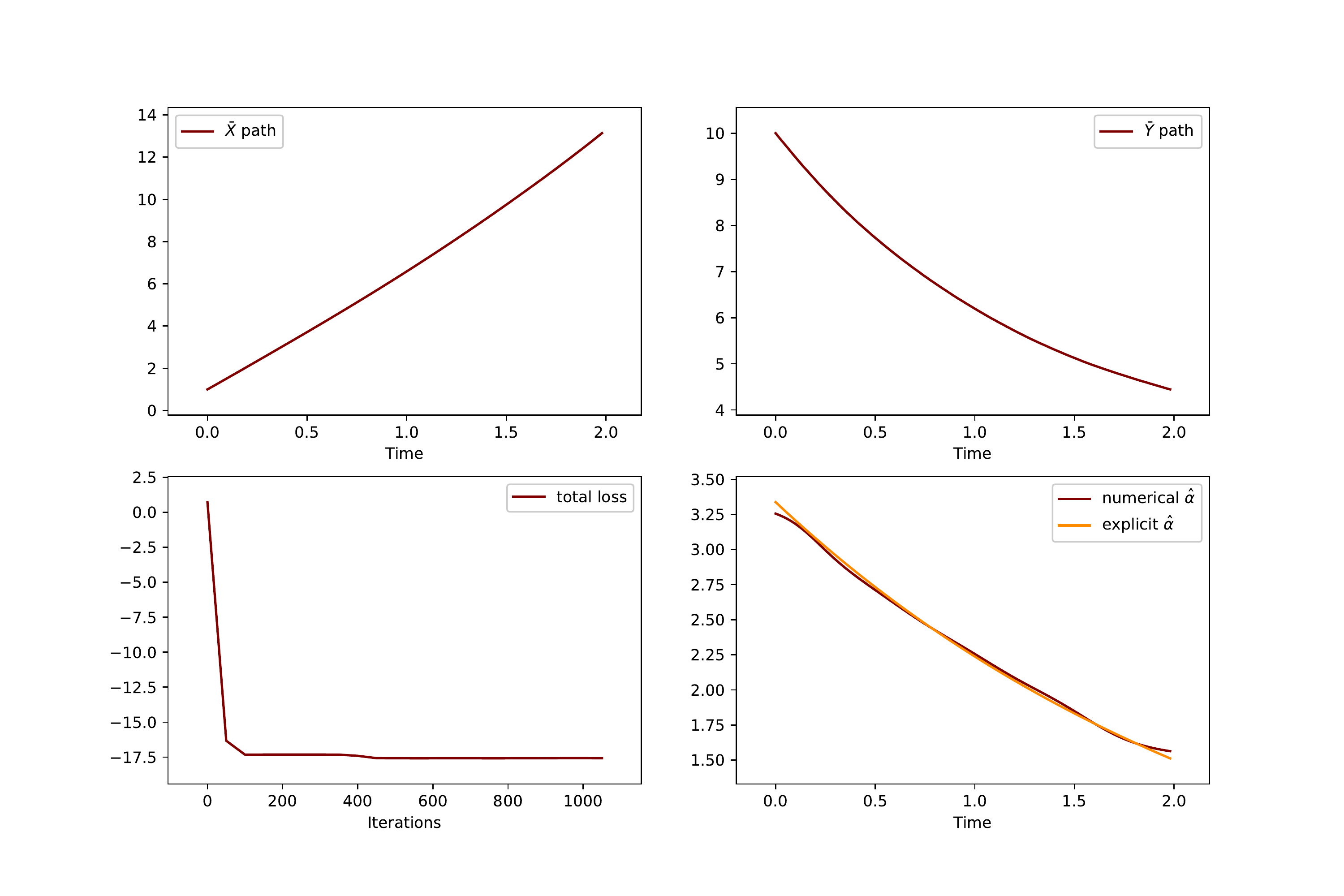}
    \end{subfigure}
    \caption{Results of the contract theory experiment where we have interactions through mean of the controls and where we use the algorithm proposed for the special case in Section~\ref{subsubsec:numerics_extended_special}. Top: the paths of $\bar X$ (left) and $\bar Y$ (right), bottom: total loss (left) and the comparison of the explicit effort level to the mean of the effort levels found by using the special numerical approach presented in Section~\ref{subsubsec:numerics_extended_special} (right). It is seen that the special numerical approach is very good at mimicking the explicit results.}
    \label{fig:contract_exp3_sannikov}
\end{figure}

\begin{figure}[h]
     \centering
     \begin{subfigure}{1\textwidth}
         \centering
         \includegraphics[width=\textwidth]{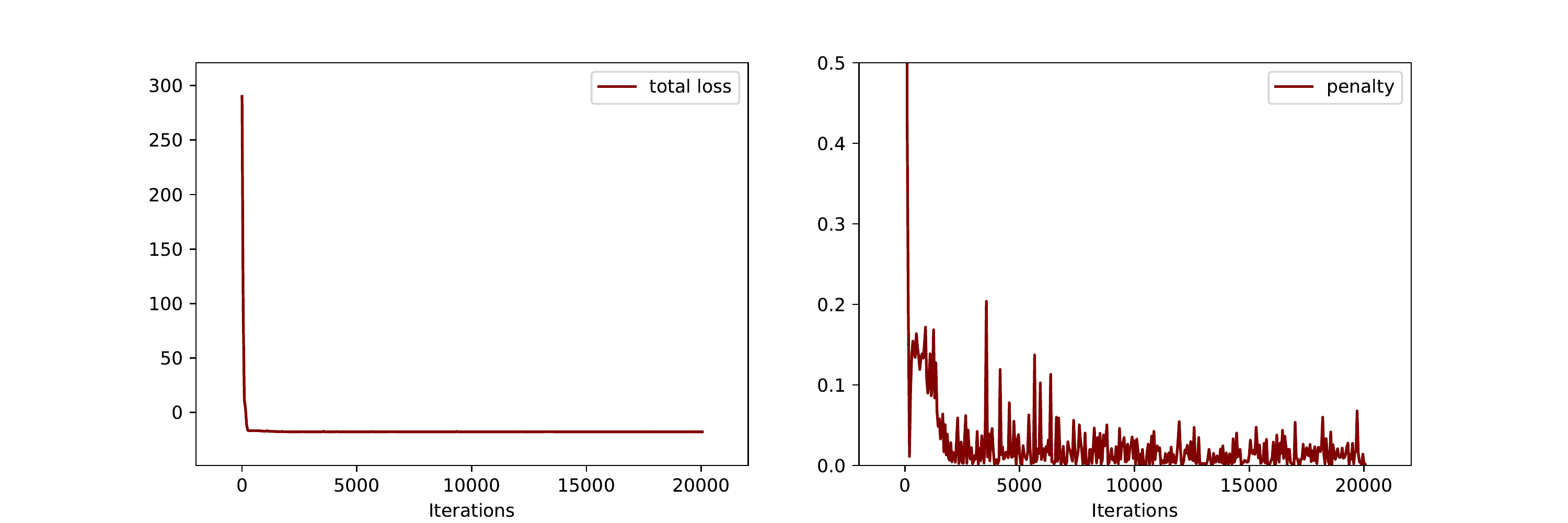}
     \end{subfigure}
     \hfill
     \begin{subfigure}{1\textwidth}
         \centering
         \includegraphics[width=\textwidth]{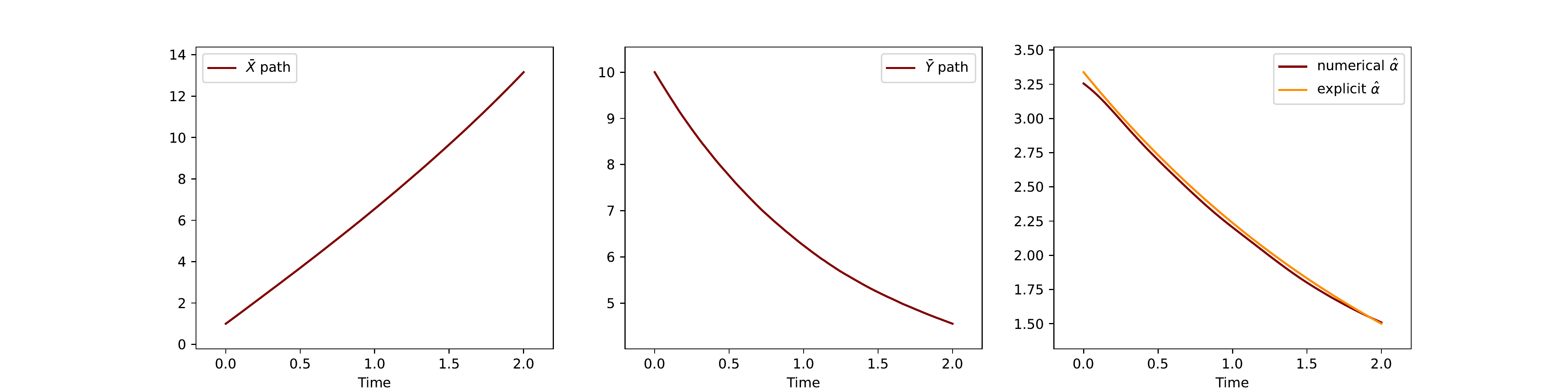}
     \end{subfigure}
        \caption{Results of the contract theory experiment
where we have interactions through the mean of the controls and where we used the generalized penalty method. Top: total loss including the penalty cost (left), penalty cost (right), bottom: the paths of $\bar X$ (left) and $\bar Y$ (middle), and the comparison of the explicit effort level to the mean of the effort levels found by using the generalized penalty method (right). We see that the penalty cost is converging to levels close to $0$. It is also seen that the numerical approach is
very good at mimicking the explicit results when the generalized numerical approach is used.}
        \label{fig:contract_exp3_penalty}
\end{figure}

As we did in the experiments 1 and 2, we test both numerical approaches (the one for the special case introduced in Section~\ref{subsubsec:numerics_extended_special} and the generalized penalty method) against the explicit solutions and the results can be seen in Figures~\ref{fig:contract_exp3_sannikov} and~\ref{fig:contract_exp3_penalty}, respectively. We see that in both algorithms, the numerical solutions are very good at mimicking the explicit solutions even when we have \textit{extended} mean field games setup.

\subsection{Systemic risk model under the influence of a principal.}
\label{subsec:experiments_systemicrisk}

In this section, we focus on the systemic risk model that is analyzed in ~\cite[Ch. 3]{carmona2018probabilistic} and extend it by adding a regulator (i.e., principal). This model is introduced for modeling the borrowing and lending between banks. The state of the agents i.e., the banks is the logarithm of their cash reserve and the state of the representative agent at time $t$ is denoted as $X_t$. The dynamics of the log-cash reserve is given as:
\begin{equation}
    dX_t = \big[a(\bar X_t - X_t)+\alpha_t\big]dt +  dW_t \,,
\end{equation}
where $W_t$ is the idiosyncratic noise and $a>0$. The control of the representative bank i.e., the lending and borrowing amounts at time $t$ is denoted with $\alpha_t$. The objective of the representative bank is given as 
\begin{equation}
    \inf_{\alpha} \ \mathbb{E}\Big[\int_0^T \big[\dfrac{\alpha_t^2}{2} - \lambda(t)\alpha_t(\bar X_t- X_t) + \dfrac{\epsilon}{2}(\bar X_t - X_t)^2\big]dt + \dfrac{c}{2} (\bar X_T - X_T)^2\Big] \,,
\end{equation}
where $\epsilon, c, \lambda >0$ are coefficients and $\bar X_t =\int_{\mathbb{R}} x d\mu_t (x)$. Here, $\epsilon$ and $c$ balance the individual bank's behavior with the average behavior of the other banks. $\lambda$ weighs the contribution of the components and helps to determine the sign of the control i.e., whether to borrow or to lend. Here, the model is extended by adding a regulator and the regulator's control is the $\boldsymbol\lambda$. Here, if the log cash reserve of the individual bank is smaller than the empirical mean, then the bank will want to borrow and choose $\alpha_t>0$, and vice versa. We should assume $\lambda(t)^2 \leq \epsilon,\ \forall t\in[0,T]$ in order to guarantee the convexity of the running cost functions while analyzing the explicit solution.

The regulator's aim is to keep the coefficient $\lambda(t)$ around the aimed level, $\lambda^{\rm aim}(t)$, $\forall t\in[0,T]$; furthermore, they want to minimize the proportion of banks defaulted. Mathematically the regulator's objective is written as follows:
\begin{equation}
    \inf_{\boldsymbol\lambda} \ \int_0^T (\lambda(t) - \lambda^{\rm aim}(t))^2 dt +\gamma \mathbb{E}\Big[\mathbbm{1}_{X_T<D}\Big],
\end{equation}
where $D<0$ is a constant that represents the critical liability threshold. Below this level the bank is considered in a state of \textit{default}.

Following the model analysis given in Section~\ref{subsec:model_mfganalysis}, the Hamiltonian for the representative bank can be written as follows
\begin{equation}
    H(t, X, \mu, z, \alpha; \lambda) = \left(a\left(\bar X- X\right) +\alpha\right)z +\dfrac{\alpha^2}{2} - \lambda \alpha \left(\bar X - X\right) + \dfrac{\epsilon}{2} \left(\bar X- X\right)^2,
\end{equation}
where $\bar X = \int_{\mathbb{R}} x d\mu (x)$.
Then, the minimizer of the Hamiltonian is given as
\begin{equation}
    \hat \alpha_t = \lambda(t) (\bar X_t - X_t) - Z_t.
\end{equation}
Finally, the system of FBSDE that characterize the mean field Nash equilibrium is given as
\begin{equation}
    \begin{aligned}
        dX_t &= \big[a(\bar X_t - X_t) + \lambda(t) (\bar X_t - X_t) - Z_t\big]dt + dW_t, && X_0 \sim\mu_0\\[3mm]
        dY_t &= -\Big[\dfrac{(\lambda(t) (\bar X_t - X_t)-Z_t)^2}{2}-\lambda(t) (\lambda(t) (\bar X_t - X_t) - Z_t) (\bar X_t - X_t) +&&\\&\hspace{6cm} \epsilon (\bar X_t - X_t)^2  \Big]dt + Z_t dW_t, \quad &&Y_T = \dfrac{c}{2} (\bar X_T - X_T)^2 \,.
    \end{aligned}
\end{equation}

\begin{table}[H]
	\centering
	{\small
		\begin{tabular}{@{}cccccccc@{}}
			\toprule
			$T$ & $\Delta t$ & $\mu_0$ & $a$ & $c$ & $\epsilon$ & $\lambda^{\rm aim}$ &  $D$\\
			\midrule
			$2.0$ & $0.02$ & $\delta_1$ & $1.0$ & $1.0$ & $1.0$ & $0.5$ &  $-0.001$ \\
			\bottomrule
	\end{tabular}}
	\caption{\small Parameters used in the systemic risk experiment. }
	\label{table:systemicrisk}
\end{table}

\begin{figure}[h]
    \centering
    \includegraphics[width=1\textwidth]{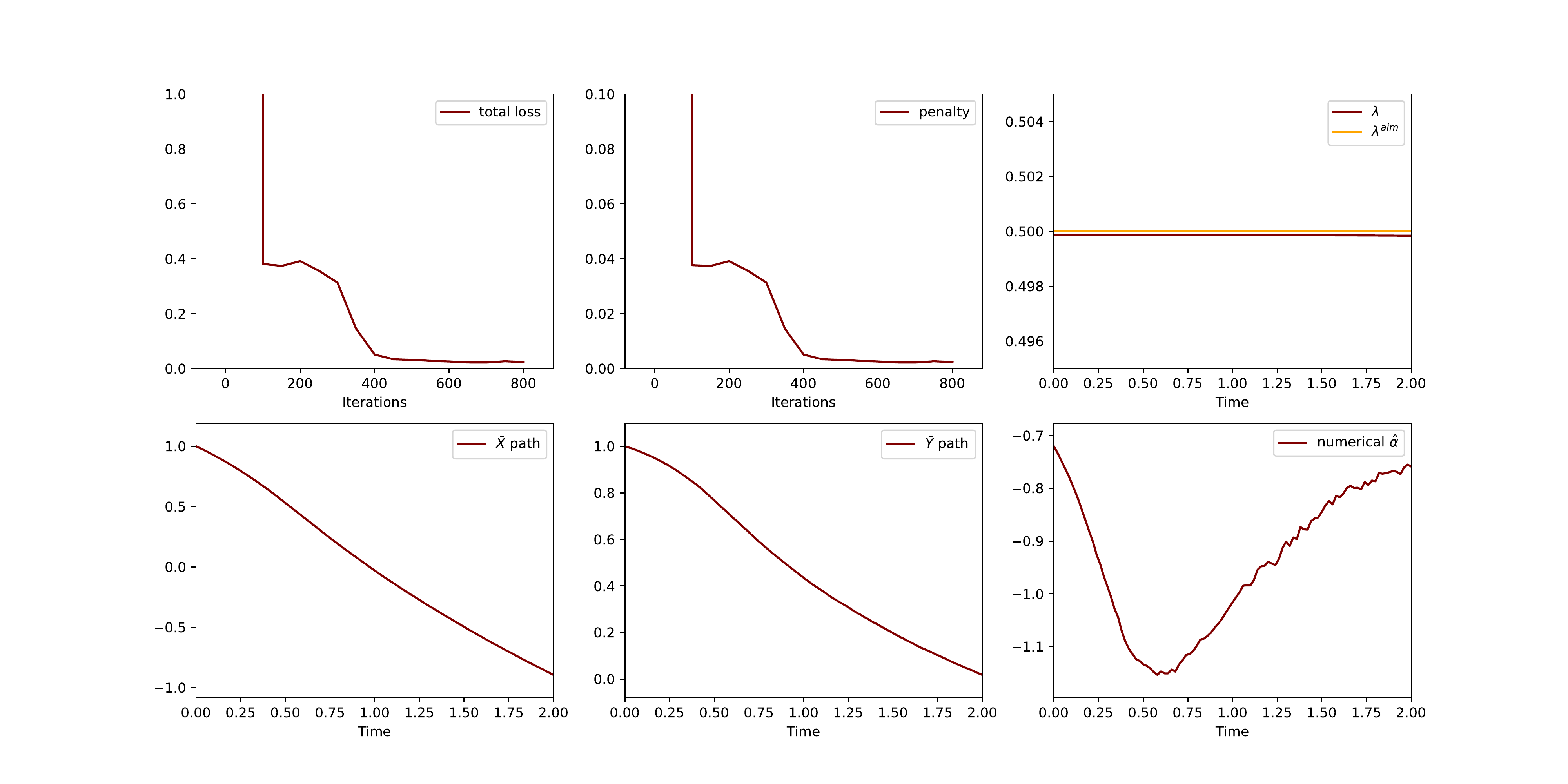}
    \caption{Results of the systemic risk experiment with the addition of a regulator that controls $\boldsymbol\lambda$. In this experiment, the regulator does not have a cost for the proportion of the banks defaulted, i.e., $\gamma=0$.) Top: total loss including the penalty cost (left), penalty cost (middle), comparison of the regulator's aimed policy level and the numerical optimal policy level (right), Bottom: the paths of $\bar X$ (left) and $\bar Y$ (middle), the mean of the control of the banks found by using the generalized penalty method (right).}
    \label{fig:systemicrisk_nodefault}
\end{figure}

\begin{figure}[h]
    \centering
    \includegraphics[width=1\textwidth]{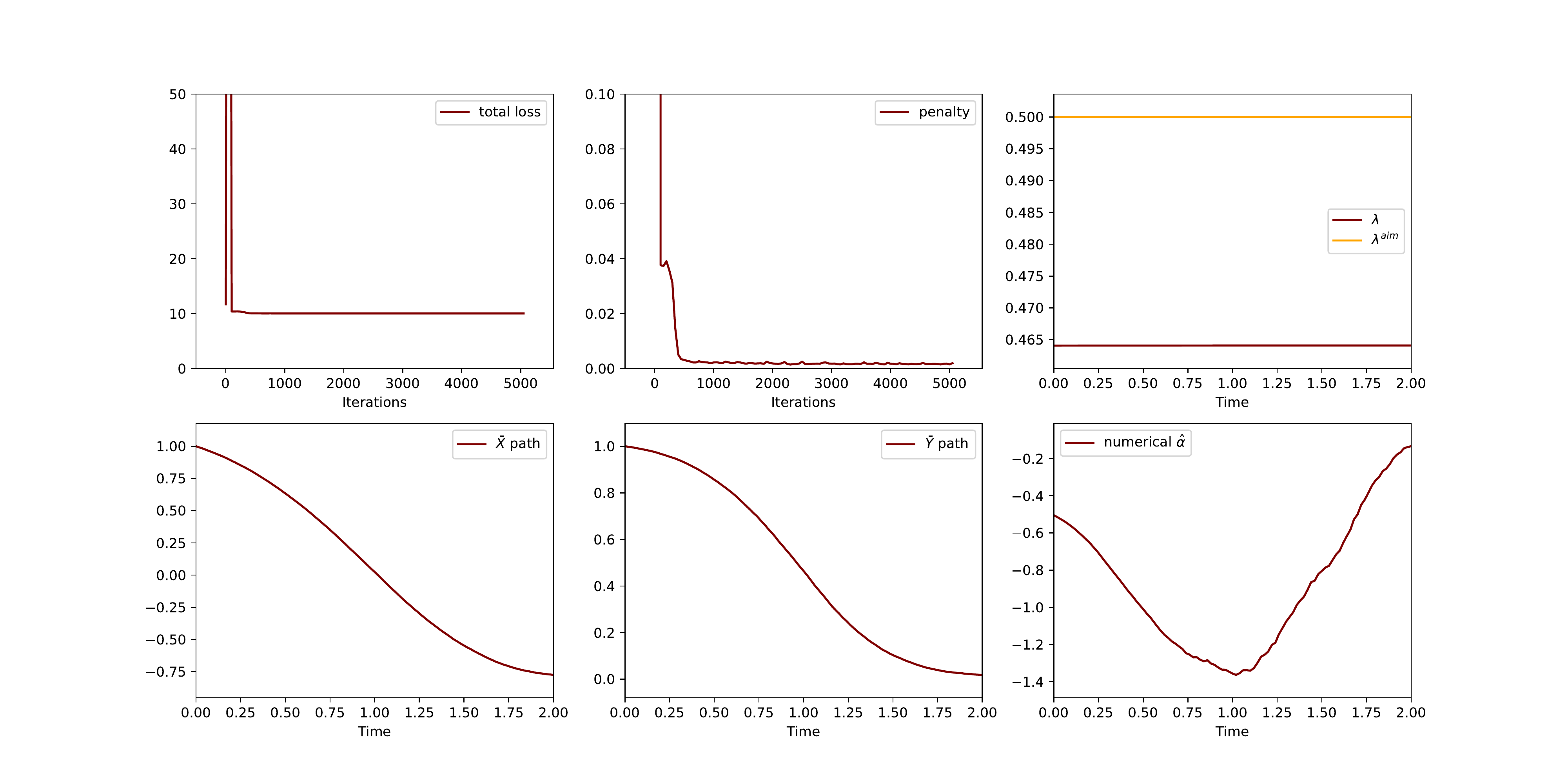}
    \caption{Results of the systemic risk experiment with the addition of a regulator that controls $\boldsymbol\lambda$. In this experiment, the regulator has a cost for the proportion of the banks defaulted. ($\gamma=10$) Top: total loss including the penalty cost (left), penalty cost (middle), comparison of the regulator's aimed policy level and the numerical optimal policy level (right), Bottom: the paths of $\bar X$ (left) and $\bar Y$ (middle), the mean of the control of the banks found by using the generalized penalty method (right).}
    \label{fig:systemicrisk_withdefault}
\end{figure}

The model in this example does not fit in the special case form discussed in~\ref{subsubsec:numerics_extended_special}; therefore, we will be only presenting our results where we used the generalized penalty method. We first consider an example where the regulator does not have a cost for defaulting banks (i.e., $\gamma=0$). Here, we expect the regulator to choose $\hat \lambda(t) = \lambda^{\rm aim}(t)$. As it can be seen from Figure~\ref{fig:systemicrisk_nodefault}, our numerical approach is able to capture this phenomenon. Later, we also present results from an experiment where we include the cost related to the defaulted banks for the regulator (i.e., $\gamma\neq 0$). From Figure~\ref{fig:systemicrisk_withdefault}, we can see that the optimal policy of the regulator is different than the $\lambda^{\rm aim}(t)$ in order to decrease the proportion of banks defaulting.

\begin{figure}[h]
    \centering
    \begin{subfigure}[b]{\textwidth}
         \centering
         \includegraphics[width=\textwidth]{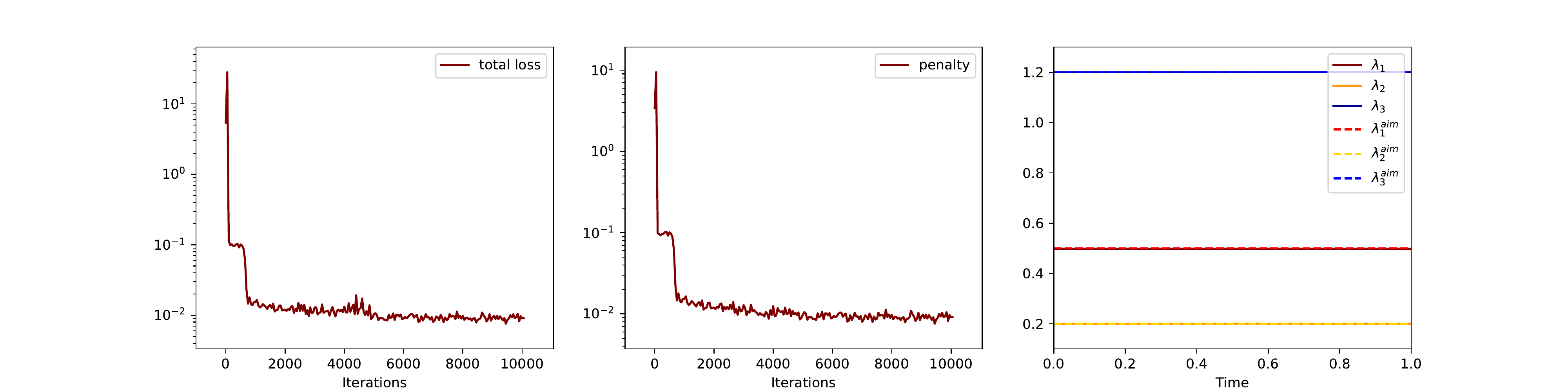}
     \end{subfigure}
     \hfill
     \begin{subfigure}[b]{\textwidth}
         \centering
         \includegraphics[width=\textwidth]{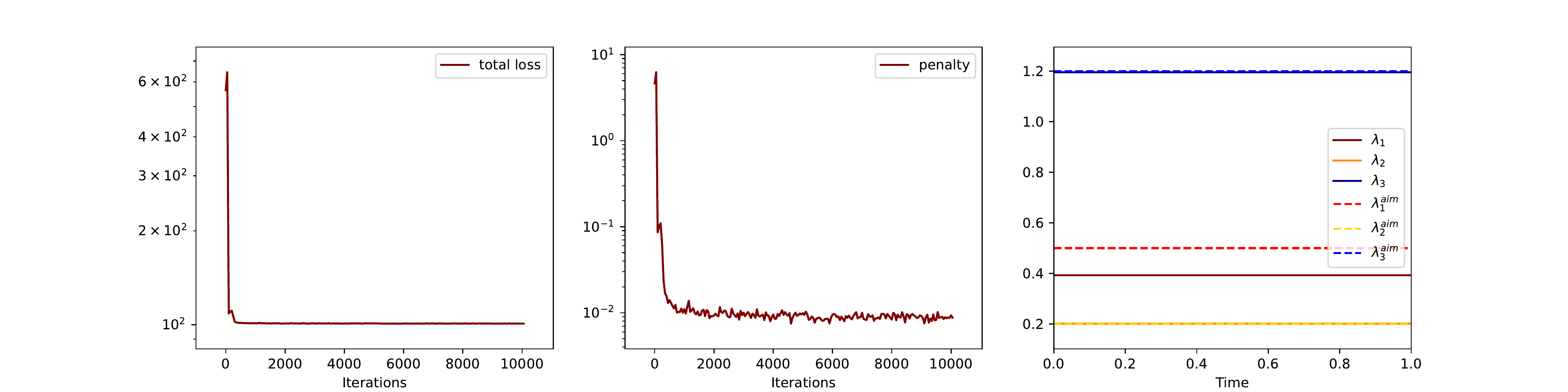}
     \end{subfigure}
    \caption{Results of the systemic risk experiment where $m=3$ and with the addition of a regulator that controls $\boldsymbol\lambda$. Top: total loss including the penalty cost (left), penalty cost (middle), comparison of the regulator's aimed policy level and the numerical optimal policy level (right) in the experiment where $\gamma = 0$, Bottom: total loss including the penalty cost (left), penalty cost (middle), comparison of the regulator's aimed policy level and the numerical optimal policy level (right) in the experiment where $\gamma = 50$.}
    \label{fig:systemicrisk_3D}
\end{figure}

\begin{remark}
    In order to show the scalability of the proposed numerical approach in higher dimensional problems, we also consider an extended version of the systemic risk model. We assume that banks have $m$-dimensional log cash reserve as their state which can represent the cash reserves held as different currencies. The control of the representative bank is also designed $m$-dimensional to represent the borrowing and lending from the central bank in the considered currencies and denoted by $\boldsymbol{\alpha}=(\alpha_t)_{t\in[0,T]}$ where $\alpha_t = (\alpha_t^1, \alpha_t^2, \dots, \alpha_t^m)$. Then the representative bank would like to minimize the following cost
    \begin{equation*}
        \inf_{\boldsymbol{\alpha}} \ \mathbb{E}\Big[\int_0^T \sum_{l=1}^m\Big[\dfrac{(\alpha_t^l)^2}{2} - \lambda_l(t)\alpha^l_t(\bar X^l_t- X^l_t) + \dfrac{\epsilon_l}{2}(\bar X^l_t - X^l_t)^2\Big]dt + \sum_{l=1}^m\dfrac{c_l}{2} (\bar X^l_T - X^l_T)^2\Big] \,.
    \end{equation*}
    The state of the representative bank, $\bsX = (X_t)_{t\in[0,T]}$ where $X_t = (X_t^1, X_t^2, \dots, X_t^m)$, has the following dynamics
    \begin{equation*}
        dX^l_t = \big[a_l(\bar X^l_t - X^l_t)+\alpha^l_t\big]dt +  dW^l_t\, ,\qquad \forall l \in \{1,2,\dots,m\}.
    \end{equation*}
    Finally the regulator would like to minimize the following cost
    \begin{equation*}
    \inf_{\boldsymbol\lambda} \ \int_0^T \sum_{l=1}^m (\lambda_l(t) - \lambda^{aim}_l(t))^2 dt +\gamma \mathbb{E}\Big[\mathbbm{1}_{\sum_{l=1}^m c^0_l X^l_T<D}\Big],
    \end{equation*}
    where $c_l^0$ for $l \in \{1, 2, \dots, m\}$ are nonnegative coefficients. In Figure~\ref{fig:systemicrisk_3D}, we see the results of two experiments where $m$ is chosen to be 3. The top plots correspond to the case where $\gamma$ is chosen to be 0 and bottom plots corresponding to the case where $\gamma$ is chosen to be 50. As expected from the form of the model, in the case where $\gamma=0$, we can see that the optimal policy of the regulator is given as $\hat{\lambda}_l(t) = \lambda^{\rm aim}_l (t)$. Furthermore, we observe that when $\gamma$ is nonzero, the optimal policy of the regulator is differing from the aimed levels to decrease the proportion of the banks defaulting similar to the results we observed in the one-dimensional case. We realize that now the principal is choosing $\hat{\lambda}_2(t)<\lambda^{\rm aim}_2(t)$, this lets the representative player to be able to borrow the second currency type at higher levels which helps with lower levels of defaulting. Our numerical approach was capable of capturing the higher dimensional results without any significant change in the computational time and the implementation was directly executed by increasing the output dimensions of the neural networks approximating the principal's policy, $\lambda(t)$, and $Z(t, X_t)$.
\end{remark}

\section{Conclusion.}
\label{sec:conclusion}
In this paper, we first write the Stakelberg mean field game problem as a constrained stochastic optimal problem and use a penalty method to add the Nash equilibrium constraint in the problem of the principal. Then, we propose a numerical approach to solve Stackelberg MFG. In such problems, the principal aims to optimize her policy while taking into account the reaction of the strategic agents to this policy. This creates a bi-level problem where, at the population level, the agents solve a mean field Nash equilibrium given the policy and, at the upper level, the principal looks for the optimal policy. We first rewrite this problem as a single-level problem through writing it as a first constrained and second as a penalized stochastic optimal control problem to solve it efficiently and introduce a numerical approach that uses neural networks and Monte Carlo simulation. 
On the theoretical side, we first show the convergence of the reformulated, penalized, problem to the original problem. Second, we analyze the errors due to the approximations used in our numerical implementation: empirical population approximation, time discretization, and neural network approximation. For each, we obtain an upper bound on the rate of convergence. 
We also show that the numerical approach can be extended to more complex model forms where the principal has a contract theory type of problem with an additional $\mathcal{F}_T$-measurable incentive.
Lastly, we show the effectiveness of the proposed numerical approach with two different examples.

From here, several directions could be investigated. 
Theoretically one of the most challenging direction is the proof of existence, uniqueness, and regularity of the solution of Stackelberg MFG, in settings beyond the ones studied here (for instance in the presence common noise). This is also relevant from the numerical viewpoint because the regularity of the principal's is important to know whether it will be possible to approximate it using a neural network. It may inform the choice of the architecture (and in particular the activation function).
Another interesting direction is about selecting a specific Nash equilibrium when there is no uniqueness. The principal might be able to select the best equilibrium from their viewpoint, for example to reduce the price of anarchy. However, this is a broader question, since to the best of our knowledge, most of the MFG literature focuses on situations in which the Nash equilibrium is unique. 
Last but not least, from the numerical viewpoint, it would be interesting to solve more realistic Stackelberg MFGs. Indeed, in this work we focused on showcasing the performance of the proposed method on well-known benchmarks from the literature. But deep learning offers the hope to solve much more complex models, and we this method will open new possibilities in terms of solving complex Stackelberg MFGs.

\subsubsection*{Acknowledgments.}
The authors would like to thank Ren\'e Carmona for helpful discussions. The authors also would like to thank the anonymous reviewers for their insightful comments and suggestions which helped improving the paper.

\appendix

\section{Auxiliary Proofs}

\subsection{Penalization analysis} 
\label{sec:proofs-lemma-penalization}

We provide auxiliary lemmas that are used in the proof of Theorem~\ref{thm:main-cv-thm-penalization}.

\begin{lemma}
\label{lem:eps-min-ineq}
Let $\bsbeta \in\UU$. 
Let $\bsbeta' \in\VV$ be an $\epsilon$-minimizer for the penalized problem~\eqref{eq:penalized-pb-J0P} with penalty parameter $\nu>0$.
Then, $J^0(\bsbeta) +\epsilon\geq J^{0}_{\nu}(\bsbeta')\geq J^0(\bsbeta')$.
\end{lemma}
\proof[Proof of Lemma~\ref{lem:eps-min-ineq}. ]
Let $\bsbeta$ be as in the statement. In particular, it is in the feasible set of problem~\eqref{eq:original-pb-J0} so the penalty function is equal to 0. We have:
\begin{align*}
        J^0(\bsbeta) + \epsilon
        &= J^0(\bsbeta) +\nu \overline{\mbsP}(\bsbeta)  + \epsilon
        = J^0_\nu(\bsbeta)  + \epsilon
        \ge J^0_{\nu}(\bsbeta')
        \geq J^0(\bsbeta').
\end{align*}
\endproof

\begin{lemma}
\label{lem:cont-J-P}
    Suppose Assumption~\ref{hyp2} holds. Then, $J^0$ and $\overline{\mbsP}$ are continuous on $\VV$.
\end{lemma} 
\proof[Proof of Lemma~\ref{lem:cont-J-P}. ]
\noindent \textbf{Step 1:} We show that the state process $\bsX^{\bsbeta}$ is continuous with respect to control $\bsbeta$. Let $\tilde\bsbeta, \bsbeta \in \VV$. We bound $\EE|\bsX^{\tilde{\bsbeta}}- \bsX^{{\bsbeta}}|^2$ in terms of $\|\tilde\bsbeta -  \bsbeta\|_{\VV}^2$, where we recall that the norm on $\VV$ is the sup norm. 
\vskip6pt
Let $\Gamma = \|\tilde\bsbeta -  \bsbeta\|_{\VV}$. 
Below, $C$ denotes a generic constant whose value can change from one line to the other but is independent of $\Gamma$, $L_z$, and $L_\lambda$. 
By Lipschitz continuity of $\tilde{z}$ and $\tilde{\lambda}$ with respect to $x$, we have:
\begin{equation}\label{eq:cont_bound_betas_gamma_new2}
\begin{aligned}
&\EE\Big[|\tilde\beta(t, X_0, X^{\tilde{\bsbeta}}_t) - \beta(t,X_0, X^{\bsbeta}_t)|^2\Big]
    \\[2mm]
    &\leq \EE\Big[|\tilde{\beta}(t, X_0, X_t^{\tilde{\bsbeta}}) - \tilde{\beta}(t, X_0, X_t^{\bsbeta})|^2 +|\tilde{\beta}(t, X_0, X_t^{\bsbeta}) - \beta(t, X_0, X_t^{\bsbeta})|^2\Big]
    \\[2mm]
    &\leq {(}L^{{2}}_z{+L^2_\lambda)}\EE\Big[|X_t^{\tilde{\bsbeta}} - X_t^{\bsbeta}|^2\Big]+ \Gamma^2
    \end{aligned}
\end{equation}

To simplify the notation, introduce $\varepsilon_t^{\check\bsbeta}:= (t, X_t^{\check\bsbeta}, \beta(t, X_0, X^{\check\bsbeta}_t), \mu_t^{\check\bsbeta})$ where $\check\bsbeta \in \{\boldsymbol{\tilde\beta},\bsbeta\}$. Recall that we introduced the notations $\mathcal{b}$ and $\mathcal{f}$ defined in~\eqref{eq:def-calb-calf}. With these notations we have the SDEs:
\begin{equation}
\label{eq:SDE_newnotation_new}
    \begin{aligned}
        dX^{\bsbeta}_t &= \mathcal{b}(\varepsilon_t^{\bsbeta})dt +\sigma dW_t, \qquad &&X_0 \sim \mu_0,\\[2mm]
        dY^{\bsbeta}_t &= -\mathcal{f}(\varepsilon_t^{\bsbeta})dt + z(t,X^{\bsbeta}_t) \cdot dW_t, && Y_0 = y_0(X_0).
    \end{aligned}
\end{equation}
Using these dynamics, for all $r\in[0,T]$ we have
\begin{equation}
\label{eq:cost_cont_X_bound_new}
    \begin{aligned}
        &\mathbb{E} \Big[\sup_{s\in[0,r]} |X_s^{\tilde{\bsbeta}}-X_s^{\bsbeta}|^2 \Big]\leq C\Big( \mathbb{E}\sup_{s \in [0,r]} \big|\int_0^s [\mathcal{b}(\varepsilon_t^{\tilde{\bsbeta}})-\mathcal{b}(\varepsilon_t^{\bsbeta})]dt\big|^2 \Big).
    \end{aligned}
\end{equation}

     By Cauchy-Schwarz inequality and the Lipschitz property of the drift $\mathcal{b}$ (Assumption~\ref{hyp2:drift-lip}) we obtain that
\begin{align}
         & \mathbb{E}\sup_{s \in [0,r]}\big|\int_0^s [\mathcal{b}(\varepsilon_t^{\tilde{\bsbeta}})-\mathcal{b}(\varepsilon_t^{\bsbeta})]dt\big|^2& 
         \notag 
         \\[1mm]
         &\hskip1cm\leq \mathbb{E}\sup_{s \in [0,r]} s \int_0^s \big[\mathcal{b}(\varepsilon_t^{\tilde{\bsbeta}})-\mathcal{b}(\varepsilon_t^{\bsbeta})\big]^2dt
         \notag
         \\[1mm]
         &\hskip1cm\leq r \int_0^r\mathbb{E} \big[\mathcal{b}(\varepsilon_t^{\tilde{\bsbeta}})-\mathcal{b}(\varepsilon_t^{\bsbeta})\big]^2dt
         \notag 
         \\[1mm]
         &\hskip1cm\leq C \int_0^r \Big(\mathbb{E}\big|X_t^{\tilde{\bsbeta}} - X_t^{\bsbeta}\big|^2  + \mathbb{E}\big|\tilde\beta(t, X_0, X^{\tilde{\bsbeta}}_t) - \beta(t, X_0, X^{\bsbeta}_t)\big|^2+ W^2_2(\mu^{\tilde{\bsbeta}}_t, \mu^{\bsbeta}_t)\Big) dt
         \notag 
         \\[1mm]
         &\hskip1cm\leq C(1+L^{{2}}_z{+L^2_\lambda}) \int_0^r \mathbb{E} \sup_{0\leq s\leq t} \big|X_s^{\tilde{\bsbeta}} - X_s^{\bsbeta}\big|^2 dt + C \Gamma^2 + C\int_0^r W^2_2(\mu^{\tilde{\bsbeta}}_t, \mu^{\bsbeta}_t) dt,
    \label{eq:b_cont_X_bound_new}
\end{align}
where the last inequality holds by~\eqref{eq:cont_bound_betas_gamma_new2}. %

Combining~\eqref{eq:cost_cont_X_bound_new} and~\eqref{eq:b_cont_X_bound_new}, we obtain:
\begin{equation*}
    \begin{aligned}
        &\mathbb{E} \Big[\sup_{s\in[0,r]} |X_s^{\tilde{\bsbeta}}-X_s^{\bsbeta}|^2  \Big]
        \\
        & \hskip0.5cm \leq C\Big(  (1+L^{{2}}_z{+L^2_\lambda})\int_0^r \mathbb{E}\sup_{s\in[0,r]}\big|X_s^{\tilde{\bsbeta}} - X_s^{\bsbeta}\big|^2 dt + \Gamma^2 + \int_0^r W^2_2(\mu^{\tilde{\bsbeta}}_t, \mu^{\bsbeta}_t) dt  \Big)
        \\[1mm]
        & \hskip0.5cm \leq C\Big(  (1+L^{{2}}_z{+L^2_\lambda})  \int_0^r \mathbb{E}\sup_{s\in[0,r]}\big|X_s^{\tilde{\bsbeta}} - X_s^{\bsbeta}\big|^2 dt + \Gamma^2+ \int_0^r  \mathbb{E}\big|X_t^{\tilde{\bsbeta}} - X_t^{\bsbeta}\big|^2 dt \Big)
        \\[1mm]
        & \hskip0.5cm \leq C\Big(  (1+L^{{2}}_z{+L^2_\lambda}) \int_0^r \mathbb{E}\sup_{s\in[0,r]}\big|X_s^{\tilde{\bsbeta}} - X_s^{\bsbeta}\big|^2 dt  + \Gamma^2 \Big).
    \end{aligned}
\end{equation*}
By Gr\"onwall's inequality, we conclude that
\begin{equation}
\label{eq:bound-supp-diff-X-Gamma2_new}
    \mathbb{E} \Big[\sup_{s\in[0,T]} |X_s^{\tilde{\bsbeta}}-X_s^{\bsbeta}|^2  \Big] \leq C e^{C (1+L^{{2}}_z{+L^2_\lambda})}\Gamma^2,
\end{equation}
which implies that $\bsX^{\bsbeta}$ is continuous with respect to control $\bsbeta$.

\vskip12pt
\noindent \textbf{Step 2 (Continuity of $J^0$):} We now bound $|J^0(\tilde{\bsbeta}) - J^0(\bsbeta)|$ in terms of $\Gamma$.\\
By the Lipschitzness assumption on $f_0$ and $g_0$ (Assumption~\ref{hyp2:cost0-lip}), we infer that:
\begin{equation*}
    \begin{aligned}
        &\Big|J^0(\tilde{\bsbeta}) - J^0(\bsbeta) \Big| 
        \\[1mm]
        &\hskip0.5cm \leq {\EE}\int_0^T\Big|f_0(t,\mu_t^{\tilde{\bsbeta}}, \tilde{\lambda}(t{,X_t^{\tilde{\bsbeta}}}))-f_0(t,\mu_t^{\bsbeta}, \lambda(t{,X_t^{{\bsbeta}}}))\Big|dt + \Big|g_0(\mu_T^{\tilde{\bsbeta}}, \tilde{\lambda}(t{,X_t^{\tilde{\bsbeta}}})) -g_0(\mu_T^{\bsbeta}, \lambda(T{,X_t^{{\bsbeta}}}))\Big| 
        \\[1mm]
        &\hskip0.5cm\leq  C \left(\int_0^T \Big|W_2(\mu_t^{\tilde{\bsbeta}}, \mu_t^{\bsbeta})\Big| dt +{\EE} \int_0^T \Big|\tilde{\lambda}(T{,X_t^{\tilde{\bsbeta}}}) - \lambda(T{,X_t^{{\bsbeta}}})\Big|dt  + \Big|W_2(\mu_T^{\tilde{\bsbeta}}, \mu_T^{\bsbeta})\Big| + {\EE}\Big|\tilde{\lambda}(T{,X_t^{\tilde{\bsbeta}}}) - {\lambda}(T{,X_t^{{\bsbeta}}})\Big| \right)
        \\[1mm]
         &\hskip0.5cm\leq C \Big( \int_0^T \Big( \mathbb{E}\Big[\sup_{0\leq s \leq t}|X_s^{\tilde{\bsbeta}} - X_s^{\bsbeta}|^2\Big]\Big)^{1/2}dt + \Gamma{+T \mathbb{E}  \Big[\sup_{0\leq s \leq t}|X_s^{\tilde{\bsbeta}} - X_s^{\bsbeta}|\Big]\Big)}
         \\[1mm]
          &\hskip0.5cm\leq C (1 +  e^{C (1+L^{{2}}_z{+L^2_\lambda})})\Gamma,
    \end{aligned}
\end{equation*}
where the second to last inequality is by the fact that $\Gamma = \|\tilde\bsbeta -  \bsbeta\|_{\VV}$, and the last inequality is by~\eqref{eq:bound-supp-diff-X-Gamma2_new}.
\vskip12pt
\noindent\textbf{Step 3 (Continuity of $\overline{\mbsP}$):} In order to bound $|\overline{\mbsP}(\tilde{\bsbeta}) - \overline{\mbsP}(\bsbeta)|$ in terms of $\Gamma$, we need to show that $Y_T^{\bsbeta}$, $X_T^{\bsbeta}$ and $\mu_T^{\bsbeta} =\mathcal{L}(X_T^{\bsbeta})$ are continuous with respect to $\bsbeta$. We start with $Y_T^{\bsbeta}$. Remember that our SDE system with the new notations is given in \eqref{eq:SDE_newnotation_new}. We follow the ideas used in {\bf Step 1} of this proof: 
\begin{equation}
\label{eq:cost_cont_YT_bound_new_2}
    \begin{aligned}
        &\mathbb{E} \Big[ |Y_T^{\tilde{\bsbeta}}-Y_T^{\bsbeta}|^2 \Big]\\
        &  \leq C \Big( \mathbb{E}\big|\tilde{y}_0(X_0)- y_0(X_0)\big|^2 +  \mathbb{E} \big|\int_0^T [\mathcal{f}(\varepsilon_t^{\tilde{\bsbeta}})-\mathcal{f}(\varepsilon_t^{\bsbeta})]dt\big|^2 +\mathbb{E}  \big|\int_0^s [\tilde{z}_t(t,X_t^{\tilde\bsbeta})- z_t(t,X_t^{\bsbeta})]dW_t\big|^2\Big)
    \end{aligned}
\end{equation} 
In order to bound the second term on the right hand side, we use the Lipschitz property of the function $\mathcal{f}$ (Assumption~\ref{hyp2:f-lip}), and the bound given in \eqref{eq:cont_bound_betas_gamma_new2}:
\begin{align}
    &\mathbb{E} \big|\int_0^T [\mathcal{f}(\varepsilon_t^{\tilde{\bsbeta}})-\mathcal{f}(\varepsilon_t^{\bsbeta})]dt\big|^2
    \notag \\[1mm]
    &\hskip1cm \leq T\int_0^T \mathbb{E}\big[\mathcal{f}(\varepsilon_t^{\tilde{\bsbeta}})-f(\varepsilon_t^{\bsbeta})\big]^2dt
    \notag \\[1mm]
    &\hskip1cm \leq CT\int_0^T \Big(\mathbb{E}\big|X_t^{\tilde{\bsbeta}} - X_t^{\bsbeta}\big|^2  + \mathbb{E}\big|\tilde\beta(t, X_0, X^{\tilde{\bsbeta}}_t) - \beta(t, X_0, X^{\bsbeta}_t)\big|^2+ W^2_2(\mu^{\tilde{\bsbeta}}_t, \mu^{\bsbeta}_t)\Big) dt
    \notag \\[1mm]
    &\hskip1cm \leq C  \left((1+L^{{2}}_z{+L^2_\lambda})\mathbb{E} \sup_{0\leq s\leq T} \big|X_s^{\tilde{\bsbeta}} - X_s^{\bsbeta}\big|^2 + \Gamma^2\right).   
    \label{eq:cost_cont_YT_bound_new_part1}
\end{align}
For the term involving the Wiener process we use the Ito isometry, and then~\eqref{eq:cont_bound_betas_gamma_new2}:
    \begin{align}
            &\mathbb{E} \big|\int_0^T [\tilde{z}(t, X_t^{\tilde{\bsbeta}})-z(t, X_t^{\bsbeta})] \cdot dW_t\big|^2
            \nonumber \\[1mm]
            &\hskip1cm {=}\mathbb{E}\int_0^T \big|\tilde{z}(t, X_t^{\tilde{\bsbeta}})-z(t, X_t^{\bsbeta})\big|^2dt
            \nonumber \\[1mm]
            &\hskip1cm \leq C \left(L^{{2}}_z \int_0^T \mathbb{E}\big|X_t^{\tilde{\bsbeta}} - X_t^{\bsbeta}\big|^2 dt + \Gamma^2\right)
            \nonumber \\[1mm]
            &\hskip1cm \leq C  \Big(L^{{2}}_z \mathbb{E}\sup_{s\in[0,T]}\big|X_s^{\tilde{\bsbeta}} - X_s^{\bsbeta}\big|^2 + \Gamma^2  \Big).
            \label{eq:cost_cont_YT_bound_new_part2}
    \end{align}
    
Going back to~\eqref{eq:cost_cont_YT_bound_new_2}, we now combine ~\eqref{eq:cost_cont_YT_bound_new_part1} with~\eqref{eq:cost_cont_YT_bound_new_part2}, and we use~\eqref{eq:bound-supp-diff-X-Gamma2_new} as well as the fact that $\Gamma = \|\tilde\bsbeta -  \bsbeta\|_{\VV}$, which yields:
\begin{align*}
        \mathbb{E} \Big[ |Y_T^{\tilde{\bsbeta}}-Y_T^{\bsbeta}|^2 \Big] 
        &\leq C\Big( (1 + L^{{2}}_z{+L^2_\lambda})\mathbb{E}\sup_{s\in[0,T]}\big|X_s^{\tilde{\bsbeta}} - X_s^{\bsbeta}\big|^2 +  \Gamma^2+  \mathbb{E}\big|\tilde{y}_0(X_0)- y_0(X_0)\big|^2 \Big)
        \\[1mm]
        &\leq C ((1 + L^{{2}}_z{+L^2_\lambda})C e^{C (1+L^{{2}}_z{+L^2_\lambda})} + 1) \Gamma^2.
\end{align*}
\vskip10pt
Next we note that the continuity of $\mu_T^{\bsbeta} = \mathcal{L}(X_T^{\bsbeta})$ is a direct consequence of~\eqref{eq:bound-supp-diff-X-Gamma2_new}: 
\begin{equation}
\begin{aligned}
    W^2_2(\mu^{\tilde{\bsbeta}}_T, \mu^{\bsbeta}_T) 
    &\leq\mathbb{E}\big|X_T^{\tilde{\bsbeta}} - X_T^{\bsbeta}\big|^2 \leq  C e^{C (1+L^{{2}}_z{+L^2_\lambda})} \Gamma^2.
\end{aligned}
\end{equation}
\vskip12pt
Finally, by using the  Lipschitz property of $\mbsP$  (Assumption~\ref{hyp2:Pbar}) and $g$ (Assumption~\ref{hyp2:g-lip}) in their arguments and also by using the previous bounds, we conclude that
\begin{align*}
    |\overline{\mbsP}(\tilde{\bsbeta}) - \overline{\mbsP}(\bsbeta)|^2
    &= \EE\left|\mbsP\big(Y_T^{\tilde{\bsbeta}} - g(X_T^{\tilde{\bsbeta}}, \mu_T^{\tilde{\bsbeta}}; \tilde{\lambda}(T{,X_T^{\tilde{\bsbeta}}}))\big) - \mbsP\big(Y_T^{{\bsbeta}} - g(X_T^{{\bsbeta}}, \mu_T^{{\bsbeta}}; {\lambda}(T{,X_T^{{\bsbeta}}}))\big) \right|^2
    \\[1mm]
    &\leq C \EE \left(\left|Y_T^{\tilde{\bsbeta}} - Y_T^{{\bsbeta}}\right|^2
    +
    \left| g(X_T^{\tilde{\bsbeta}}, \mu_T^{\tilde{\bsbeta}}; \tilde{\lambda}(T{,X_T^{\tilde{\bsbeta}}})) - g(X_T^{{\bsbeta}}, \mu_T^{{\bsbeta}}; {\lambda}(T{,X_T^{\tilde{\bsbeta}}}))\right|^2 \right)
    \\[1mm]
    &\leq C\left( ((1 + L^{{2}}_z{+L^2_\lambda})C e^{C (1+L^{{2}}_z{+L^2_\lambda})} + 1) + e^{C (1+L^{{2}}_z{+L^2_\lambda})}\right)\Gamma^2
    \\[1mm]
    &\leq C\left( 1 + (1 + L^{{2}}_z{+L^2_\lambda})e^{C (1+L^{{2}}_z{+L^2_\lambda})}\right)\Gamma^2.
\end{align*}
\endproof

\subsection{Finite population approximation}
\label{sec:proof-finite-pop-approx}

In this section we prove Theorem~\ref{thm:finite-pop-approx}.

\proof[Proof of Theorem~\ref{thm:finite-pop-approx}. ]
Let $\beta \in \VV$. We first analyze the difference between the mean field dynamics and the finite-population dynamics. We then evaluate the difference between the costs. 

\vskip 6pt
    \textbf{Step 1: Propagation of chaos. } 
    We will use the notation: $Q^{N,i}_t = (X^{N,i}_t,Y^{N,i}_t)$, $i=1,\dots,N$. 
    We have that $Q^{N,i}_0$ are i.i.d. with distribution $\nu_0 := \mathcal{L}(X_0,y_0(X_0))$ with $X_0 \sim \mu_0$, and moreover the family of $Q^{N,i}, i=1,\dots,N$ satisfies:
    \begin{equation}
    \label{eq:particle-dynamics-CNi}
        dQ^{N,i}_t = \boldsymbol{B}(t, Q^{N,i}_t, \nu^N_t) dt + \boldsymbol{\Sigma}(t, Q^{N,i}_t) d W^i_t
    \end{equation}
    where  $\nu^{N}_t = \frac{1}{N} \sum_{i=1}^N \delta_{Q^{N,i}_t}$ with:
    \begin{equation}\label{eq:def-B-Sigma-for-C}
        \boldsymbol{B}(t, (x,y), \nu) = \begin{pmatrix}
        b(t, x, \hat{\alpha}_t(x, \nu_x, {\sigma^{-1}}^{\top}z(t,x); \lambda(t,x)), \nu_x;\lambda(t,x))
        \\
        -f(t, x, \hat{\alpha}_t(x, \nu_x, {\sigma^{-1}}^{\top}z(t,x); \lambda(t,x)), \nu_x;\lambda(t,x))
        \end{pmatrix},
        \qquad
        \boldsymbol{\Sigma}(t, (x,y)) = \begin{pmatrix}
        \sigma
        \\
        z(t,x)
        \end{pmatrix},
    \end{equation}
    where $\nu_x$ denotes the first marginal of $\nu$.

    We now consider an analogous system but with decoupled particles, whose dynamics involve their own (true) distribution instead of the empirical one. Let $(Q^i)_{i=1,\dots,N}$ be $N$ independent copies of the MKV process $(X,Y)$ solving~\eqref{eq:FFSDE_y0z} respectively with Brownian motions $W^i, i=1,\dots,N$, and let $\nu_t = \mathcal{L}(Q^1_t)$ by their (common) law at time $t$. To be specific, $Q^i_0 = Q^{N,i}_0$, and for $t \in [0,T]$,  
        \[
        		dQ^i_t = \boldsymbol{B}(t, Q^i_t, \nu_t) dt + \boldsymbol{\Sigma}(t, Q^i_t) d W^i_t,
    	\]
	with the same notations as above.

    As we argue below, we can check that $\boldsymbol{B}$ and $\boldsymbol{\Sigma}$ are bounded on bounded subsets of $[0,T] \times \RR^{m+1} \times \mathcal{P}_2(\RR^{m+1})$, and are Lipschitz continuous in $q=(x,y)$ and $\mu$, uniformly in $t \in [0,T]$,  $\RR^{m+1}$ and $\mathcal{P}_2(\RR^{m+1})$  being equipped with the Euclidean norm and the $2$-Wasserstein distance respectively. Indeed, we have:
    \begin{itemize}
        \item For $\boldsymbol{\Sigma}$: the desired property stems from the fact that $\sigma$ is constant and, by definition of $\VV$, $x \mapsto z(t,x)$ is  Lipschitz continuous and of linear growth, hence it is bounded on bounded subsets.
        \item  For $\boldsymbol{B}$: 
        \begin{itemize}
            \item[$\circ$] $\lambda$ is Lipschitz in $x$ uniformly $t \in [0,T]$ by definition of $\VV$.
            \item[$\circ$] $\mathcal{b}$ and $\mathcal{f}$ are Lipschitz continuous with respect to their inputs, by Assumptions~\ref{hyp2:drift-lip}, \ref{hyp2:cost0-lip}, and \ref{hyp2:g-lip}. 
            \item[$\circ$] As a consequence, the functions $ (t,(x,y),\nu) \mapsto b(t, x, \hat{\alpha}_t(x, \nu_x, {\sigma^{-1}}^{\top}z(t,x); \lambda(t,x)), \nu_x;\lambda(t,x))$ and $ (t,(x,y),\nu) \mapsto b(t, x, \hat{\alpha}_t(x, \nu_x, {\sigma^{-1}}^{\top}z(t,x); \lambda(t,x)), \nu_x;\lambda(t,x))$ are bounded on bounded subsets of $[0,T] \times \RR^{m+1} \times \mathcal{P}_2(\RR^{m+1})$, and are Lipschitz continuous in $q=(x,y)$ and $\mu$, uniformly in $t \in [0,T]$,  $\RR^{m+1}$.
        \end{itemize}  
    \end{itemize}

    By~\cite[Theorem 2.12]{carmona2018probabilistic2}, the system~\eqref{eq:particle-dynamics-CNi} of $N$ SDEs in dimension $(m+1)$ has a unique solution and it satisfies:
    \begin{equation}
    \label{eq:prop-chaos-CNi}
        \max_{1 \le i \le N} \EE\left[ \sup_{0 \le t \le T} |Q^{N,i}_t - Q^i_t|^2 \right] + \sup_{0 \le t \le T} \EE \left[ W_2\left( \nu^N_t, \nu_t\right)^2\right] \le C \epsilon_{m+1}(N), 
    \end{equation}
    where $\epsilon_{d}(N)$ is defined by~\eqref{eq:chaos-epsilon-N}.
    
    \vskip 6pt
    \textbf{Step 2: Analysis of the cost. } %
    We want to upper bound the difference between the total costs $J^0_{\nu,N}(\beta)$ and $J^0_{\nu}(\beta)$. We will bound the square of their difference and we start by noticing that, by Jensen's inequality:
    \begin{equation}
    \label{eq:proof-pop-approx-Jdiff}
        |J^0_{\nu,N}(\beta) - J^0_{\nu}(\beta)|^2
        \le \frac{1}{N} \sum_{i=1}^N  \EE \left[T \int_0^T |\Delta f_0^i(t)|^2 dt + \Delta g_0^i] + \nu |\Delta \mbsP|^2 \right] 
    \end{equation}
    where
    \[
        \Delta f_0^i(t) = f_0(t, \mu_t, \lambda(t{,X^{i}_t})) - f_0(t, \mu_t^{N}, \lambda(t{,X^{N,i}_t})),
        \quad
        \Delta g_0^i = g_0(\mu_T,\lambda(T{,X^{i}_T})) - g_0(\mu_T^{N},\lambda(T{,X^{N,i}_T}))
    \]
    and
    \[
        \Delta \mbsP^i = \mbsP\Big(Y^{i}_T-g\big(X^{i}_T, \mu_T; \lambda(T{,X^{i}_T})\big)\Big) - \mbsP\Big(Y^{N,i}_T-g\big(X^{N,i}_T, \mu_T^{N}; \lambda(T{,X^{N,i}_T})\big)\Big).
    \]
    
    We first compare the terminal costs for the finite-particle setting and the mean field regime. We have:
    \begin{align*}
    	&\EE |g_0(\mu_T, \lambda(T{,X^{i}_T})) - g_0(\mu^{N}_T, \lambda(T{,X^{i,N}_T}))|^2 
    	\\[2mm]
        &\le C \EE \left[ W_2(\mu_T, \mu^{N}_T)^2 + |\lambda(T{,X^{i}_T}) + \lambda(T{,X^{i,N}_T})|^2\right] 
    	\\[2mm]
        &\le C \sup_{0 \le t \le T} \left(\EE \left[ W_2(\nu_T, \nu^{N}_T)^2 \right] +  \max_{1 \le i \le N} \EE\left[ \sup_{0 \le t \le T} |Q^{N,i}_t - Q^i_t|^2 \right] \right)
    	\\[2mm]
        &\le C \epsilon_{m+1}(N),
    \end{align*}
    where we used the Lipschitz continuity of $g_0$ in the first inequality, the Lipschitz continuity of $\lambda$ in the second inequality and~\eqref{eq:prop-chaos-CNi} in the last inequality.
  
    We then consider the running costs:
    \begin{align*}
    	&\EE\left[\int_0^T |f_0(t, \mu_t, \lambda(t{,X^i_t})) - f_0(t, \mu^N_t, \lambda(t{,X^{N,i}_t})) |^2 dt \right]
	\\
	\le
	\, &C\EE\left[\int_0^T (W_2(\mu_t, \mu^N_t)^2 + |\lambda(t{,X^i_t}) - \lambda(t{,X^{N,i}_t})|^2) dt \right]
	\\
	\le
	\, & C  T \left( \sup_{0 \le t \le T}  \EE \left[ W_2(\mu_t, \mu^{N}_t)^2 \right] +   \EE\left[\sup_{0 \le t \le T} |X^i_t - X^{N,i}_t|^2 \right] \right)
	\\
	\le
	\, & C  T \left(  \sup_{0 \le t \le T} \EE \left[ W_2(\nu_t, \nu^{N}_t)^2 \right] +   \EE\left[\sup_{0 \le t \le T} |Q^i_t - Q^{N,i}_t|^2 \right] \right)
	\\
	\le
	\, & C \epsilon_{m+1}(N), 
    \end{align*}
    where we used the Lipschitz continuity of $f_0$ in the first inequality, the Lipschitz continuity of $\lambda$ in the second inequality, and~\eqref{eq:prop-chaos-CNi} in the last inequality. 
    
    Finally, we consider the penalty:
    \begin{align*}
    	&\EE\left[ |\mbsP\Big(Y^{i}_T-g\big(X^{i}_T, \mu_T; \lambda(T{,X^{i}_T})\big)\Big) - \mbsP\Big(Y^{N,i}_T-g\big(X^{N,i}_T, \mu^{N}_T; \lambda(T{,X^{N,i}_T})\big)\Big)|^2  \right]
	\\[2mm]
	\le
	\, &C \EE\left[ |Y^{i}_T - Y^{N,i}_T|^2
	+ |g\big(X^{i}_T, \mu_T; \lambda(T{,X^{i}_T})\big) - g\big(X^{N,i}_T, \mu^{N}_T; \lambda(T{,X^{N,i}_T})\big)|^2  \right]
	\\[2mm]
	\le
	\, &C \EE\left[ |Y^{i}_T - Y^{N,i}_T|^2
	+  (|X^{i}_T-X^{N,i}_T|^2
		+ W_2(\mu_T, \mu^{N}_T)^2
		+  |X^{i}_T - X^{N,i}_T|^2) \right]
	\\[2mm]
	\le
	\, & C \left(  \sup_{0 \le t \le T} \EE \left[ W_2(\nu_t, \nu^{N}_t)^2 \right] +   \EE\left[\sup_{0 \le t \le T} |Q^i_t - Q^{N,i}_t|^2 \right] \right)
	\\[2mm]
	\le
	\, & C \epsilon_{m+1}(N),
    \end{align*}
    where we used the Lipschitz continuity of $\mbsP$ in the first inequality, the Lipschitz continuity of $g$ and $\lambda$ in the second inequality, and~\eqref{eq:prop-chaos-CNi} in the last inequality.

    Going back to~\eqref{eq:proof-pop-approx-Jdiff} and collecting the terms, the difference between the total expected cost is upper bounded by $C \sqrt{\epsilon_{m+1}(N)}$. 
    
\hfill\qed\endproof

\subsection{Discrete time approximation}
\label{sec:proof-discrete-time-approx}

In this section, we prove Theorem~\ref{thm:discrete-time-approximation}. 

\proof[Proof of Theorem~\ref{thm:discrete-time-approximation}. ]

We split the proof into two steps: first we compare continuous-time and discrete-time dynamics. Then we compare the costs. 

\vskip 6pt
{\bf Step 1: Euler-Maruyama scheme approximation. }

	We now discretize time and consider an Euler-Maruyama scheme for the dynamics~\eqref{eq:particle-dynamics-CNi} of the particle system $(Q^{N,i})_{i=1,\dots,N}$. More precisely, let us consider a uniform discretization $t_n = n \Delta t$ of $[0,T]$ with $N_T+1$ points and $\Delta t = T/N_T$. 
	Consider the discrete time process $\check{Q}^{N,i}$ defined as follows: 
	$\check{Q}^{N,i}_0 = Q^{N,i}_0$ and the family of $\check{Q}^{N,i}, i=1,\dots,N$, satisfies:
    \begin{equation*}
        \check{Q}^{N,i}_{t_{n+1}}  = \check{Q}^{N,i}_{t_{n}} + \boldsymbol{B}(t_n, \check{Q}^{N,i}_{t_n}, \check{\nu}^N_{t_n}) \Delta t + \boldsymbol{\Sigma}(t_n, \check{Q}^{N,i}_{t_n}) \Delta W^i_{t_{n+1}}\,,
    \end{equation*}
    where  $\check{\nu}^{N}_t = \frac{1}{N} \sum_{i=1}^N \delta_{\check{C}^i_t}$ and $\Delta W^i_{t_{n+1}} = W^i_{t_{n+1}} - W^i_{t_{n}}$, and $\boldsymbol{B}$ and $\boldsymbol{\Sigma}$ are as in~\eqref{eq:def-B-Sigma-for-C}.

    By~\cite[Lemma 14]{carmona2022convergence}, we obtain: 
	\begin{equation}
    \label{eq:checkC-C-bound}
    	\frac{1}{N} \sum_{i=1}^N \EE[|\check{Q}^{N,i}_{t_{n}} - Q^{N,i}_{t_{n}}|^2] \le C \Delta t,
	\end{equation}
 where $C$ is independent of $N$ but may depend on the Lipschitz constants of $b(t)$, $f(t)$ and $\beta(t)$, as well as $\beta(0,0,0)$. 

\vskip 6pt
{\bf Step 2: Analysis of the cost. }  

Our goal is now compare the difference between the discrete-time cost and the continuous-time one. We cannot apply directly~\cite[Proposition 15]{carmona2022convergence} because our problem involve the control in the terminal cost (through $\lambda$), which was not the case in the setting of~\cite[Proposition 15]{carmona2022convergence}. However, by using Assumption~\ref{assm:fC0-gC0}, it is possible to follow the proof of that result and adapt it to our setting.

We rewrite $J^0_{\nu,N}$ in terms of the particle system $(Q^{N,i})_i$ and its empirical distribution $\nu^N$ as:
\begin{equation*}
    \begin{split}
        J^0_{\nu,N}(\bsbeta) 
        &= \frac{1}{N} \sum_{i=1}^N \EE\Big[\int_0^T f^Q_0(t, \nu_t^{N}, \Lambda^Q(t,Q^{N,i}_t))dt + g^Q_0(\nu_T^{N},\Lambda^Q(T,Q^{N,i}_T))\Big]
    \end{split}
    \end{equation*}
where $\Lambda^Q(t, (x,y)) = \lambda(t,x)$ and we recall that $f^Q_0$ and $g^Q_0$ are defined in~\eqref{eq:def-fC0-gC0}. We can similarly rewrite $\check{J}^0_{\nu,N,N_T}$ as:
\[
    \check{J}^0_{\nu,N,N_T}(\ubsbeta) 
    = \frac{1}{N} \sum_{i=1}^N \EE\Big[\sum_{n=0}^{N_T-1} f^Q_0(t_n, \check{\nu}_{t_n}^{N}, \Lambda^Q(t_n, \check{Q}^{N,i}_{t_n})) \Delta t + g^Q_0(\check{\nu}_T^{N},\Lambda^Q(T,\check{Q}^{N,i}_T)) \Big].
\]

We have:
\[
    |J^0_{\nu,N}(\beta) - \check{J}^0_{\nu,N,N_T}(\beta)| \le \delta_{f} + \delta_{g}
\]
with 
\[
    \delta_{f} = \left|\frac{1}{N} \sum_{i=1}^N \EE\Big[\sum_{n=0}^{N_T-1} \int_{t_n}^{t_{n+1}} \left(f^Q_0(t, \nu_t^{N}, \Lambda^Q(t,Q^{N,i}_t)) - f^Q_0(t_n, \check{\nu}_{t_n}^{N}, \Lambda^Q(t_n, \check{Q}^{N,i}_{t_n})) \right) dt \Big]\right|
\]
and
\[
    \delta_{g} = \left| \frac{1}{N} \sum_{i=1}^N \EE\Big[g^Q_0(\nu_T^{N},\Lambda(T,Q^{N,i}_T)) - g^Q_0(\check{\nu}_T^{N},\Lambda^Q(T,\check{Q}^{N,i}_T))\Big]\right|.
\]
For the term $\delta_{f}$, we can follow the analysis of the corresponding term in the proof of ~\cite[Proposition 15]{carmona2022convergence}, which yields $\delta_f \le C\sqrt{\Delta t}$.  

For the term $\delta_{g}$, we have:
\begin{align*}
    \delta_{g} = \left| \frac{1}{N} \sum_{i=1}^N \EE\Big[G^Q(Q^{N,i}_T, \nu_T^{N}) - G^Q(\check{Q}^{N,i}_T, \check{\nu}_T^{N})\Big]\right|, 
\end{align*}
where we used the notation:
$
    G^Q(q, \nu) = G^Q_0(\nu, \Lambda^Q(T, q)).
$ 

Following the first part of the proof of~\cite[Proposition 15]{carmona2022convergence}, we obtain:
\[
    \delta_{g} \le \left( \frac{1}{N} \sum_{i=1}^N \EE\left[ \left| \int_0^1 \partial_\nu \boldsymbol{G}^Q(\nu_\theta)(Q^{N,i}_\theta) d\theta \right|^2 \right]\right)^{1/2} \left( \frac{1}{N} \sum_{i=1}^N \EE \left[ |\check{Q}^{N,i}_T - Q^{N,i}_T|^2 \right] \right)^{1/2},
\]
where $\nu_\theta = \frac{1}{N} \sum_{i=1}^N \delta_{Q^{N,i}_\theta}$ with $Q^{N,i}_\theta = Q^{N,i}_T + \theta (\check{Q}^{N,i}_T - Q^{N,i}_T) $,  and 
$
    \boldsymbol{G}^Q(\nu) = \int G^Q(q, \nu) \nu(dq).
$ 
By~\eqref{eq:checkC-C-bound}, the second term in the product is bounded by $C\sqrt{\Delta t}$. Using the fact that the partial derivatives $\partial_c G^Q(q,\nu)$ and $\partial_\nu G^Q(q,\nu)$ are of linear growth, the first term can also be bounded by $C\sqrt{\Delta t}$. 

\hfill\qed\endproof

\subsection{Neural network approximation}
\label{sec:proof-main-thm-approx-NN} 

The goal of this section is to analyze the approximation due to the neural network control. We first give the proof of Theorem~\ref{thm:aprpox-NN-timestep} by using Proposition~\ref{prop:approx-fct-NN} and Lemma~\ref{lem:J0_P_regularity_compact-discrete}, which are stated and proved in the next subsections.  

\subsubsection{Proof of Theorem~\ref{thm:aprpox-NN-timestep} using auxiliary results}

\proof[Proof of Theorem~\ref{thm:aprpox-NN-timestep}. ]
    Fix $n_{\mathrm{in}}>n_0$ where $n_0$ is as in Proposition~\ref{prop:approx-fct-NN} below. Let $R = n_{\mathrm{in}}^{1/(3(m+1))}$. We want to apply this result for each time $t_n$, $n=0,\dots,N_T-1$. Note that for every $n$, $\beta(t_n,\cdot) \in \cC_{m+1,d+2}(R, K, L)$ so, by Proposition~\ref{prop:approx-fct-NN} below, there exists a one-hidden layer neural network $\beta_{t_n,\theta_n} \in \bN^{\psi}_{m+1, n_{\mathrm{in}}, {d+2}}$ with parameters $\theta_n$ such that
	$$
		\|\beta(t_n,\cdot) - \beta_{t_n,\theta_n}\|_{\cC^0(\bar{B}_{m+1}(0,R); \RR^{d+2})} \le C_0 (1+R)n_{\mathrm{in}}^{-1/(2 (m+1)}, 
	$$
and such that the Lipschitz constant of $\beta_{t_n,\theta_n}$ is at most $C_0 (1+Rn_{\mathrm{in}}^{-1/(2 (m+1))})$.  Here 
 $C_0$ is a constant depending only on the data of the problem, on $N_T$, and on $\psi$ through $\hat\psi_1$ and $\|\psi'\|_{\cC^0(\TT^{m+1})}$. 
Furthermore, by the above bound and since $\beta \in \VV$,
\[
    |z_{t_n,\theta_n}(x)| 
    \le |z(t_n,x)| + C_0 (1+R)n_{\mathrm{in}}^{-1/(2 (m+1))} 
    \le K'(1 + |x|).
\] 
As a consequence, $\hat\beta_{t_n,\theta_n} \in \VV_{N_T}(L_z', L_\lambda')$ with  $L_z = L_\lambda = \max\{C_0(1+Rn_{\mathrm{in}}^{-1/(2 (m+1))}), K'  \}$, which is smaller than $\max\{2 C_0, K'\}$ by our choice of $R$.

Let $\hat\ubsbeta_\theta = (\hat\beta_{t_n,\theta_n})_{n=0,\dots,N_T}$. Then the condition in Lemma~\ref{lem:J0_P_regularity_compact-discrete} below is satisfied with $\ubsbeta = \hat\ubsbeta$, $\tilde\ubsbeta = \hat\ubsbeta_\theta$ and $\Gamma = C_0 (1+R)n_{\mathrm{in}}^{-1/(2 (m+1))}$. We deduce that:
\begin{equation*}
    |J^0_{N,N_T}(\hat\ubsbeta_\theta) - J^0_{N,N_T}(\hat\ubsbeta)| 
    \leq C_B\Big(\Gamma^2+\frac{1}{R}\Big)^{1/2}\quad 
    \text{and}\quad |\overline{\mbsP}(\hat\ubsbeta_\theta) - \overline{\mbsP}(\hat\ubsbeta)| 
    \leq C_B\Big(\Gamma^2+\frac{1}{R}\Big)^{1/2},
\end{equation*}
where $C_B$ can depend on the data of the problem and $C_0$. 
So:
\begin{equation*}
    |J^0_{\nu,N,N_T}(\hat\ubsbeta_\theta) - J^0_{\nu,N,N_T}(\hat\ubsbeta)| 
    \leq C \Big(\Gamma^2+\frac{1}{R}\Big)^{1/2},
\end{equation*}
where $C$ depends on $C_B$ and also on the penalty parameter $\nu$.

To conclude, we recall that $R = n_{\mathrm{in}}^{1/(3(m+1))}$, hence
 \begin{align*}
 	\Gamma^2 + \frac{1}{R} 
	&= C_0^2  \left(1 + n_{\mathrm{in}}^{1/(3(m+1))}  \right)^2 n_{\mathrm{in}}^{-1/(m+1)} + n_{\mathrm{in}}^{-1/(3(m+1))} 
	\\
	&\le C_1 n_{\mathrm{in}}^{-1/(3(m+1))},
 \end{align*}
where $C_1$ depends only on $C_0$. This completes the proof.
\hfill\qed\endproof

\subsubsection{Neural network approximation for one time step}

To prove Theorem~\ref{thm:aprpox-NN-timestep}, we make use of an approximation result for neural networks. More precisely, we will use~\cite[Theorems 2.3 and 6.1]{MhaskarMicchelli}. These results are more precise than simple density results such as~\cite{cybenko1989approximation,hornik1990universal,hornik1991approximation,pinkus1999approximation} because they provide rates of convergence. Furthermore, they provide an explicit construction of a neural network achieving a given accuracy.

These results are stated for functions on a compact set, and more precisely on the torus (i.e., periodic function). So we use them to approximate control functions locally in space. This is sufficient for our purpose because the processes remain bounded with high probability. Similar techniques were used in~\cite{carmona2022convergence} to approximate the whole control, as a function of the time and space variables because the neural network approximation was used \emph{before} the time discretization. In contrast, here we are going to use the neural network approximation after the time discretization, and we are going to approximate the control function at each time step by a different neural network. This will allow us to work with less regular neural networks and hence to use less restrictive assumptions on the regularity of the control to be approximated.

We will use the following notations. For a positive integer $d$, $\TT^d$ denotes the $2\pi-$torus in dimension $d$.
For positive integers $n$ and $d$, and a function $g \in \cC^0(\TT^d) $, $E_n^d(g)$ denotes the trigonometric degree of approximation of $g$ defined by
$
	E_n^d(g) = \inf_{P} \|g - P\|_{\cC^0(\TT^d)},
$
where the infimum is over trigonometric polynomials of degree at most $n$ in each of its $d$ variables.

In the sequel, we denote by $\cC_{d,k}(R, K, L)$ the class of functions $\ff: \bar{B}_d(0,R) \ni x \mapsto \ff(x) \in \RR^{k}$ such that $\ff$ is Lipschitz continuous with $\cC^0(\bar{B}_d(0,R))-$norm bounded by $K$ and Lipschitz constant bounded by $L$.
We will sometimes use the notation $\nabla \ff = (\partial_{x_1} \ff, \dots, \partial_{x_d} \ff)$ for the gradient.

We state a reformulation of results from Mhaskar and Micchelli~\cite{MhaskarMicchelli}. 
\begin {theorem}[Theorems 2.3 and 6.1 in~\cite{MhaskarMicchelli}]
\label{th:MM}
	There exists an absolute constant $C_{mm}$ such that for every positive integers $d, n$ and $N$, there exists a positive integer $n_{\mathrm{in}} \le C_{mm} N n^d$ with the following property.
	 For any $\ff: \TT^{d} \to \RR$ of class $\cC^{1}(\TT^d)$, i.e., continuously differentiable, there exists $\varphi_\ff \in \bN^\psi_{d, n_{\mathrm{in}}, 1}$ such that:
	\begin{align}
		\label{eq:MhaskarMicchelli-bdd-f}
		\|\ff - \varphi_\ff\|_{\cC^0(\TT^{d})} &\leq c\left[E^{d}_n(\ff) + E^1_N( \psi) n^{d/2} \|\ff\|_{\cC^0(\TT^{d})} \right],
	\end{align}
	where the constant $c$ depends only on $d$ and $\hat\psi_1$, and where the degree of trigonometric approximation $E^m_n(\cdot)$ was defined in the text above.
\end{theorem}

Based on the above result, one can obtain the following approximation property, that we use to approximate the control by a different neural network at each time step in the proof of Theorem~\ref{thm:aprpox-NN-timestep}.

\begin{proposition}
\label{prop:approx-fct-NN}
For every real numbers $K>0$ and $L>0$, there exists a constant $C$ depending only on the above constants, on $d,k,T,$ and on the activation function through $\hat\psi_1$ and $\|\psi'\|_{\cC^0}$, and there exists a constant $n_0$ depending only on $C_{mm}$ and $d$ with the following property. For every $R>0$, for every $\ff \in \cC_{d,k}(R, K, L)$ and for every integer $n_{\mathrm{in}}>n_0$, there exists a one-hidden layer neural network $\varphi_{\ff} \in \bN^{\psi}_{{d}, n_{\mathrm{in}}, {d+2}}$ such that
	$$
		\|\ff - \varphi_{\ff}\|_{\cC^0(\bar{B}_d(0,R); \RR^{k})} \le C (1+R)n_{\mathrm{in}}^{-1/(2 d)}, 
	$$
and such that the Lipschitz constant of $\varphi_{\ff}$ is at most $C (1+Rn_{\mathrm{in}}^{-1/(2 d)})$.
\end{proposition}
We stress that the constant $C$ in the above statement does not depend on $R$.

Proposition~\ref{prop:approx-fct-NN} is proved by adapting the arguments in the proof of Proposition~\cite[Proposition 10]{carmona2022convergence}, which was for the case of $\mathfrak{f}$ of class $\cC^3$. Since the proof is very similar, we omit it for the sake of brevity.

\subsubsection{Discrete time analysis of the cost sensitivity}

Next, we prove that the (discrete-time, discrete-population) cost does change too much when we approximate the control by another control on a compact set.

Let $\VV_{N_T}$ be the set of discrete time controls defined as:
\begin{align*}
    \VV_{N_T} &= \Big\{ \underline\bsbeta = (\boldsymbol \lambda_{t_n}, \boldsymbol z_{t_n}, y_0)_{n=0,\dots,N_T} \in (\cC^0({\RR^m}; \RR) \times \cC^0(\RR^m; \RR^d)  \times \cC^0(\RR^m; \RR))^{N_T+1} \,:\, 
    \\
    &\qquad\qquad\qquad{\hbox{for every $n=0,\dots,N_T-1$, } \lambda_{t_n} \hbox{ is $L_{\lambda}$-Lipschitz continuous,}}
    \\
    &\qquad\qquad\qquad  {\hbox { and has linear growth } |\lambda_{t_n}(x)|<L_{\lambda}(1+|x|), x \in \RR^m, }
    \\
    &\qquad\qquad\qquad  \hbox{for every $n=0,\dots,N_T-1$, } z_{t_n} \hbox{ is $L_z$-Lipschitz continuous }
    \\
    &\qquad\qquad\qquad  \hbox { and has linear growth } |z_{t_n}(x)|<L_z(1+|x|), x \in \RR^m  \Big\}.
\end{align*}

\begin{lemma}
\label{lem:J0_P_regularity_compact-discrete}
Suppose Assumption~\ref{hyp2} holds. 
Let $\tilde{\ubsbeta}$ and ${\ubsbeta} \in \VV_{N_T}$ be two discrete-time controls. Then, there exists a constant $C_B$ depending only on the parameters of the model, on $L_\lambda$, on $L_z$ and on $N_T$, such that for all $\Gamma>0$ and $R>0$, if 
\begin{equation}
      \max_{n=0,\dots,N_T}\lVert \tilde{\beta}_{t_n} - {\beta}_{t_n} \rVert_{\mathcal{C}^0(\overline{B}_m(0,R)\times\overline{B}_m(0,R))} \leq \Gamma,
\end{equation}
where $\overline{B}_m(0,R)$ denotes the closed m-dimensional ball with radius R, then:
\begin{equation}
\label{eq:J0_P_bound_compact}
    |J^0_{N,N_T}(\tilde{\ubsbeta}) - J^0_{N,N_T}({\ubsbeta})| 
    \leq C_B\Big(\Gamma^2+\frac{1}{R}\Big)^{1/2}\quad 
    \text{and}\quad |\overline{\mbsP}(\tilde{\ubsbeta}) - \overline{\mbsP}({\ubsbeta})| 
    \leq C_B\Big(\Gamma^2+\frac{1}{R}\Big)^{1/2}. 
\end{equation}
\end{lemma}

\begin{remark} 
We stress that in Lemma~\ref{lem:J0_P_regularity_compact-discrete}, we only assume that the controls are close on the compact set $\overline{B}_m(0,R)\times\overline{B}_m(0,R)$ and not necessarily outside of this compact set. This is justified by the fact we want to apply this result to a neural network approximation of a given control $\bsbeta$.
\end{remark}

\proof[Proof of Lemma~\ref{lem:J0_P_regularity_compact-discrete}. ]
At a high level, the proof of Lemma~\ref{lem:J0_P_regularity_compact-discrete} shares some similarities with the proof of Lemma~\ref{lem:cont-J-P} and draws inspiration from \cite[Prop. 13]{carmona2022convergence} but there are several important differences. First, we need to treat differently the situation inside and outside the compact set $\overline{B}_m(0,R)\times\overline{B}_m(0,R)$. Second, the population distribution and the time are discrete, which requires ad-hoc estimates. 

Below $C$ denotes a generic constant whose value can change from one line to the other and is independent of $\Gamma$ and $R$ (but it can depend on $N_T$).

\vskip6pt
\noindent \textbf{Step 1:} We first start by bounding the probability that a particle of the interacting system exits a bounded domain.

Denote by $\check\bsX^{i,\ubsbeta}$ the solution of the following system controlled by $\ubsbeta$:
\begin{equation}
    d \check{X}^{i,\ubsbeta}_{t_{n+1}} = \mathcal{b}(\check\varepsilon^{i,\ubsbeta}_{t_n}) \Delta t +\sigma \Delta W^i_{t_{n+1}}, \quad n=0,\dots,N_T-1, \qquad \check{X}^{N,i}_0 \sim \mu_0,
\end{equation}
where $\mathcal{b}(\check\varepsilon^{i, \bsbeta}_{t_n})=\mathcal{b}(t, \check{X}^{i, \bsbeta}_{t_n}, \beta_{t_n}(\check{X}^{i,\ubsbeta}_0, \check{X}^{i,\ubsbeta}_{t_n}), \check\mu^{N,\bsbeta}_{t_n})$.
Since we assume $\int |x|^4\mu_0(dx)< \infty$ (Assumption~\ref{hyp2:mu0}), there exists a constant $C$ only depending on the data of the problem and the Lipschitz constant of the controls $\tilde{\ubsbeta}$ and ${\ubsbeta}$ such that, for $\check{\ubsbeta} \in\{\tilde{\ubsbeta},\ubsbeta\}$,
\begin{equation}
    \mathbb{E}\Big[\max_{n=0,\dots,N_T} |\check{X}_{t_n}^{i, \check{\ubsbeta}}|^2\Big]\leq C, \qquad \mathbb{E}\Big[\max_{n=0,\dots,N_T} |\check{X}_{t_n}^{i, \check{\ubsbeta}}|^4\Big]\leq C.
\end{equation}
Let $\overline{\mcE}^{i,\check{\ubsbeta}}_R$ denote the complement of the following event
\begin{equation}
    {\mcE}^{i,\check{\ubsbeta}}_R = \Big\{\max_{n} |X_{t_n}^{i, \check{\ubsbeta}}|\leq R\Big\}.
\end{equation}
By using Markov's inequality, we conclude that:
\begin{equation}
\label{eq:bound-complement-CR2-discrete}
    \mathbb{P}\big[\overline{\mcE}^{i, \check{\ubsbeta}}_R\big]\leq \dfrac{C}{R^2}.
\end{equation}
\vskip6pt

\noindent \textbf{Step 2:} Next, we show that $\check{X}^{i}$ is continuous with respect to the control, i.e., we bound $\max_{n}|\check{X}^{i, \tilde{\ubsbeta}}_{t_n}- \check{X}^{i, {\ubsbeta}}_{t_n}|$ in terms of $\Gamma$ and the probability that a particle of the interacting system exits a bounded domain.
We assume that $\tilde{\ubsbeta}$ and ${\ubsbeta}$ are in $\VV_{N_T}$ (and therefore Lipschitz continuous at each time $t_n$) and that, for each $n=0,\dots,N_T$, $\tilde{\beta}_{t_n}$ and ${\beta}_{t_n}$ are close to each other on a compact set $\overline{B}_d(0,R)\times\overline{B}_d(0,R)$ such that 
\begin{equation}
\label{eq:control_dist-discrete}
    \max_{n=0,\dots,N_T} \lVert \tilde{\beta}_{t_n} - {\beta}_{t_n} \rVert_{\mathcal{C}^0(\overline{B}_d(0,R)\times\overline{B}_d(0,R))} \leq \Gamma.
\end{equation}

By using \eqref{eq:control_dist-discrete}, triangle inequality, Lipschitz continuity and linear growth of feedback controls  (see definition of $\VV_{N_T}$) we have:
\begin{align*}
    &|\tilde\beta_{t_n}(X_0^{i}, X^{i,\tilde{\ubsbeta}}_{t_n}) - \beta_{t_n}(X_0^{i}, X^{i,{\ubsbeta}}_{t_n})|\\[2mm]
    &\hskip2cm\leq \mathbbm{1}_{\mcE^{i,{\ubsbeta}}_R}\Big[|\tilde{\beta}_{t_n}(X_0^{i}, X_{t_n}^{i,\tilde{\ubsbeta}}) - \tilde{\beta}_{t_n}(X_0^{i}, X_{t_n}^{i,{\ubsbeta}})| +|\tilde{\beta}_{t_n}(X_0^{i}, X_{t_n}^{i,{\ubsbeta}}) - {\beta}_{t_n}(X_0^{i}, X_{t_n}^{i,{\ubsbeta}})|\Big]\\
    &\hskip4cm  + \mathbbm{1}_{\overline{\mcE}^{i,{\ubsbeta}}_R}\Big[|\tilde{\beta}_{t_n}(X_0^{i}, X_{t_n}^{i,\tilde{\ubsbeta}})| + |{\beta}_{t_n}(X_0^{i}, X_{t_n}^{i,{\ubsbeta}})|\Big]\\[2mm]
    &\hskip2cm\leq C \Big\{\mathbbm{1}_{\mcE^{i,{\ubsbeta}}_R}\Big[ |X_{t_n}^{i,\tilde{\ubsbeta}} - X_{t_n}^{i,{\ubsbeta}}|+ \Gamma\Big] + \mathbbm{1}_{\overline{\mcE}^{i,{\ubsbeta}}_R}\Big[|X_{t_n}^{i,\tilde{\ubsbeta}}| + |X_{t_n}^{i,{\ubsbeta}}|+1\Big]\Big\}\\[2mm]
    &\hskip2cm\leq C \Big\{\Big[ |X_{t_n}^{i,\tilde{\ubsbeta}} - X_{t_n}^{i,{\ubsbeta}}|+ \Gamma\Big] + \mathbbm{1}_{\overline{\mcE}^{i,{\ubsbeta}}_R}\Big[|X_{t_n}^{i,\tilde{\ubsbeta}}| + |X_{t_n}^{i,{\ubsbeta}}|+1\Big]\Big\}.
\end{align*}

We focus on bounding the following expression and make use of Cauchy-Schwarz inequality:
\begin{equation}
\label{eq:bound-diff-betas-GammaR2-discrete}
    \begin{aligned}
         &\mathbb{E}[|\tilde\beta_{t_n}(X_0^{i}, X^{i,\tilde{\ubsbeta}}_{t_n}) - \beta_{t_n}(X_0^{i}, X^{i,{\ubsbeta}}_{t_n})|^2]
         \\[2mm]
         & \hskip1cm \leq C\Big\{\mathbb{E}\big[|X_{t_n}^{i,\tilde{\ubsbeta}} - X_{t_n}^{i,{\ubsbeta}}|^2\big] +\Gamma^2 + \sqrt{\mathbb{P}(\overline{\mcE}_R^{i,{\ubsbeta}})}\big(\mathbb{E}|X_{t_n}^{i,\tilde{\ubsbeta}}|^2 +\mathbb{E}| X_{t_n}^{i,{\ubsbeta}}|^2+1\big)\Big\}
         \\[2mm]
         & \hskip1cm \leq C\Big\{\mathbb{E}\big[|X_{t_n}^{i,\tilde{\ubsbeta}} - X_{t_n}^{i,{\ubsbeta}}|^2\big] +\Gamma^2 + \dfrac{1}{R}\Big\},
    \end{aligned}
\end{equation}
where the last inequality holds by~\eqref{eq:bound-complement-CR2-discrete}.

By using the notation introduced in the proof of Lemma~\ref{lem:cont-J-P}, remember we have the following Euler-Maruyama scheme for any $\check\ubsbeta\in \{\tilde\ubsbeta, \ubsbeta\}$ (with a slight abuse of notation since now the controls are in discrete time only), 
\begin{equation}
\label{eq:SDE_newnotation-discrete}
    \begin{aligned}
        \check{X}^{i,\check\ubsbeta}_{t_{n+1}} &= \mathcal{b}(\check\varepsilon_{t_{n}}^{i,\check\ubsbeta}) \Delta t +\sigma \Delta \check W^{i}_{t_{n+1}}, \qquad &&\check{X}^{i,\check\ubsbeta}_0 = X^i_0,
        \\
        \check{Y}^{i,\check\bsbeta}_{t_{n+1}} &= -\mathcal{f}(\check\varepsilon_{t_{n}}^{i,\check\ubsbeta}) \Delta t+ \check{z}_{t_n}( \check{X}_{t_n}^{i,\check{\bsbeta}}) \cdot \Delta \check W^{i}_{t_{n+1}}, && \check{Y}^i_0 = \check{y}_0(X^{i,\check\bsbeta}_0),
    \end{aligned}
\end{equation}
where the $X^i_0$ are i.i.d with distribution $\mu_0$.

Using these dynamics, for all $n \in \{0, \dots, N_T\}$ we have 
\begin{equation}
\label{eq:cost_cont_X_bound-discrete}
    \begin{aligned}
        &\mathbb{E} \Big[\max_{n' \in \{0,\dots,n\}} |\check{X}_{t_{n'}}^{i,\tilde{\ubsbeta}}-\check{X}_{t_{n'}}^{i,{\ubsbeta}}|^2 \Big]
        =
        \mathbb{E}\max_{n' \in \{0,\dots,n\}} \left|\sum_{\ell=0}^{n'} [\mathcal{b}(\check\varepsilon_{t_{\ell}}^{i,\tilde\ubsbeta}) - \mathcal{b}(\check\varepsilon_{t_{\ell}}^{i,\ubsbeta})] \Delta t\right|^2 .
    \end{aligned}
\end{equation}

     By using Cauchy-Schwarz inequality and the Lipschitz property of the drift (Assumption~\ref{hyp2:drift-lip}) as we did in the proof of Lemma~\ref{lem:cont-J-P}, we can conclude that
\begin{equation}
    \begin{aligned}
         & \mathbb{E}\max_{n' \in \{0,\dots,n\}} \left|\sum_{\ell=0}^{n'} [\mathcal{b}(\check\varepsilon_{t_{\ell}}^{i,\tilde\ubsbeta}) - \mathcal{b}(\check\varepsilon_{t_{\ell}}^{i,\ubsbeta})] \Delta t\right|^2&
         \\
         &\hskip1cm\leq C \sum_{n'=0}^{n} \mathbb{E} \max_{0\leq \ell \leq n'} \big|\check{X}_{t_\ell}^{i,\tilde{\ubsbeta}} - \check{X}_{t_\ell}^{i,{\ubsbeta}}\big|^2 \Delta t  + C\left(\Gamma^2+ \dfrac{1}{R}\right)+ C \sum_{n'=0}^{n} \EE [W^2_2(\check\mu^{N,\tilde{\bsbeta}}_{t_{n'}}, \check\mu^{N,{\bsbeta}}_{t_{n'}})] \Delta t,
    \end{aligned}
\end{equation}
where we used the bound in~\eqref{eq:bound-diff-betas-GammaR2-discrete}, which differs from the proof of Lemma~\ref{lem:cont-J-P}.

Note that: 
\begin{equation}
    \label{eq:proof-discrete-approx-measure-continuity-beta}
    W^2_2(\check\mu^{N,\tilde{\ubsbeta}}_{t_{n'}}, \check\mu^{N,{\ubsbeta}}_{t_{n'}})
    \leq \frac{1}{N} \sum_{i=1}^N |\check{X}_{t_{n'}}^{i,\tilde{\ubsbeta}}-\check{X}_{t_{n'}}^{i,{\ubsbeta}}|^2 
    \leq \frac{1}{N} \sum_{i=1}^N \max_{0\leq \ell \leq n'} \big|\check{X}_{t_\ell}^{i,\tilde{\ubsbeta}} - \check{X}_{t_\ell}^{i,{\ubsbeta}}\big|^2,
\end{equation}
and $\big|\check{X}_{t_\ell}^{i,\tilde{\ubsbeta}} - \check{X}_{t_\ell}^{i,{\ubsbeta}}\big|^2$, $i=1,\dots,N$, are i.i.d.. 
In this way, the right hand side in \eqref{eq:cost_cont_X_bound-discrete} can be bounded as follows:
\begin{equation}
    \begin{aligned}
         \mathbb{E} \Big[\max_{n' \in \{0,\dots,n\}} |\check{X}_{t_{n'}}^{i,\tilde{\ubsbeta}}-\check{X}_{t_{n'}}^{i,{\ubsbeta}}|^2 \Big] %
        \hskip0.5cm \leq C\Bigg(  \sum_{n'=0}^{n-1}  \mathbb{E} \max_{0\leq \ell \leq n'} \big|\check{X}_{t_\ell}^{i,\tilde{\ubsbeta}} - \check{X}_{t_\ell}^{i,{\ubsbeta}}\big|^2 \Delta t  + \Big(\Gamma^2+ \dfrac{1}{R}\Big) \Bigg).
    \end{aligned}
\end{equation}
We recall the discrete-time Gr\"onwall's inequality: if $(\alpha_n)_n$ and $(w_n)_n$ are sequences of nonnegative numbers satisfying, for all $n \ge 0$, $\alpha_n \le c + \sum_{\ell=0}^{n-1} w_\ell \alpha_\ell$, then
    $\alpha_n \le c \exp(\sum_{\ell=0}^{n-1} w_\ell)$. 
We apply it with $\alpha_n = \mathbb{E} \Big[\max_{n' \in \{0,\dots,n\}} |\check{X}_{t_{n'}}^{i,\tilde{\ubsbeta}}-\check{X}_{t_{n'}}^{i,{\ubsbeta}}|^2 \Big]$ and $w_n = \Delta t$, and we  conclude that
\begin{equation}
\label{eq:bound-supp-diff-X-GammaR2-discrete}
    \mathbb{E} \Big[\max_{n' \in \{0,\dots,n\}} |\check{X}_{n'}^{i,\tilde{\ubsbeta}}-\check{X}_{n'}^{i,{\ubsbeta}}|^2 \Big] \leq C \Big(\Gamma^2+ \dfrac{1}{R}\Big) .
\end{equation}
\vskip6pt
\noindent \textbf{Step 3:} Finally, we can bound $|J^0_{N,N_T}(\tilde{\ubsbeta}) - J^0_{N,N_T}(\ubsbeta)|$ in terms of $\Gamma$ and $R$.

By using the Lipschitz continuity assumption on $f_0$ and $g_0$ (Assumption~\ref{hyp2:cost0-lip}), we infer that
\begin{equation}
    \begin{aligned}
        &\Big|J^0_{N,N_T}(\tilde{\ubsbeta}) - J^0_{N,N_T}(\ubsbeta) \Big| \\[1mm]
        &\hskip0.5cm\leq \frac{1}{N}\sum_{i=1}^N \mathbb{E}\Big[
        \sum_{n=0}^{N_T} \left( W^2_2(\check\mu^{N,\tilde{\bsbeta}}_{t_n}, \check\mu^{N,{\bsbeta}}_{t_n}) + \Big|\tilde{\lambda}_{t_n}(\check{X}^{i,\tilde\ubsbeta}_{t_n}) - {\lambda}_{t_n}(\check{X}^{i,\ubsbeta}_{t_n})\Big| \right) \Delta t 
          \\
          &\qquad\qquad\qquad\qquad + \Big|W_2(\check\mu_T^{N,\tilde{\ubsbeta}}, \check\mu_T^{N,{\ubsbeta}})\Big| + \Big|\tilde{\lambda}_{T}(\check{X}^{i,\tilde\ubsbeta}_{T}) - {\lambda}_{T}(\check{X}^{i,\ubsbeta}_{T})\Big| 
        \Big]
        \\[1mm]
         &\hskip0.5cm\leq C \Big( \sum_{n=0}^{N_T} \Big( \mathbb{E}\Big[\max_{0\leq \ell \leq N_T}|\check{X}_{t_\ell}^{1,\tilde{\ubsbeta}} - \check{X}_{t_\ell}^{1,\ubsbeta}|^2\Big]\Big)^{1/2} \Delta t + \big(\Gamma^2 + \dfrac{1}{R}\big)^{1/2}\Big)\\[1mm]
          &\hskip0.5cm\leq C \Big(\Gamma^2 + \dfrac{1}{R}\Big)^{1/2}.
    \end{aligned}
\end{equation}
\vskip6pt
\textbf{Step 4:} In order to bound $|\overline{\mbsP}(\tilde{\bsbeta}) - \overline{\mbsP}(\bsbeta)|$ in terms of $\Gamma$, we will show that $\check{Y}_T^{i,\ubsbeta}$, $\check{Y}_T^{i,\ubsbeta}$ and $\check{\mu}_T^{N, \ubsbeta}$ are continuous with respect to control $\ubsbeta$. We first start by showing that $\check{Y}^{i,\ubsbeta}_T$ is continuous with respect to control $\ubsbeta$. Remember that our SDE system with the new notations is given in \eqref{eq:SDE_newnotation_new}. 
We find a bound for $\mathbb{E}[|\check{Y}^{i,\ubsbeta}_T(\tilde{\ubsbeta}) - \check{Y}^{i,\ubsbeta}_T({\ubsbeta})|^2]$ (for any $i$ since $|\check{Y}^{i,\ubsbeta}_T(\tilde{\ubsbeta}) - \check{Y}^{i,\ubsbeta}_T({\ubsbeta})|$ are i.i.d.) as follows:
\begin{equation}
\label{eq:cost_cont_YT_bound_new-discrete}
    \begin{aligned}
        &\mathbb{E} \Big[ |Y_T^{i,\tilde{\ubsbeta}}-Y_T^{i,\ubsbeta}|^2 \Big]
        \\
        & \leq C \Big( \mathbb{E}\big|\tilde{y}_0(X^{i}_0)- y_0(X^{i}_0)\big|^2 +  \mathbb{E} \big|\sum_{n=0}^{N_T-1} [\mathcal{f}(\check\varepsilon_{t_n}^{i,\tilde{\ubsbeta}})-\mathcal{f}(\check\varepsilon_{t_n}^{i,\ubsbeta})] \Delta t\big|^2 
        + \mathbb{E}  \big|\sum_{n=0}^{N_T-1} [\tilde{z}_{t_n}(\check{X}_{t_n}^{i,\tilde{\ubsbeta}})-z_{t_n}(\check{X}_{t_n}^{i,\ubsbeta})]\Delta W_{t_{n+1}}^{i}\big|^2\Big).
    \end{aligned}
\end{equation} 

The first term in the right hand side can be bounded by using~\eqref{eq:bound-diff-betas-GammaR2-discrete} and~\eqref{eq:bound-supp-diff-X-GammaR2-discrete}. In order to bound the second term, we use the Lipschitz property of the function $\mathcal{f}$ (Assumption~\ref{hyp2:f-lip}), and the bound given in \eqref{eq:bound-diff-betas-GammaR2-discrete}: 
\begin{equation}
\begin{aligned}
    \mathbb{E} \big|\sum_{n=0}^{N_T-1} [\mathcal{f}(\check\varepsilon_{t_n}^{i,\tilde{\ubsbeta}})-\mathcal{f}(\check\varepsilon_{t_n}^{i,\ubsbeta})] \Delta t \big|^2
    \leq C \left(\mathbb{E} \max_{0 \le n \le N_T} \big|\check{X}_{t_n}^{i,\tilde{\ubsbeta}} - \check{X}_{t_n}^{i,\ubsbeta}\big|^2 + \Gamma^2 +\dfrac{1}{R}\right)   
\end{aligned}
\end{equation}
For the term involving increments of Wiener process, we use Cauchy-Schwarz inequality and the bounds given in~\eqref{eq:bound-diff-betas-GammaR2-discrete} and~\eqref{eq:bound-supp-diff-X-GammaR2-discrete}. 
    \begin{equation}
        \begin{aligned}
            &\mathbb{E}  \big|\sum_{n=0}^{N_T-1} [\tilde{z}_{t_n}(\check{X}_{t_n}^{i,\tilde{\ubsbeta}})-z_{t_n}(\check{X}_{t_n}^{i,\ubsbeta})]\Delta W_{t_{n+1}}^{i}\big|^2
            \\
            &\le 
             N_T \sum_{n=0}^{N_T-1} \mathbb{E} [\big|\tilde{z}_{t_n}(\check{X}_{t_n}^{i,\tilde{\ubsbeta}})-z_{t_n}(\check{X}_{t_n}^{i,\ubsbeta})\big|^2] \, \mathbb{E} [\big|\Delta W_{t_{n+1}}^{i}\big|^2]
            \\
            &=
             N_T L_z^2\sum_{n=0}^{N_T-1} \mathbb{E} [\big|\check{X}_{t_n}^{i,\tilde{\ubsbeta}} - \check{X}_{t_n}^{i,\ubsbeta}\big|^2] \Delta t
            \\
            &\leq C  \Big(\mathbb{E}\max_{n=0,\dots,N_T}\big|X_{t_n}^{i,\tilde{\ubsbeta}} - X_{t_n}^{i,\ubsbeta}\big|^2+\Gamma^2 +\dfrac{1}{R}\Big)
            \\
            &\leq C\left(\Gamma^2 +\dfrac{1}{R}\right),
        \end{aligned}
    \end{equation}
    where $C$ may depend on the number of time steps $N_T$. 
    
In this way, the expression in \eqref{eq:cost_cont_YT_bound_new-discrete} can be bounded as follows
\begin{equation}
    \begin{aligned}
        \mathbb{E} \Big[ |\check{Y}_T^{i,\tilde{\ubsbeta}} - \check{Y}_T^{i,\ubsbeta}|^2 \Big]
        \leq C\Bigg( \mathbb{E}\max_{n=0,\dots,N_T}\big|\check{X}_{t_n}^{i,\tilde{\ubsbeta}} - \check{X}_{t_n}^{i,\ubsbeta}\big|^2 + \left( \Gamma^2+\dfrac{1}{R} \right)+  \mathbb{E}\big|\tilde{y}_0(X^{i}_0)- y_0(X^{i}_0)\big|^2 \Bigg).
    \end{aligned}
\end{equation}
By using the bound in \eqref{eq:bound-supp-diff-X-GammaR2-discrete}, we conclude that
\begin{equation}
\begin{aligned}
    \mathbb{E} \Big[ |\check{Y}_T^{i,\tilde{\ubsbeta}} - \check{Y}_T^{i,\bsbeta}|^2 \Big] &\leq C\Bigg( \mathbb{E}\big|\tilde{y}_0(X^{i}_0)- {y}_0(X^{i}_0)\big|^2+\left(\Gamma^2+\dfrac{1}{R}\right) \Bigg)\\
    &\leq C\left(\Gamma^2+\dfrac{1}{R}\right),
\end{aligned}
\end{equation}
where the last inequality follows from \eqref{eq:bound-diff-betas-GammaR2-discrete} at $t=0$.
\vskip10pt
The continuity of $\ubsbeta \mapsto \check{X}^{i,\ubsbeta}_T$ with respect to the control $\ubsbeta$ follows the same ideas given in the {\bf Step~2} of this proof and for the sake of space, it is omitted. 

The Lipschitz continuity of 
    $\ubsbeta \mapsto \check\mu^{N,\ubsbeta}_{t_{n}}$ steps from~\eqref{eq:proof-discrete-approx-measure-continuity-beta} and the continuity of $\check{X}^{i,\ubsbeta}_{t_n}$ for all $i=1,\dots,N$.

\noindent \textbf{Step 5:} By using the Lipschitz continuity of $ { \mbsP}$ (Assumption~\ref{hyp2:Pbar}) and $g$ (Assumption~\ref{hyp2:g-lip}) in its arguments and also by using the results from {\bf Step 4}, we conclude that 
\begin{equation*}
    |\overline{\mbsP}(\tilde{\ubsbeta}) - \overline{\mbsP}({\ubsbeta})| \leq C\Big(\Gamma^2+\frac{1}{R}\Big)^{1/2}.    
\end{equation*}
\endproof

\end{document}